\documentclass[a4paper,11pt]{amsart}

\usepackage{amsmath,amsthm,amscd,amssymb,amsfonts, amsbsy}
\usepackage{latexsym}
\usepackage{txfonts}
\usepackage[colorlinks,citecolor=red,pagebackref,hypertexnames=false]{hyperref}
\numberwithin{equation}{section}\usepackage{graphicx}
\usepackage{epstopdf}
\usepackage{geometry}                % See geometry.pdf to learn the layout options. There are lots.
\numberwithin{equation}{section}

\theoremstyle{plain}
\newtheorem{theorem}[equation]{Theorem}
\newtheorem{lemma}[equation]{Lemma}
\newtheorem{corollary}[equation]{Corollary}
\newtheorem{proposition}[equation]{Proposition}

\theoremstyle{definition}
\newtheorem{definition}[equation]{Definition}

\theoremstyle{remark}
\newtheorem{remark}[equation]{Remark}

\newcommand{\fiint}{\operatornamewithlimits{\fint\!\!\!\!\fint}}

\newcommand{\supp}{\operatorname{supp}}

\newcommand{\vt}{\vartheta}

\newcommand{\RR}{\mathbb{R}}
\newcommand{\NN}{\mathbb{N}}
\newcommand{\ZZ}{\mathbb{Z}}

\newcommand{\rn}{\mathbb{R}^n}

\newcommand{\A}{\mathcal{A}}
\newcommand{\Q}{\mathcal{Q}}

\newcommand{\eps}{\varepsilon}

\newcommand{\po}{{\partial\Omega}}

\newcommand{\bp}{\noindent {\it Proof}.\,\,}
\newcommand{\ep}{\hfill$\Box$ \vskip 0.08in}

\newcommand{\PP}{{\mathbb{P}}}

\newcommand{\Lmn}{{\mathcal L}^{m,n}}
\newcommand{\Lmno}{{\mathcal L}^{m,n}_o}

\def\ring{\mathaccent"0017 }

\begin{document}

\author{Svitlana Mayboroda and  Vladimir Maz'ya}

\title[Polyharmonic capacity and Wiener test]{Polyharmonic capacity and Wiener test of higher order}

%Boundedness and continuity at the boundary of the gradient
%for the solutions of the biharmonic equation

%MathSubj 		31B30, 31A30, 35G15, 35J40, 31A15

%\date{ }

\begin{abstract} 

In the present paper we establish the Wiener test for boundary regularity of the solutions to the polyharmonic operator. We introduce a new notion of polyharmonic capacity and demonstrate necessary and sufficient conditions on the capacity of the domain responsible for the regularity of a polyharmonic function near a boundary point. 

In the case of the Laplacian the test for regularity of a boundary point is the celebrated Wiener criterion of 1924. It was extended to the biharmonic case in dimension three by [Mayboroda, Maz'ya, Invent. Math. 2009]. As a preliminary stage of this work, in [Mayboroda, Maz'ya, Invent. Math. 2013] we demonstrated boundedness of the appropriate derivatives of solutions to the polyharmonic problem in arbitrary domains, accompanied by sharp estimates on the Green function. The present work pioneers a new version of capacity and establishes the Wiener test in the full generality of the polyharmonic equation  of arbitrary order. 

\end{abstract}

\maketitle

\tableofcontents

\section{Introduction}
\setcounter{equation}{0} 

The maximum principle for harmonic functions guarantees that a solution to the Laplace's equation with bounded Dirichlet data is always bounded, and further regularity theory assures that it is infinitely differentiable at any interior point. These results hold on arbitrary bounded open sets and require no additional geometrical restrictions. However, continuity of the solutions at the boundary has proven to be a far more delicate problem. For instance, harmonic functions are always continuous at a vertex of a cone (Poincar\'e \cite{P}, Zaremba \cite{Z}), while due to the Lebesgue's counterexample \cite{L} in a complement of a sufficiently thin cusp this property may actually fail.  

In 1924 Wiener introduced the harmonic (Wiener) capacity and established his famous criterion for regularity of a boundary point \cite{Wiener}. Roughly speaking, it states that a point $O\in\po$ is regular (i.e., every solution to the
Dirichlet problem for the Laplacian, with continuous data, is
continuous at $O$) if and only if the complement of the domain near
the point $O$, measured in terms of the Wiener capacity,
is sufficiently massive. More specifically, the harmonic capacity of
a compactum $K\subset\RR^n$ can be defined as
\begin{equation}\label{eq1.5-1}
{\rm cap}\,(K):=\inf\Bigl\{\|\nabla u\|_{L^2(\RR^n)}^2:\,\,u\in
C_0^\infty(\RR^n),\,\,u=1\mbox{ in a neighborhood of }K\Bigr\},
\end{equation}

\noindent where $n\geq 3$, and the regularity of the point $O$ is
equivalent to the condition
\begin{equation}\label{eq1.6-1}
\sum_{j=0}^\infty2^{j(n-2)}\, {\rm
cap}\,(\overline{B_{2^{-j}}}\setminus\Omega)=+\infty,
\end{equation}
where $B_{2^{-j}}$ is the ball of radius $2^{-j}$ centered at the origin. An appropriate version of this condition is also available in dimension $n=2$.
 
Wiener's result became one of the pillars of elliptic theory. It gave the first necessary and sufficient conditions characterizing the properties of the solutions in terms of the geometric features of the boundary.  The notion of capacity provided a non-linear analogue to the Lebesgue measure, suitable for the context of Sobolev spaces, and offered an adequate language to address many important problems in probability, potential theory, function spaces.  Over the years, Wiener's test found a large array of applications, and has been extended to a variety of second order differential equations, including, but not restricted to the general divergence form elliptic equations \cite{LSW}, degenerate elliptic equations \cite{FJK}, parabolic operators \cite{FGL}, \cite{EG}, Schr\"dinger operator \cite{DMM}, and various quasilinear and nonlinear contexts \cite{MZ},  \cite{TW}, \cite{La}, \cite{AHe}. One can consult, e.g., \cite{MWiener}, \cite{AWiener}, for an excellent survey of related results.

Despite all these successes, the higher order operators remained out of reach. Elliptic equations of order greater than two, particularly based on the powers of the Laplacian, are common in physics and in engineering design, with applications ranging from standard models of elasticity \cite{Meleshko} to cutting-edge research of Bose-Einstein condensation in graphene and similar materials \cite{Kamenev}. They naturally appear in many areas of mathematics too, including conformal geometry (Paneitz operator, $Q$-curvature \cite{AliceChang}, \cite{AliceChang2}), free boundary problems \cite{AdamsPotAnal}, non-linear elasticity \cite{Skrypnik},  \cite{CiarletElasticity123}, \cite{Antman2005}, and have enjoyed increasing attention in the past several decades.   Yet, none of the previously devised methods could handle the Wiener criterion or Wiener capacity in the higher order context, and moreover, even modest analogues of the maximum principle remained an open problem. Let us discuss those in more details,  for to address the Wiener criterion for continuity of the higher order derivatives, one has to establish their boundedness first. 

The study of the higher order PDEs on {\it smooth} domains went hand-in-hand with the second order theory and, in particular, in 1960 the weak maximum principle has been established (\cite{Agmon}, see also \cite{Mir48}, \cite{Mir58}). Roughly speaking, it
amounts to the estimate
\begin{equation}\label{eq1.1}
\sum_{k\leq m-1}\|\nabla^{k}u\|_{L^\infty(\overline\Omega)}\leq
C\sum_{k\leq m-1} \|\nabla^{k}u\|_{L^\infty(\po)},
\end{equation}
\noindent where $u$ is a solution of an elliptic differential
equation of order $2m$ with smooth coefficients such that $|\nabla^{m-1}u|$ is continuous
up to the boundary, $\nabla^{k}u=\{\partial^\alpha
u\}_{|\alpha|=k}$ is a vector of all partial derivatives of $u$
of order $k$ and we adopt the usual convention that the zeroth order derivative of $u$ is $u$ itself.  Here, $\Omega$ was of class $C^{2m}$.

More recently, with the breakthrough in understanding of boundary value problems on Lipschitz domains, \eqref{eq1.1} was proved for the $m$-Laplacian, $(-\Delta)^m$, on three-dimensional domains with Lipschitz boundary   (\cite{PVmax}, \cite{PVpoly}; see also \cite{DJK}, \cite{PVLp}, \cite{Shen1}, \cite{Shen2} for related work), and, by different methods, in three-dimensional domains diffeomorphic to a polyhedron (\cite{KMR}, \cite{MR1}). It was also established that in both cases  $\nabla^{m-1}u$ is continuous whenever $u$ is an $m$-harmonic function with nice data. However, it turned out that {\it for every elliptic operator
of order greater than two the maximum principle can be violated}, in a four-dimensional cone  (\cite{MR}, see also \cite{MazNP83},  \cite{PVLp}, \cite{PVpoly}). In particular, in dimensions $n\geq 4$ there are solutions to the polyharmonic equation with unbounded $(m-1)$-st derivatives (cf. \eqref{eq1.1}).  

%Since those discoveries, increasing interest to higher order PDEs brought many delicate and powerful results (see, e.g., \cite{DKV}, \cite{Ver90}, \cite{PV91}, \cite{PVLp}, \cite{PVmax}, \cite{KenigBook}, \cite{PVpoly},  \cite{PVdilation}, \cite{Ver96}, \cite{DKPV97}, \cite{AdPi}, \cite{Ver05}, \cite{She06A},  \cite{She06B}, \cite{She06C}, \cite{She07A}, \cite{Mit10}, \cite{MazMS10}, \cite{Ver10}, \cite{KS11A}, \cite{MitMW11}), in the half-space, in smooth and Lipschitz domains. 

These results and counterexamples raise a number of fundamental questions: whether \eqref{eq1.1} could be extended to arbitrary domains in dimension 3, on par with the maximum principle for the Laplacian; whether in higher dimensions one can establish results of similar magnitude, possibly for lower order derivatives; and finally, if any of these answers is positive, whether one could aspire to get the capacitory conditions governing continuity of the appropriate derivatives.

{\it The series of papers \cite{MayMaz2}, \cite{MayMazGr}, \cite{MayMazInvent2}, culminating at the present manuscript, achieves a complete description of the boundary regularity of polyharmonic functions in arbitrary domains, provides sharp dimensional restrictions on boundedness of the derivatives of the solution, and establishes geometric conditions on the domain necessary and sufficient for their continuity, an analogue of the Wiener test.} The methods rely on intricate weighted integral inequalities and a new notion of higher-order capacities, see a more detailed discussion below teh statements of the main results.

In \cite{MayMazInvent2} we have established the exact order of smoothness for polyharmonic functions on domains with no geometrical restrictions. The principal result reads as follows.
\begin{theorem}\label{t1.1-I2}
Let $\Omega$ be a bounded domain in $\RR^n$, $2\leq n \leq 2m+1$,
and
\begin{equation}\label{eq1.2-I2}
(-\Delta)^m u=f \,\,{\mbox{in}}\,\,\Omega, \quad f\in
C_0^{\infty}(\Omega),\quad u\in \ring
W^{m,2}(\Omega).
\end{equation}

\noindent Then  the
solution to the boundary value problem \eqref{eq1.2-I2}
satisfies
\begin{equation}\label{eq1.3-I2}
\nabla^{m-n/2+1/2} u\in L^\infty(\Omega)\,\,\mbox{when $n$ is odd\quad and \quad} \nabla^{m-n/2} u\in L^\infty(\Omega)\,\,\mbox{when $n$ is even}.
\end{equation}

\noindent In particular, 
\begin{equation}\label{eq1.4-I2}
\nabla^{m-1} u\in L^\infty(\Omega)\,\,\mbox{when $n=2,3$}.
\end{equation}

\end{theorem}

\noindent Here the space $\ring W^{m,2}(\Omega)$, is, as usually, a
completion of $C_0^\infty(\Omega)$ in the norm given by $\|u\|_{\ring
W^{m,2}(\Omega)}=\|\nabla^m u\|_{L^2(\Omega)}$. We note that $\ring W^{m,2}(\Omega)$ embeds into $C^k(\Omega)$ only when $k$ is strictly smaller than $m-\frac n2$, $n<2m$. Thus, whether the dimension is even or odd, Theorem~\ref{t1.1-I2} gains one derivative over the outcome of Sobolev embedding.

The results of Theorem~\ref{t1.1-I2} are sharp, in the sense that the solutions  do not exhibit higher smoothness than warranted by \eqref{eq1.3-I2}--\eqref{eq1.4-I2} in general domains. To be more precise, when the dimension $n\in [3,2m+1]\cap\NN$ is odd, one can find a solution in a punctured ball  
for which  the derivatives of the order $m-\frac n2+\frac 32$ fail to be bounded, and moreover, $\nabla^{m-\frac n2+\frac 12} u$
is not 
continuous at the origin \cite{MayMazInvent2}. On the other hand, when $n$ is even, the results in \cite[Section~10.4]{KMR} demonstrate  that in an exterior of a ray there is an $m$-harmonic function behaving as $|x|^{m-\frac n2+\frac 12}$. Thus, upon a suitable truncation, one obtains a solution to \eqref{eq1.2-I2} in $B_1\setminus\{x_1=0, ..., x_{n-1}=0, 0\leq x_n<1\}$, whose derivatives of order $m-\frac n2+1$ are not bounded, confirming sharpness of \eqref{eq1.3-I2} in even dimensions too, but for now at the level of full (rather than fractional) derivatives. We shall return to this topic below. 

At this point Theorem~\ref{t1.1-I2} finally sets the stage for a discussion of the {\it Wiener test} for continuity of the corresponding derivatives of the solution, which brings us to the main results of the present paper.  

Assume that $m\in\NN$ and $n\in [2,2m+1]\cap \NN$. Let us denote by $Z$ the following set of indices: 
\begin{eqnarray}
\label{z1} && Z=\{0,1,...,m-n/2+1/2\}, \quad \mbox{if $n$ is odd}, \\[4pt]
\label{z2} && Z=\{-n/2+2, -n/2+4,..., m-n/2-2, m-n/2\}\cap(\NN\cup\{0\}), \,\, \mbox{if $n$ is even, $m$ is even}, \\[4pt]
\label{z3} && Z=\{-n/2+1,-n/2+3,..., m-n/2-2, m-n/2\}\cap(\NN\cup\{0\}), \,\, \mbox{if $n$ is even, $m$ is odd}.
\end{eqnarray}
Now let $\Pi$ be the space of linear combinations of spherical harmonics
\begin{equation}\label{cap6.1}
P(x)=\sum_{p\in Z} \sum_{l=-p}^{p} b_{pl}
Y_l^p(x/|x|), \qquad b_{pl}\in\RR, \quad x\in\RR^n\setminus \{O\},
\end{equation}

\noindent with the norm 
\begin{equation}\label{cap6.2}
\|P\|_{\Pi}:=\left(\sum_{p\in Z} \sum_{l=-p}^{p} b_{pl}^2\right)^{\frac 12} \quad \mbox{and}\quad \Pi_1:=\{P\in\Pi:\,\|P\|_{\Pi}=1\}.
\end{equation}

Then, given $P\in\Pi_1$, an open set $D$ in $\RR^n$ such that $O\in
\RR^n\setminus D$, and a compactum $K$ in $D$, we define
\begin{equation}\label{cap8}
{\rm Cap}_P\,(K,D):=\inf\Bigg\{\int_D |\nabla^m u(x)|^2\,dx:\,\,u\in \ring W^{m,2}(D),\,\,u=P\mbox{ in a
neighborhood of }K\Bigg\},
\end{equation}

\noindent with 
\begin{equation}\label{cap9}
{\rm Cap}\,(K,D):=\inf_{P\in \Pi_1} {\rm Cap}_P\,(K,D).
\end{equation}
In the context of the Wiener test, we will be working extensively with the capacity of the complement of a domain $\Omega\subset \RR^n$ in the balls $B_{2^{-j}}$, $j\in\NN$, and even more so, in dyadic annuli, $C_{2^{-j}, 2^{-j+2}}$, $j\in\NN$, where $C_{s,as}:=\{x\in\RR^n:\,s<|x|<as\}$, $s,a>0$. Following the custom, it will be convenient to abbreviate dropping the reference to the ``ambient" set
\begin{equation}\label{drop_amb}{{\rm
Cap}_P\,(\overline{C_{2^{-j},2^{-j+2}}}\setminus\Omega)}:= {{\rm
Cap}_P\,(\overline{C_{2^{-j},2^{-j+2}}}\setminus\Omega, C_{2^{-j-2},2^{-j+4}})},\quad j\in\NN,
\end{equation}
and similarly for ${\rm Cap}$. In fact, it will be proved below that there are several equivalent definitions of capacity, in particular, 
for any $n\in [2, 2m+1]$  and for any $s>0$, $a>0$, $K\subset \overline{C_{s,as}}$, we have  
\begin{multline}\label{eq5.12-2-intro}
\inf\Bigg\{\sum_{k=0}^m\int_{\RR^n} \frac{|\nabla^ku(x)|^2}{|x|^{2m-2k}}\,dx:\,\,u\in \ring W^{m,2}(\RR^n\setminus\{O\}),\,\,u=P\mbox{ in a
neighborhood of }K\Bigg\}\\[4pt]\approx {\rm Cap}_P(K,
C_{s/2,2as}).
\end{multline}
In the case when the dimension is odd, also
$${{\rm
Cap}_P\,(\overline{C_{s,as}}\setminus\Omega, C_{s/2, 2as}})\approx {{\rm
Cap}_P\,(\overline{C_{s,as}}\setminus\Omega, \RR^n\setminus \{O\})}.$$
Thus, either of the above can be used in \eqref{drop_amb} conditions as convenient. 
%Using this observation and monotonicity properties of capacity, one can also work with $a$-adic rather than dyadic annuli and various other choices of the test sets in the capacitory conditions in Theorem~\ref{to1.2}.

Let $\Omega$ be a domain in $\RR^n$, $n\geq 2$. The point $Q\in\po$ is {\it $k$-regular} with
respect to the domain $\Omega$ and the operator $(-\Delta)^m$, $m\in\NN$, if the solution to the boundary problem
\begin{equation}\label{eq6.1}
(-\Delta)^m u=f \,\,{\mbox{in}}\,\,\Omega, \quad f\in
C_0^{\infty}(\Omega),\quad u\in \ring W^{m,2}(\Omega),
\end{equation}

\noindent satisfies the condition
\begin{equation}\label{eq6.2}
\nabla^k u(x)\to 0\mbox{ as } x\to Q,\,x\in\Omega,
\end{equation}

\noindent that is, all partial derivatives of $u$ of order $k$ are continuous.
Otherwise, we say that $Q\in\po$ is $k$-irregular.

The main result of this paper is as follows. 

\begin{theorem}\label{to1.2} Let $\Omega$ be an arbitrary open set in $\RR^n$, $m\in\NN$, $2\leq n \leq 2m+1$. Let $\lambda$ be given by 
\begin{equation}\label{eqo7.5}
\lambda=\left\{\begin{array}{l} m-n/2+1/2 \quad \,\,\mbox{when $n$ is odd},\\[8pt] m-n/2\qquad\qquad\mbox{when $n$ is even}.\end{array}
\right.
\end{equation}
If 
\begin{equation}\label{eqo1.9}
\sum_{j=0}^\infty 2^{-j(2m-n)} \,\inf_{P\in\Pi_1}{{\rm
Cap}_P\,(\overline{C_{2^{-j},2^{-j+2}}}\setminus\Omega)} =+\infty, \quad \mbox{when $n$ is odd},
\end{equation}
and 
\begin{equation}\label{eqo1.9-even}
\sum_{j=0}^\infty j\,2^{-j(2m-n)} \,\inf_{P\in\Pi_1}{{\rm
Cap}_P\,(\overline{C_{2^{-j},2^{-j+2}}}\setminus\Omega)}=+\infty, \quad \mbox{when $n$ is even},
\end{equation}

\noindent then the point $O$ is $\lambda$-regular with
respect to the domain $\Omega$ and the operator $(-\Delta)^m$.

Conversely, if the point $O\in\po$ is $\lambda$-regular with
respect to the domain $\Omega$ and the operator $(-\Delta)^m$ then 
\begin{equation}\label{eqo1.10}
\inf_{P\in\Pi_1} \sum_{j=0}^\infty 2^{-j(2m-n)} \,{{\rm
Cap}_P\,(\overline{C_{2^{-j},2^{-j+2}}}\setminus\Omega)}=+\infty, \quad \mbox{when $n$ is odd},
\end{equation}
and 
\begin{equation}\label{eqo1.10-even}
\inf_{P\in\,\Pi_1} \sum_{j=0}^\infty j\,2^{-j(2m-n)} \,{{\rm
Cap}_P\,(\overline{C_{2^{-j},2^{-j+2}}}\setminus\Omega)}=+\infty, \quad \mbox{when $n$ is even}. \end{equation}

\noindent Here, as before,  $C_{2^{-j},2^{-j+2}}$ is the annulus $\{x\in\RR^n:\,2^{-j}<|x|<2^{-j+2}\}$, $j\in\NN\cup \{0\}$.
\end{theorem}

Let us now discuss the results of Theorem~\ref{to1.2} in more details. This is the first treatment of the continuity of derivatives of an elliptic equation of order $m>2$ at the boundary, and the first time the capacity \eqref{cap8} appears in literature. When applied to the case $m=1$, $n=3$, it yields the classical Wiener criterion for continuity of a harmonic function (cf. \eqref{eq1.5-1}--\eqref{eq1.6-1}). Furthermore, continuity of the solution itself (rather than  its derivatives) has been previously treated for the polyharmonic equation, and for $(-\Delta)^m$ the resulting criterion also follows from Theorem~\ref{to1.2}, in particular, when $m=2n$, the new notion of capacity \eqref{z1}--\eqref{cap6.2}. coincides with the potential-theoretical Bessel capacity used in \cite{M2}. In the case $\lambda=0$, covering both of the above, necessary and sufficient condition in Theorem~\ref{to1.2} are trivially the same, as $P\equiv 1$ when $n=2m$ in even dimensions and $n=2m+1$ in odd ones. 
For lower dimensions $n$ the discrepancy is not artificial, for, e.g., \eqref{eqo1.9} may fail to be necessary as was shown in \cite{MayMaz2}. 
Finally, as we pointed out already, the bilaplacian in dimension three was our first result pioneering this line of work, and was addressed in \cite{MayMaz2}. 

It is not difficult to verify that we also recover aforementioned bounds in Lipschitz and in smooth domains, as the capacity of a cone and hence, capacity of an intersection with a complement of a Lipschitz domains, assures divergence of the series in \eqref{eqo1.9}--\eqref{eqo1.9-even}. On the other hand, given Theorem~\ref{to1.2} and following traditional in this context considerations (choosing sufficiently small balls in the consecutive annuli to constitute a complement of the domain), we can build a set with a convergent capacitory integral and, respectively, an irregular solution with discontinuous derivatives of order $\lambda$ at the point $O$. Note that this yields further sharpness of the results of Theorem~\ref{t1.1-I2}. In particular, in even dimensions, it is a stronger counterexample than that of a continuum  discussed above (not only $m-n/2+1$ derivatives are not bounded, but $m-n/2$ derivatives might be discontinuous). 
   
Before the proof, let us say a few words about our methods and highlight new challenges, particularly  in comparison with the biharmonic case. First of all, odd and even dimensions prove to yield very different problems. Our approach is rooted in weighted integral inequalities whose nature heavily depends on the parity of the dimension. In particular, when $n$ is even, the situation is significantly influenced by an additional logarithmic term. Moreover, the case of even dimensions splits further depending on the parity of $m-n/2$ and the parity of $m$. The underlying effects can already be glimpsed from the definition of $\Pi$, the space of ``boundary data" of the new capacity (cf. \eqref{z1}--\eqref{z3}). The treatment of the biharmonic problem in \cite{MayMaz2} is restricted to dimension three and somewhat resonates with the ideas used in odd dimensions here. However, even for $n$ odd, we need to use  heavy machinery of \cite{MayMazInvent2} in order to assure positivity of various terms in resulting weighted inequalities and to develop tools to control the others: expressions which could be explicitly computed for $\Delta^2$ in dimension 3 now lead to severe technical obstructions.

One of the most difficult aspects of proof of Theorem~\ref{to1.2} is finding a correct notion of polyharmonic capacity and understanding its key properties. A peculiar choice of linear combinations of spherical harmonics (see \eqref{z1}--\eqref{z3} and \eqref{cap6.1}) is crucial at several stages of the argument, specific to the problem at hand, and no alterations would lead to reasonable necessary and sufficient conditions. At the same time, the new capacity and the notion of higher-order regularity sometimes exhibit surprising properties, such as for instance sensitivity to the affine changes of coordinates \cite{MayMaz2}, or the fact that in sharp contrast with the second order case \cite{LSW}, one does not expect same geometric conditions to be responsible for regularity of solutions to all higher order elliptic equations. For instance, the solution to the biharmonic equation $(-\Delta)^2u=0$ is continuous at a vertex of a cone in any dimension, while this property fails for $[(-\Delta^2)+a (\partial/\partial x_n)^4]u=0$ in dimensions $n\geq 8$ for all $a>5+2\sqrt 5$ \cite{MWiener}. This underlines the delicacy of the analysis of new capacitory conditions: recall that in the second order case the regularity of the solutions to all divergence form elliptic equations is governed by the same capacity \eqref{eq1.5-1}. Other  features will be discussed in the body of the paper.

\section{Regularity of solutions to the polyharmonic equation}\label{srecall}
\setcounter{equation}{0}

In the present section we set the notation and recall the main results of \cite{MayMazInvent2} which will be extensively used in the present paper. 

Let us start with a list of notation and conventions used throughout the paper.

For any domain $\Omega\subset \RR^n$ a function  $u\in
C_0^\infty(\Omega)$ can be extended by zero to $\RR^n$ and we will
write $u\in C_0^\infty(\RR^n)$ whenever convenient. Similarly, the
functions in $\ring W^{m,2}(\Omega)$, $m\in\NN$, will be extended by zero and
treated as functions on $\RR^n$ or other open sets containing $\Omega$ without further comments.

The symbols $B_r(Q)$ and $S_r(Q)$ denote,
respectively, the ball and the sphere with radius $r$ centered at
$Q$ and $C_{r,R}(Q)=B_R(Q)\setminus\overline{B_r(Q)}$. When the
center is at the origin, we write $B_r$ in place of $B_r(O)$, and
similarly $S_r:=S_r(O)$ and $C_{r,R}:=C_{r,R}(O)$.

Let $(r,\omega)$ be spherical coordinates in $\RR^n$, $n\geq 2$, i.e. $r=|x|\in
(0,\infty)$ and $\omega=x/|x|$ is a point of the unit sphere $S^{n-1}$.
In fact, it will be more convenient to use $e^{-t}$, $t\in\RR$, in place of $r$, so that $t=\log r^{-1}=\log |x|^{-1}$. Then by  
$\varkappa$ we denote the mapping
\begin{equation}\label{eqo2.2}
\RR^n \ni
x\,\stackrel{\varkappa}{\longrightarrow} \,(t,\omega)\in\RR\times
S^{n-1},\quad n\geq 2.
\end{equation}

\noindent The symbols $\delta$ and $\nabla_\omega$ refer,
respectively, to the Laplace-Beltrami operator and the gradient on
$S^{n-1}$.

Finally, by $C$, $c$, $C_i$ and $c_i$, $i\in\NN$, we generally denote some
constants, possibly depending on the order of operator $m$ and the dimension $n$ but not on any other variables and not on the domain, unless explicitly stated otherwise. Their exact values are of no importance and can change from line to line. Also, we write
$A\approx B$, if $C^{-1}\,A\leq B\leq C\,A$ for some $C>0$.

\begin{theorem}\label{to2.1} Assume that $m\in\NN$ and $n\in [3,2m+1]\cap \NN$
is odd. Let $\Omega$ be a bounded domain in $\RR^n$, $O\in
\RR^n\setminus\Omega$, $u\in C_0^\infty(\Omega)$ and
$v=e^{\left(m-\frac n2+\frac 12\right)t}(u\circ \varkappa^{-1})$.
Then
\begin{eqnarray}\label{eqo2.3}
\hskip -1cm &&\int_{\RR^n}(-\Delta)^m u(x)\,u(x)|x|^{-1}\,dx\nonumber\\[4pt]
\hskip -1cm &&\quad \geq C\sum_{\stackrel{k\geq 1, \,i\geq 0}{i+k\leq m}} \int_{\RR}\int_{S^{n-1}}
\left(\partial_t^k \nabla_\omega^i v\right)^2\,d\omega dt+C
\int_{\RR}\int_{S^{n-1}} v\,\prod_{p=-\frac n2+\frac 32}^{m-\frac
n2+\frac 12}\left(-\delta-p\,(p+n-2)\right) v\,d\omega dt,
\end{eqnarray}

\noindent where $C>0$ is some constant depending on $m$ and $n$
only.
\end{theorem}

\bp The inequality 
\begin{eqnarray}\label{eqo2.3.1}
\hskip -1cm &&\int_{\RR^n}(-\Delta)^m u(x)\,u(x)|x|^{-1}\,dx\nonumber\\[4pt]
\hskip -1cm &&\quad \geq C\sum_{k=1}^m \int_{\RR}\int_{S^{n-1}}
\left(\partial_t^k v\right)^2\,d\omega dt+C
\int_{\RR}\int_{S^{n-1}} v\,\prod_{p=-\frac n2+\frac 32}^{m-\frac
n2+\frac 12}\left(-\delta-p\,(p+n-2)\right) v\,d\omega dt,
\end{eqnarray}
is the statement of Theorem~2.1 in \cite{MayMazInvent2}. It only remains to show that a part of the sum corresponding to $i\neq 0$ in \eqref{eqo2.3} is controlled by the remaining terms. 

To this end, let us recall the notation from the proof of Theorem~2.1 in \cite{MayMazInvent2}. In the system of coordinates $(t,\omega)$
the polyharmonic operator can be written as
\begin{equation}\label{eqo2.4}
(-\Delta)^m=(-1)^me^{2mt}\prod_{j=0}^{m-1}\Bigl((-\partial_t-2j)(-\partial_t-2j+n-2)+\delta\Bigr).
\end{equation}

\noindent Then
\begin{equation}\label{eqo2.5}
\int_{\RR^n}(-\Delta)^m
u(x)\,u(x)|x|^{-1}\,dx=\int_{\RR}\int_{S^{n-1}} {\mathcal
L}^{m,n}(\partial_t,\delta)v(t,\omega)\,v(t,\omega)\,d\omega dt,
\end{equation}

\noindent with
\begin{equation}\label{eqo2.6}
{\mathcal L}^{m,n}(\partial_t,\delta)=(-1)^m\prod_{j=0}^{m-1}
\Bigg(\Bigl(-\partial_t+m-\frac n2+\frac
12-2j\Bigr)\Bigl(-\partial_t+m+\frac n2-\frac
32-2j\Bigr)+\delta\Bigg).
\end{equation}

\noindent Denote by $v_{pl}$ the coefficients of the expansion of
$v$ into spherical harmonics:
\begin{equation}\label{eqo2.7}
v(t,\omega)=\sum_{p=0}^\infty \sum_{l=-p}^{p} v_{pl}(t)
Y_l^p(\omega), \qquad t\in\RR, \,\,\omega\in S^{n-1}.
\end{equation}

\noindent Then we can write the expression on the right-hand side
of (\ref{eqo2.5}) as
\begin{equation}\label{eqo2.8}
\sum_{p=0}^\infty \sum_{l=-p}^{p} \int_{\RR} {\mathcal
L}^{m,n}(\partial_t,-p\,(p+n-2)) v_{pl}(t) \,v_{pl}(t)\,dt.
\end{equation}

\noindent We remark for future reference that, denoting by $\widehat v$ the Fourier transform
of $v$, i.e.
\begin{equation}\label{eqo2.9}
\widehat v(\gamma)=\frac{1}{\sqrt{2\pi}} \int_\RR e^{-i\gamma\,
t}v(t)\,dt,\qquad\gamma\in\RR.
\end{equation}

\noindent we have by the Plancherel's identity that (\ref{eqo2.8}) is equal to
\begin{eqnarray}\label{eqo2.10}\nonumber
&&\sum_{p=0}^\infty \sum_{l=-p}^{p} \int_{\RR} {\mathcal
L}^{m,n}(i\gamma,-p\,(p+n-2))
\left|\widehat{v_{pl}}(\gamma)\right|^2\,d\gamma\\[4pt]
&&\qquad\qquad = \sum_{p=0}^\infty \sum_{l=-p}^{p} \int_{\RR} \Re
e\, {\mathcal L}^{m,n}(i\gamma,-p\,(p+n-2))
\left|\widehat{v_{pl}}(\gamma)\right|^2\,d\gamma.
\end{eqnarray}

\noindent See \cite{MayMazInvent2} for details. 

Now, for $1\leq k$, $0\leq i$, $i+k\leq m$, we have

\begin{eqnarray}\label{eqoo4.35}
&&\int_{\RR}\int_{S^{n-1}} (\partial_t^k\nabla_\omega^i
v(t,\omega))^2 \,d\omega dt =\sum_{p=0}^\infty \sum_{l=-p}^{p}
\Bigl(p\,(p+n-2)\Bigr)^i \int_{\RR}(\partial_t^k v_{pl}(t))^2 \,dt.
\end{eqnarray}

\noindent We break down the sum above into two parts,
corresponding to the cases $p\leq m-\frac n2+\frac 12$ and $p\geq m-\frac n2+\frac 32$, respectively. In the first case,
\begin{eqnarray}\label{eqoo4.36}
&&\sum_{p=0}^{m-\frac n2+\frac 12} \sum_{l=-p}^{p}
\Bigl(p\,(p+n-2)\Bigr)^i \int_{\RR}(\partial_t^k v_{pl}(t))^2
\,dt\leq C_{m,n} \sum_{p=0}^{m-\frac n2+\frac 12} \sum_{l=-p}^{p}
\int_{\RR}(\partial_t^k v_{pl}(t))^2 \,dt\nonumber\\[4pt]
&&\quad\leq C_{m,n}
\int_{\RR}\int_{S^{n-1}} (\partial_t^k v(t,\omega))^2 \,d\omega dt, 
\end{eqnarray}

\noindent which is bounded by the right-hand side of \eqref{eqo2.3.1}. As for $p\geq m-\frac n2+\frac 32$, since $i+k\leq m$ then by Young's inequality
\begin{eqnarray}\label{eqoo4.39}
&& \Bigl(p\,(p+n-2)\Bigr)^i \int_{\RR}(\partial_t^k v_{pl}(t))^2
\,dt= \int_{\RR}\Bigl(p\,(p+n-2)\Bigr)^i \gamma^{2k}
|\widehat{v_{pl}}(\gamma)|^2 \,d\gamma\nonumber\\[4pt] &&\qquad \leq
\int_{\RR}\Bigl(\Bigl(p\,(p+n-2)\Bigr)^m +\gamma^{2m}\Bigr)
|\widehat{v_{pl}}(\gamma)|^2 \,d\gamma\nonumber\\[4pt] &&\qquad \leq\Bigl(p\,(p+n-2)\Bigr)^m
\int_{\RR}(v_{pl}(t))^2 \,dt +\int_{\RR}(\partial_t^m v_{pl}(t))^2
\,dt.
\end{eqnarray}

\noindent However, 
\begin{equation}\label{eqoo4.37}
\Bigl(p\,(p+n-2)\Bigr)^m\leq C \prod_{s=-\frac n2+\frac
32}^{m-\frac n2+\frac
12}\Bigl(p\,(p+n-2)-s\,(s+n-2)\Bigr),\quad\mbox{ for every $p\geq
m-\frac n2+\frac 32$},
\end{equation}

\noindent where $C>0$ depends on $m$ and $n$ only. This follows from the fact that  one can choose $C$ such that 
\begin{equation}\label{eqoo4.38}
\frac 1C\leq \left(1-\frac{\left(m-\frac n2+\frac 12\right)\left(m+\frac n2-\frac 32\right)}{\left(m-\frac n2+\frac 32\right)\left(m+\frac n2-\frac 12\right)}\right)^m\leq \prod_{s=-\frac n2+\frac 32}^{m-\frac n2+\frac
12}\left(1-\frac{s\,(s+n-2)}{p\,(p+n-2)}\right),
\end{equation}

\noindent for every $p\geq
m-\frac n2+\frac 32$.

Therefore, by \eqref{eqoo4.37}
\begin{eqnarray}\label{eqoo4.40}
&&\sum_{p=m-\frac n2+\frac 32}^{\infty} \sum_{l=-p}^{p} \Bigl(p\,(p+n-2)\Bigr)^i
\int_{\RR}(\partial_t^k v_{pl}(t))^2 \,dt \nonumber\\[4pt] &&\qquad \leq\sum_{p=m-\frac n2+\frac 32}^{\infty} \sum_{l=-p}^{p}
\,\,\prod_{s=-\frac n2+\frac 32}^{m-\frac n2+\frac
12}\Bigl(p\,(p+n-2)-s\,(s+n-2)\Bigr)
\int_{\RR}(v_{pl}(t))^2 \,dt \nonumber\\[4pt]&&\qquad\qquad +\sum_{p=m-\frac n2+\frac 32}^{\infty}
\sum_{l=-p}^{p}\int_{\RR}(\partial_t^m v_{pl}(t))^2 \,dt,
\end{eqnarray}

\noindent  which is also bounded by the right-hand side of \eqref{eqo2.3.1}, now invoking both the first and the second term in \eqref{eqo2.3.1}.  This finishes the proof of \eqref{cap6.0.1}.
 \ep

\begin{lemma}\label{lo4.1} Assume that $m\in\NN$ and $n\in [3,2m+1]\cap \NN$
is odd. Consider the equation
\begin{equation}
\Lmn(-\partial_t,0)\,h
 =\delta,  \label{eqo4.1}
\end{equation}

\noindent where $\delta$ stands for the Dirac delta function. A
unique solution to {\rm (\ref{eqo4.1})} which is bounded and
vanishes at $+\infty$ has a form
\begin{equation}\label{eqo4.2}
h(t)=\left\{\begin{array}{l}
\sum_{j=1}^m \nu_j \,e^{-\alpha_j t},\qquad \qquad\qquad t>0,\\[4pt]
\sum_{j=1}^m \mu_j \,e^{\beta_j t},\qquad \qquad\qquad\,\,t<0.\\[4pt]
\end{array}
\right.
\end{equation}

\noindent Here $\alpha_j>0$, $j=1,2,...,m$, $\beta_j>0$ for
$j=2,...,m$ and $\beta_1=0$ are such that

\begin{equation}\label{eqo4.3}
\{-\alpha_j\}_{j=1}^{m}\bigcup
\{\beta_j\}_{j=1}^{m}=\Bigl\{-m+\frac n2-\frac
12+2j\Bigr\}_{j=0}^{m-1}\cup \Bigl\{-\frac n2 -\frac
12+m-2j\Bigr\}_{j=0}^{m-1},
\end{equation}

\noindent and with the notation
\begin{equation}\label{eqo4.4}
\vec\gamma=(-\alpha_1,...,-\alpha_m, \beta_1,...,\beta_m),\qquad
\vec\kappa=(\nu_1,...,\nu_m, -\mu_1,...,-\mu_m)
\end{equation}
\noindent the coefficients $\nu_j,\mu_j\in\RR$ satisfy
\begin{equation}\label{eqo4.5}
\kappa_i=(-1)^{m+1}\Bigl(\prod_{j\neq i}(\gamma_j-\gamma_i)\Bigr)^{-1}.
\end{equation}
\end{lemma}

\begin{theorem}\label{to4.2} Assume that $m\in\NN$ and  $n\in [3,2m+1]\cap \NN$
is odd. Let $\Omega$ be a bounded domain in $\RR^n$, $O\in
\RR^n\setminus\Omega$, $u\in C_0^\infty(\Omega)$ and
$v=e^{\left(m-\frac n2+\frac 12\right)t}(u\circ \varkappa^{-1})$. 
Then for
every $\xi\in\Omega$ and $\tau=\log |\xi|^{-1}$ 
\begin{multline}\label{eqo4.11}
\int_{S^{n-1}}v^2(\tau,\omega)\,d\omega \\[4pt] + \sum_{\stackrel{k\geq 1, \,i\geq 0}{i+k\leq m}} \int_{\RR}\int_{S^{n-1}}
\left(\partial_t^k \nabla_\omega^i v\right)^2\,d\omega dt+
\int_{\RR}\int_{S^{n-1}} v\,\prod_{p=-\frac n2+\frac 32}^{m-\frac
n2+\frac 12}\left(-\delta-p\,(p+n-2)\right) v\,d\omega dt \\[4pt] \leq  C
\int_{\RR^n}(-\Delta)^m u(x)\,u(x)g(\log |x|^{-1},\log |\xi|^{-1})\,dx,
\end{multline}

\noindent with
\begin{equation}\label{eqo4.12}
g(t,\tau)=e^{t}\left(C_1 h(t-\tau)+C_2\right),\qquad t,\tau\in\RR,
\end{equation}

\noindent and $h$ given by Lemma~\ref{lo4.1}. Here  $C,C_1,C_2$ are some positive constants depending on $m$ and $n$ only.
\end{theorem}

\bp This is the result of Theorem~4.2 from \cite{MayMazInvent2}, {\it loc.cit.} There, in the statement we only display the first term on the left-hand side of \eqref{eqo4.11}, but it is clear that one can add the terms appearing in the second line at the expense of possibly augmenting the constant $C_2$ in \eqref{eqo4.12}, owing to Theorem~\ref{to2.1}. \ep

\begin{remark}\label{ro1} We mention for the future record that the proof of Theorem~4.2 from \cite{MayMazInvent2}, {\it loc.cit.}, demonstrates an additional term in the lower estimate for the right-hand side of \eqref{eqo4.11}. Namely, in addition to \eqref{eqo4.11}, we have
\begin{multline}\label{eqo4.11.1}
\sum_{p=0}^{m-\frac n2 +\frac 12} \sum_{l=-p}^{p} \int_{\RR} 
v_{pl}^2(t) \left({\mathcal L}^{m,n}(-\partial_t,-p(p+n-2))-{\mathcal L}^{m,n}(-\partial_t,0)\right) h(t-\tau)\,  dt \\[4pt] \leq  C
\int_{\RR^n}(-\Delta)^m u(x)\,u(x)g(\log |x|^{-1},\log |\xi|^{-1})\,dx.
\end{multline}
It is shown in Proposition~3.1 from \cite{MayMazInvent2}, {\it loc.cit.}, that the weight of the integral on the left-hand side of \eqref{eqo4.11.1} is positive, that is, 
\begin{equation}\label{eqo4.26}
\left({\mathcal L}^{m,n}(-\partial_t,-p(p+n-2))-{\mathcal L}^{m,n}(-\partial_t,0)\right) h(t)\geq 0, \quad \mbox{ for all } \quad t\neq 0,\quad 0\leq p\leq  m-\frac n2+\frac 12.
\end{equation}
This fact allows to dispose of the corresponding terms and eventually arrive at \eqref{eqo4.11}, but we shall need them explicitly in the proof of necessity of the capacitary condition in this paper.   
\end{remark}

\begin{theorem}\label{to5.1} Assume that $m\in\NN$ and $n\in [2,2m]\cap \NN$
is even. Let $\Omega$ be a bounded domain in $\RR^n$, $O\in
\RR^n\setminus\Omega$, $u\in C_0^\infty(\Omega)$ and
$v=e^{\left(m-\frac n2\right)t}(u\circ \varkappa^{-1})$. Furthermore, let $R$ be a positive constant such that the support of $u$ is contained in $B_{2R}$, $C_R:=\log (4R)$, and let $\psi$ be a weight function such that either $\psi(t)=C_R+t$ for all $t\in\RR$ or $\psi(t)=1$ for all $t\in\RR$.
Then whenever  $m$ is even,
\begin{eqnarray}\label{eqo5.1}
\hskip -1cm &&\int_{\RR^n}(-\Delta)^m u(x)\,u(x)\,\psi(\log |x|^{-1})\,dx\geq C\sum_{k=1}^m \sum _{i=0}^{m-k}\int_{\RR}\int_{S^{n-1}}
\left(\partial_t^k \nabla_\omega^i v\right)^2\,\psi(t)\,d\omega dt\nonumber\\[4pt]
\hskip -1cm &&\quad +\,C
\int_{\RR}\int_{S^{n-1}} v\,\prod_{p}\left(-\delta-p\,(p+n-2)\right)^2 v\,\psi(t)\,d\omega dt,
\end{eqnarray}
\noindent where the product  is over $p=-n/2+2, -n/2+4,..., m-n/2-2, m-n/2$, that is, $p=-n/2+2j$ with $j=1,2,...,m/2$. If $m$ is odd,
\begin{eqnarray}\label{eqo5.2}
\hskip -1cm &&\int_{\RR^n}(-\Delta)^m u(x)\,u(x)\,\psi(\log |x|^{-1})\,dx\geq C\sum_{k=1}^m \sum _{i=0}^{m-k}\int_{\RR}\int_{S^{n-1}}
\left(\partial_t^k \nabla_\omega^i v\right)^2\,\psi(t)\,d\omega dt\nonumber\\[4pt]
\hskip -1cm &&\quad +\,C
\int_{\RR}\int_{S^{n-1}} v\,\prod_{p}\left(-\delta-p\,(p+n-2)\right)^2 \left(-\delta +(n/2-1)^2\right) v\,\psi(t)\,d\omega dt,
\end{eqnarray}

\noindent where the product is over $p=-n/2+3, -n/2+5,..., m-n/2-2, m-n/2$, that is, $p=-n/2+1+2j$ with $j=1,2,...,(m-1)/2$. In both cases $C>0$ is some constant depending on $m$ and $n$
only.
\end{theorem}

For future reference and to set the notation, let us record a few related details. 
According to the proof of Theorem~\ref{to5.1} in \cite{MayMazInvent2} (with the same choice of $\psi$ as in the statement of the Theorem), we can write 
\begin{eqnarray}\label{eqo5.3}
&&\int_{\RR^n}(-\Delta)^m
u(x)\,u(x)\,\psi(\log |x|^{-1})\,dx=\int_{\RR}\int_{S^{n-1}} {\mathcal
L}^{m,n}_o(\partial_t,\delta_\omega)v(t,\omega)\,v(t,\omega)\,\psi(t)\,d\omega dt\nonumber\\[4pt]
&&\qquad =\sum_{q=0}^{\infty}\sum_{l=-q}^q\int_{\RR}{\mathcal
L}^{m,n}_o(\partial_t,-q(q+n-2))v_{ql}(t)\,v_{ql}(t)\,\psi(t)\, dt,
\end{eqnarray}

\noindent where
\begin{equation}\label{eqo5.6}
{\mathcal L}^{m,n}_o(\partial_t,-q(q+n-2))=\prod_{j=0}^{m-1}
\left(-\partial_t^2+{\mathcal B}_j(q)^2\right),
\end{equation}

\noindent with
\begin{equation}\label{eqo5.7}
{\mathcal B}_j(q)^2=\Bigl(\sqrt{(n/2-1)^2+q(q+n-2)}+m-2j-1\Bigr)^2=\Bigl(q+n/2+m-2j-2\Bigr)^2,\quad q\in\NN\cup\{0\}
\end{equation}
\noindent satisfying the estimates
\begin{equation}\label{eqo5.8}
{\mathcal B}_j(q)^2\geq C \,\max\{1,q(q+n-2)\}, \mbox{when $q\in\NN\cup\{0\}$ is such that}\,\, q\neq 2j-m-n/2+2,
\end{equation}

\noindent for some $C>0$ depending on $m$ and $n$  only, and 
\begin{equation}\label{eqo5.9}
{\mathcal B}_j(q)=0\quad \mbox{if}\quad  q=2j-m-n/2+2.
\end{equation}
Whenever $m$ is even,
\begin{equation}\label{eqo5.13}
{\mathcal L}^{m,n}_o(0,-q(q+n-2))=\prod_{p}\left(q(q+n-2)-p\,(p+n-2)\right)^2,\quad q\in\NN\cup\{0\},
\end{equation}

\noindent where the product above is over $p=-n/2+2, -n/2+4,..., m-n/2-2, m-n/2$, and if $m$ is odd,
\begin{equation}\label{eqo5.14}
{\mathcal L}^{m,n}_o(0,-q(q+n-2))=\prod_{p}\left(q(q+n-2)-p\,(p+n-2)\right)^2 \left(q(q+n-2)+(n/2-1)^2\right),\end{equation}

\noindent  for every $q\in\NN\cup\{0\}$, with the product above over $p=-n/2+3, -n/2+5,..., m-n/2-2, m-n/2$.

\begin{lemma}\label{lo6.1} Assume that $m\in\NN$, $n\in [2,2m]\cap \NN$
is even and $m-n/2$ is even. Recall that in this case
\begin{equation}
\Lmno(-\partial_t,0)=\prod_{j=0}^{m-1} \left(-\partial_t^2+\Bigl(m-\frac n2-2j\Bigr)^2\right)
 \label{eqo6.1}
\end{equation}

\noindent and consider the equation
\begin{equation}
\prod_{j=0}^{m-1} \left(-\partial_t^2+\Bigl(m-\frac n2-2j\Bigr)^2\right)h
 =\delta,  \label{eqo6.2}
\end{equation}

\noindent where $\delta$ stands for the Dirac delta function. A
unique solution to {\rm (\ref{eqo6.2})} which 
vanishes at $+\infty$ and has at most linear growth or decay at $-\infty$ has a form
\begin{equation}\label{eqo6.3}
h(t)=\left\{\begin{array}{l}
\sum_{i=1}^{(m-n/2)/2} \nu_i^{(1)} \,e^{-2i t}+\sum_{i=1}^{(m-n/2)/2} \nu_i^{(2)} \,te^{-2i t}  + \sum_{i=1}^{n/2-1} \nu_i^{(3)} \,e^{-(m-n/2+2i) t},\qquad \qquad\,\, t>0,\\[6pt]
\sum_{i=1}^{(m-n/2)/2} \mu_i^{(1)} \,e^{2i t}+\sum_{i=1}^{(m-n/2)/2} \mu_i^{(2)} \,te^{2i t}  + \sum_{i=1}^{n/2-1} \mu_i^{(3)} \,e^{(m-n/2+2i) t}+\mu^{(4)}t+\mu^{(5)},\,\, t<0.\\[6pt]
\end{array}
\right.
\end{equation}

\noindent Here $\nu_i^{(1)}, \nu_i^{(2)}, \mu_i^{(1)} , \mu_i^{(2)}$, $i=1,...,(m-n/2)/2$, $\nu_i^{(3)}, \mu_i^{(3)}$, $i=1,...,n/2-1$, and $\mu^{(4)}, \mu^{(5)}$ are some real numbers depending on $m$ and $n$ only.

\end{lemma}

\begin{lemma}\label{lo6.2} Assume that $m\in\NN$, $n\in [2,2m]\cap \NN$
is even and $m-n/2$ is odd. Recall that in this case
\begin{equation}
\Lmno(-\partial_t,1-n)=\Lmno(-\partial_t, -1(1+n-2))=\prod_{j=0}^{m-1} \left(-\partial_t^2+\Bigl(m+\frac n2-2j-1\Bigr)^2\right)
 \label{eqo6.8}
\end{equation}

\noindent and consider the equation
\begin{equation}
\prod_{j=0}^{m-1} \left(-\partial_t^2+\Bigl(m+\frac n2-2j-1\Bigr)^2\right)h
 =\delta,  \label{eqo6.9}
\end{equation}

\noindent where $\delta$ stands for the Dirac delta function. A
unique solution to {\rm (\ref{eqo6.9})} which 
vanishes at $+\infty$ and has at most linear growth or decay at $-\infty$ has a form
\begin{equation}\label{eqo6.10}
h(t)=\left\{\begin{array}{l}
\sum_{i=1}^{(m-n/2-1)/2} \nu_i^{(1)} \,e^{-2i t}+\sum_{i=1}^{(m-n/2-1)/2} \nu_i^{(2)} \,te^{-2i t}  + \sum_{i=1}^{n/2} \nu_i^{(3)} \,e^{-(m-n/2-1+2i) t},\qquad \qquad\,\, t>0,\\[6pt]
\sum_{i=1}^{(m-n/2-1)/2} \mu_i^{(1)} \,e^{2i t}+\sum_{i=1}^{(m-n/2-1)/2} \mu_i^{(2)} \,te^{2i t}  + \sum_{i=1}^{n/2} \mu_i^{(3)} \,e^{(m-n/2-1+2i) t}+\mu^{(4)}t+\mu^{(5)},\,\, t<0.\\[6pt]
\end{array}
\right.
\end{equation}

\noindent Here $\nu_i^{(1)}, \nu_i^{(2)}, \mu_i^{(1)} , \mu_i^{(2)}$, $i=1,...,(m-n/2-1)/2$, $\nu_i^{(3)}, \mu_i^{(3)}$, $i=1,...,n/2$, and $\mu^{(4)}, \mu^{(5)}$ are some real numbers depending on $m$ and $n$ only.

\end{lemma}

\begin{theorem}\label{to6.3} Assume that $m\in\NN$ and  $n\in [2,2m]\cap \NN$
is even. Let $\Omega$ be a bounded domain in $\RR^n$, $O\in
\RR^n\setminus\Omega$, $u\in C_0^\infty(\Omega)$ and
$v=e^{\left(m-\frac n2\right)t}(u\circ \varkappa^{-1})$. Let $R$ be a positive constant such that the support of $u$ is contained in $B_{2R}$.
Then there exist positive constants $C,$ $C'$, $C''$, depending on $m$ and $n$ only, such that for
every $\xi\in B_{2R}$ and $\tau=\log |\xi|^{-1}$ we have
\begin{eqnarray}\label{eqo6.14}
\int_{S^{n-1}}v^2(\tau,\omega)\,d\omega \leq  C
\int_{\RR^n}(-\Delta)^m u(x)\,u(x)g(\log
|x|^{-1}, \log |\xi|^{-1})\,dx
\end{eqnarray}

\noindent  where $C_R=\log (4R)$ and $g$ is defined by 
\begin{equation}\label{eqo6.15}
g(t,\tau) = h(t-\tau)+\mu^{(4)}(C_R+\tau)+C'+C''(C_R+t),
\end{equation}
where $h$ and $\mu^{(4)}$ are given by \eqref{eqo6.3} and \eqref{eqo6.10}, depending on the parity of $m-\frac n2.$

\end{theorem}

\begin{proposition}\label{po7.2}
Let $\Omega$ be a bounded domain in $\RR^n$, $2\leq n \leq 2m+1$,
$Q\in\RR^n\setminus\Omega$, and $R>0$. Suppose
\begin{equation}\label{eqo7.3}
(-\Delta)^m u=f \,\,{\mbox{in}}\,\,\Omega, \quad f\in
C_0^{\infty}(\Omega\setminus B_{4R}(Q)),\quad u\in \ring
W^{m,2}(\Omega).
\end{equation}

\noindent Then
\begin{equation}\label{eqo7.4}
\frac{1}{\rho^{2\lambda +n-1}}\int_{S_{\rho}(Q)\cap\Omega}|u(x)|^2\,d\sigma_x
\leq \frac{C}{R^{2\lambda +n}} \int_{C_{R,4R}(Q)\cap\Omega} |u(x)|^2\,dx\quad
{\mbox{ for every}}\quad \rho<R,
\end{equation}

\noindent where $C$ is a  constant depending on $m$ and $n$ only, and  $\lambda=\left[m-\frac n2+\frac 12\right]$ as in \eqref{eqo7.5}.

Moreover, for every $x\in B_{R/4}(Q)\cap\Omega$
\begin{equation}\label{eqo7.15}
|\nabla^{i} u(x)|^2\leq C\, \frac{|x-Q|^{2\lambda-2i}}{R^{n+2\lambda}}\int_{C_{R/4,4R}(Q)\cap\Omega}
|u(y)|^2\,dy,\qquad 0\leq i\leq \lambda.
\end{equation}

In particular, for every bounded domain $\Omega\subset\RR^n$ the
solution to the boundary value problem {\rm (\ref{eqo7.3})}
satisfies
\begin{equation}\label{eqo7.16}
|\nabla^{m-n/2+1/2} u|\in L^\infty(\Omega)\,\,\mbox{when $n$ is odd\quad and \quad} |\nabla^{m-n/2} u|\in L^\infty(\Omega)\,\,\mbox{when $n$ is even}.
\end{equation}

\end{proposition}

We also address  the behavior of solutions ``at infinity". 
\begin{proposition}\label{po7.4} Let $\Omega$ be a bounded domain in $\RR^n$, $2\leq n \leq 2m+1$,
$Q\in\RR^n\setminus\Omega$, $r>0$ and assume that
\begin{equation}\label{eqo7.21}
(-\Delta)^m u=f \,\,{\mbox{in}}\,\,\Omega, \quad f\in
C_0^{\infty}(B_{r/4}(Q)\cap\Omega),\quad u\in \ring W^{m,2}(\Omega).
\end{equation}

\noindent Then
\begin{equation}\label{eqo7.22}
\rho^{2\lambda+n+1-4m}\int_{S_{\rho}(Q)\cap\Omega}|u(x)|^2\,d\sigma_x
\leq C\, r^{2\lambda+n-4m}\int_{ C_{r/4,r}(Q)\cap\Omega} |u(x)|^2\,dx,
\end{equation}

\noindent for any $\rho>r$ and $\lambda$ given by \eqref{eqo7.5}.

Furthermore, for any $x\in\Omega\setminus B_{4r}(Q)$
\begin{equation}\label{eqo7.23}
|\nabla^i u(x)|^2\leq
C\,\frac{r^{2\lambda+n-4m}}{|x-Q|^{2\lambda+2n-4m+2i}}\int_{C_{r/4,4r}(Q)\cap\Omega} |u(y)|^2\,dy,\qquad 0\leq i\leq \lambda.
\end{equation}
\end{proposition}

This concludes the list of main Theorems from \cite{MayMazInvent2} that will be used in the present paper. 

Finally, for future reference we record a  well-known result following from the energy estimates for solutions
of elliptic equations. 

\begin{lemma}\label{lo3.1} Let $\Omega$ be an arbitrary
domain in $\RR^n$, $n\geq 2$, $Q\in\RR^n\setminus\Omega$ and $R>0$. Suppose
\begin{equation}\label{eqo3.1}
(-\Delta)^m u=f \,\,{\mbox{in}}\,\,\Omega, \quad f\in
C_0^{\infty}(\Omega\setminus B_{4R}(Q)),\quad u\in \ring
W^{m,2}(\Omega).
\end{equation}

\noindent Then
\begin{equation}\label{eqo3.2}
\sum_{i=1}^m\frac{1}{\rho^{2m-2i}}\int_{B_{\rho}(Q)\cap\Omega}|\nabla^i u|^2\,dx
\leq
\frac{C}{\rho^{2m}}\int_{C_{\rho,2\rho}(Q)\cap\Omega}|u|^2\,dx
\end{equation}

\noindent for every $\rho<2R$.
\end{lemma}

\section{Higher order regularity of a boundary point as a local property}
\setcounter{equation}{0}

Let us recall from the introduction, \eqref{eq6.1}--\eqref{eq6.2}, the notion of $k$-regularity of a boundary point for an arbitrary $k\in\NN$. To start, we would like to show that $k$-regularity is a local
property. In other words, the $k$-regularity of a boundary point depends exclusively on the geometry of the domain near this point, rather than the geometry of the entire domain.

\begin{proposition}\label{p6.1}
Let $\Omega$ be a bounded domain in $\RR^n$, $m\in\NN$, $2\leq n \leq 2m+1$, and the point
$Q\in\po$ be $\lambda$-regular with respect to $\Omega$ and the operator $(-\Delta)^m$, with $\lambda\in \NN$ given by \eqref{eqo7.5}. If $\Omega'$ is
another domain (possibly unbounded) with the property that $B_r(Q)\cap\Omega =
B_r(Q)\cap\Omega'$ for some $r>0$ then $Q$ is $\lambda$-regular with
respect to $\Omega'$ as well.
\end{proposition}

The proof of the proposition rests on the ideas from \cite{M2}. Let us start from the following corollaries of the results in Section~\ref{srecall}. 

\begin{lemma}\label{Lim1} 
Let $\Omega$ be a bounded domain in $\RR^n$, $m\in\NN$, $2\leq n \leq 2m+1$,
$Q\in\RR^n\setminus{\overline{\Omega}}$. Then for every $u\in \ring
W^{m,2}(\Omega)$, and every $\rho>0$
\begin{eqnarray}\label{eqLim1}
&& \frac{1}{\rho^{2\lambda +n}}\int_{C_{\rho/2,\rho}(Q)\cap\Omega}|u(x)|^2\,dx\\[4pt]
\nonumber
&& \qquad \leq C\,\sup_{\xi\in C_{\rho/2,\rho}(Q)\cap\Omega}
\int_{\Omega}(-\Delta)^m u(y)\,u(y)g(\log |y-Q|^{-1},\log |\xi-Q|^{-1})\,dy, 
\end{eqnarray}

\noindent where $C$ is a  constant depending on $m$ and $n$ only, $g$ is the
function defined in (\ref{eqo4.12}) when $n$ is odd and by \eqref{eqo6.15} when $n$ is even, and 
$\lambda$ is given by \eqref{eqo7.5}.
\end{lemma}

\bp The Lemma follows from Theorems~\ref{to4.2}, \ref{to6.3} and a limiting procedure. Indeed, if $u\in C_0^{\infty}(\Omega)$, the desired estimate is literally the result of the aforementioned lemmas. Furthermore, when $u\in \ring W^{m.2}(\Omega)$ and, respectively, $(-\Delta)^m u\in W^{-m,2}(\Omega)$, the same result can be obtained approximating $u\in \ring W^{m.2}(\Omega)$ by $C_0^\infty(\Omega)$ functions in $\ring W^{m.2}(\Omega)$ norm and using the fact that for all $\xi, y\in \Omega$ the function $g(\log |y-Q|^{-1},\log |\xi-Q|^{-1})$ together with all its derivatives is bounded by a constant depending on the distance from $\overline{\Omega}$ to $Q$, $\rho$, (and ${\rm diam}\,(\Omega)$ in the case of even dimension) only. 

Note that the constant $C$ in the resulting inequality \eqref{eqLim1} does not depend on ${\rm dist}(Q, {\overline{\Omega}})$ (nor on ${\rm diam}\,(\Omega)$), the separation is only needed to justify the limiting procedure.  \ep

\begin{lemma}\label{Lim2} 
Let $\Omega$ be a bounded domain in $\RR^n$, $2\leq n \leq 2m+1$,
$Q\in\RR^n\setminus \Omega$. Suppose $u\in \ring
W^{m,2}(\Omega)$ is such that 
\begin{equation}\label{eq6.4-0}
(-\Delta)^m u=\sum_{\alpha:\,|\alpha|\leq m} \partial^\alpha f_\alpha
\,\,{\mbox{in}}\,\,\Omega, \quad f_\alpha\in L^2(\Omega)\cap
C^\infty(\Omega),\quad f_\alpha=0\mbox{ in a neighborhood of } Q.
\end{equation}

Then the estimate \eqref{eqLim1} is valid for every $\rho>0$
with the constant $C$ depending on $m$ and $n$ only, $g$ defined in (\ref{eqo4.12}) when $n$ is odd and by \eqref{eqo6.15} when $n$ is even, and 
$\lambda$ is given by \eqref{eqo7.5}.
\end{lemma}

\bp To fix the notation, let $R>0$ be such that ${\rm supp}\,f_\alpha\cap B_{4R}(Q)=\emptyset$, $|\alpha|\leq m$. Furthermore, let
$\left\{\Omega_n\right\}_{n=1}^\infty$ be a sequence of domains approximating $\Omega$ and staying away from $Q$, i.e., such that 
\begin{equation}\label{eq3.5.1}
\bigcup_{n=1}^\infty \Omega_n=\Omega, \quad 
{\Omega}_n\subset\Omega_{n+1},  \quad \mbox{and}\quad Q\not\in {\overline{\Omega_n}} \quad \mbox{for every}
\quad n\in\NN.
\end{equation}

\noindent Moreover, $\Omega_n\cap (\RR^n\setminus B_{2R}(Q))=\Omega\cap (\RR^n\setminus B_{2R}(Q))$ for all $n$. For instance, $\Omega_n$ can be built from $\Omega$ by cutting off small balls around $Q$. Let us denote by $u_n$ the solutions to the problem 
\begin{equation}\label{eq6.4-0}
(-\Delta)^m u_n=\sum_{\alpha:\,|\alpha|\leq m} \partial^\alpha f_\alpha|_{\Omega_n},
\,\,{\mbox{in}}\,\,\Omega_n, \quad u_n\in \ring
W^{m,2}(\Omega_n), 
\end{equation}

\noindent where $f_\alpha|_{\Omega_n}$ is a restriction of $f_\alpha$ to $\Omega_n$. Hence, by Lemma~\ref{Lim1} the estimate \eqref{eqLim1} holds for each $u_n$, and moreover, due to the restriction on the support of $f_\alpha$ and the construction, we have 
\begin{eqnarray}\label{eqLim1.1}
&& \frac{1}{\rho^{2\lambda +n}}\int_{C_{\rho/2,\rho}(Q)\cap\Omega}|u_n(x)|^2\,dx\nonumber\\[4pt]
&& \qquad \leq C\,\sup_{\xi\in C_{\rho/2,\rho}(Q)\cap\Omega}
\int_{\Omega\setminus{B_{2R}(Q)}}\sum_{\alpha:\,|\alpha|\leq m}  f_\alpha(y)\,(-\partial^\alpha)\left(u_n(y)g(\log |y-Q|^{-1},\log |\xi-Q|^{-1})\right)\,dy. 
\end{eqnarray}

\noindent On such a domain, the function $g$ together with all its derivatives of order less than or equal to $m$ is bounded by a constant depending $\rho$, $R$ only (and ${\rm diam}\,(\Omega)$ in the case of even dimension), and hence, we can pass to the limit on both sides of \eqref{eqLim1.1}  and conclude the proof of the Lemma as soon as we establish that $\{u_n\}$ converges to $u$ in $\ring W^{m,2}(\Omega)$. This is a consequence of the definition \eqref{eq6.4-0}. 

Indeed, first of all, for every $\varphi \in \ring W^{m,2}(\Omega)$ there is a sequence of $C_0^\infty(\Omega)$ functions approximating $\varphi$ in $\ring W^{m,2}(\Omega)$ norm, and hence, by \eqref{eq3.5.1}, we can choose from it a subsequence $\varphi_n\in C_0^\infty(\Omega_n)$ approximating $\varphi$ in $\ring W^{m,2}(\Omega)$ norm. For every such $\varphi_n$
\begin{equation}\label{eqLim1.2}
\int_{\Omega} (-\Delta)^m (u_n(y)-u(y))\,\varphi_n(y)\,dy=\int_{\Omega_n} (-\Delta)^m (u_n(y)-u(y))\,\varphi_n(y)\,dy=0.
\end{equation}

\noindent Indeed, $\varphi_n\in C_0^\infty(\Omega_n)$, and hence, the integration in \eqref{eqLim1.2} is over a domain strictly contained in $\Omega_n$, while $(-\Delta)^m u_n=(-\Delta)^m u$ on $\Omega_n$ by definition. On the other hand, 
\begin{equation}\label{eqLim1.3}
\lim_{n\to\infty}\int_{\Omega} (-\Delta)^m (u_n(y)-u(y))\,(\varphi_n(y)-\varphi(y))\,dy=0.
\end{equation}

\noindent This follows from the fact that $\varphi_n$ converges to $\varphi$ in $\ring W^{m,2}(\Omega)$ and $u_n-u \in \ring W^{m,2}(\Omega)$, with the uniformly bounded norms: 
\begin{eqnarray}\label{eqLim1.4}\nonumber
&& \|u_n\|_{\ring W^{m,2}(\Omega)}=\|u_n\|_{\ring W^{m,2}(\Omega_n)}\leq C \left\|\sum_{\alpha:\,|\alpha|\leq m} \partial^\alpha f_\alpha|_{\Omega_n}\right\|_{W^{-m,2}(\Omega_n)} \\[4pt]
&&\qquad =\sup_{v\in \ring W^{m,2}(\Omega_n): \,\|v\|_{\ring W^{m,2}(\Omega_n)}=1}\int_{\Omega_n}v(y) \sum_{\alpha:\,|\alpha|\leq m} \partial^\alpha f_\alpha|_{\Omega_n}(y)\,dy \nonumber \\[4pt]
&&\qquad \leq \sup_{v\in \ring W^{m,2}(\Omega):\,\|v\|_{\ring W^{m,2}(\Omega)}=1}\int_{\Omega}v(y) \sum_{\alpha:\,|\alpha|\leq m} \partial^\alpha f_\alpha(y)\,dy=\|f\|_{W^{-m,2}(\Omega)}.
\end{eqnarray}

 Combining \eqref{eqLim1.3} with \eqref{eqLim1.2}, we deduce that 
 \begin{equation}\label{eqLim1.5}
\lim_{n\to\infty}\int_{\Omega} (-\Delta)^m (u_n(y)-u(y))\,\varphi(y)\,dy=0, 
\end{equation}

\noindent for every $\varphi \in \ring W^{m,2}(\Omega)$. In particular,  we have 
 \begin{eqnarray}\label{eqLim1.6}
&&\int_{\Omega} (-\Delta)^m (u_n(y)-u(y))\,(u_n(y)-u(y))\,dy\nonumber\\[4pt]
&&\qquad =\int_{\Omega} (-\Delta)^m (u_n(y)-u(y))\,u_n(y)\,dy - \int_{\Omega} (-\Delta)^m (u_n(y)-u(y))\,u(y)\,dy \nonumber\\[4pt]
&&\qquad = 0-\int_{\Omega} (-\Delta)^m (u_n(y)-u(y))\,u(y)\,dy \longrightarrow 0, \quad \mbox{as $n\to\infty$,}
\end{eqnarray}

\noindent and hence, by ellipticity, the limit in $n$ of $\|u_n-u\|_{\ring W^{m,2}(\Omega)}$ is equal to 0, as desired. \ep

%\todo{The Lemma above will be substantially reduced in the final paper, in particular, the argument proving that $\{u_n\}$ converges to $u$ in $\ring W^{m,2}(\Omega)$ can perhaps be omitted as standard (see also Necas, Chapter 3, Theorem 6.7). I just wanted to make sure that it's actually true (and in fact, haven't seen approximation in domain all that often).}

\begin{lemma}\label{l6.2}
Let $\Omega$ be a bounded domain in $\RR^n$, $m\in\NN$, $2\leq n \leq 2m+1$,  and the point
$Q\in\po$ be $\lambda$-regular with respect to $\Omega$ and the operator $(-\Delta)^m$, with  $\lambda\in \NN$ given by \eqref{eqo7.5}. Then
\begin{equation}\label{eq6.3}
\nabla^{\lambda} u(x)\to 0\mbox{ as } x\to Q,\,x\in\Omega,
\end{equation}

\noindent for every $u\in \ring W^{m,2}(\Omega)$ satisfying
\begin{equation}\label{eq6.4}
(-\Delta)^m u=\sum_{\alpha:\,|\alpha|\leq m} \partial^\alpha f_\alpha
\,\,{\mbox{in}}\,\,\Omega, \quad f_\alpha\in L^2(\Omega)\cap
C^\infty(\Omega),\quad f_\alpha=0\mbox{ in a neighborhood of } Q.
\end{equation}
\end{lemma}

\bp Let us fix $\eps>0$ and take some $\eta_{\eps}\in C_0^\infty(\Omega)$ to be specified later. Let $v_\eps$ be the
solution of the Dirichlet problem
\begin{equation}\label{eq6.5}
(-\Delta)^m v_{\eps}=\sum_{\alpha:\,|\alpha|\leq m} \partial^\alpha (\eta_\eps
f_\alpha) \,\,{\mbox{in}}\,\,\Omega, \quad v_\eps\in \ring W^{m,2}(\Omega),
\end{equation}

\noindent and $w_\eps:=u-v_\eps \in  \ring W^{m,2}(\Omega)$. Since the point $Q$ is
$\lambda$-regular and the right-hand side of \eqref{eq6.5} belongs to $C_0^{\infty}(\Omega)$, the function $v_\eps$ automatically satisfies (\ref{eq6.3}) and in particular, there exists $\delta_1=\delta_1(\eps)>0$ such that 
\begin{equation}\label{eq6.5.1}
|\nabla^{\lambda} v_\eps(x)|<\eps/2 \mbox{ whenever } |x-Q|<\delta_1.
\end{equation}

Let us now consider $w_\eps\in \ring W^{m,2}(\Omega)$. It satisfies the equation 
\begin{equation}\label{eq6.5.2}
(-\Delta)^m w_{\eps}=\sum_{\alpha:\,|\alpha|\leq m} \partial^\alpha ((1-\eta_\eps)
f_\alpha) \,\,{\mbox{in}}\,\,\Omega,\quad \mbox {with}  \quad f_\alpha=0  \mbox{ in a neighborhood of } Q.
\end{equation}

As before, let $R>0$ be such that ${\rm supp}\,f_\alpha\cap B_{4R}(Q)=\emptyset$, $|\alpha|\leq m$. Then, by interior estimates for solutions of elliptic PDE \cite{ADN},
\begin{equation}\label{eq6.6}
|\nabla^{\lambda} w_\eps (x)|^2\leq  \frac{C}{d(x)^{n+2\lambda}} \int_{B_{d(x)/2}(x)} |w_\eps (y)|^2\,dy ,\quad \forall\,\,x\in B_{R/4}(Q),
\end{equation}
$d(x)$ denoting the distance to $\po$.

Now for every $x\in B_{R/4}(Q)$  denote by $x_0$ a point on $\po$ such that $d(x)=|x-x_0|$ and note that   ${\rm supp}\,f_\alpha\cap B_{3R}(x_0)=\emptyset$, $|\alpha|\leq m$. Then, according to Lemma~\ref{Lim2},
\begin{multline}\label{eq6.7}
 \frac{C}{d(x)^{n+2{\lambda}}} \int_{B_{d(x)/2}(x)} |w_\eps (y)|^2\,dy\leq \frac{C}{|x-x_0|^{n+2{\lambda}}} \int_{C_{\frac{|x-x_0|}{2}, \frac{3|x-x_0|}{2}}(x_0)} |w_\eps (y)|^2\,dy\\[4pt]
 \quad \leq C\,\sup_{\xi\in C_{\frac{|x-x_0|}{2}, \frac{3|x-x_0|}{2}}(x_0)}
\int_{\RR^n}(-\Delta)^m w_\eps(y)\,w_\eps(y)g(\log |y-x_0|^{-1},\log |\xi-x_0|^{-1})\,dy\\[4pt]
 \quad \leq C\,\sup_{\xi\in C_{\frac{|x-x_0|}{2}, \frac{3|x-x_0|}{2}}(x_0)}
\sum_{\alpha:\,|\alpha|\leq m}
\int_{\RR^n}(1-\eta_\eps(y))f_\alpha(y)\,\,\times \\[4pt]
\qquad  \qquad\quad\times\, (-\partial_y)^\alpha
\bigg(w_\eps(y)g(\log |y-x_0|^{-1},\log |\xi-x_0|^{-1})\bigg)\,dy.
\end{multline}

Since the supports of $f_\alpha$ stay away from $B_{3R}(x_0)$, the function $g$ together with its derivatives is bounded by a constant depending on $\supp f_\alpha$ $R$ only (and ${\rm diam}\,(\Omega)$ in the case of even dimension). Moreover, the $L^2$ norms of $(1-\eta_\eps)f_\alpha$ vanish as the euclidean size of support of $1-\eta_\eps$ goes to 0, and therefore, by \eqref{eq6.5.2}, $\|w_\eps\|_{\ring W^{m,2}(\Omega)}$ vanish as well.  Then, we can choose $\eta_\eps$ so that the expression in \eqref{eq6.7} and,  hence, \eqref{eq6.6}, does not exceed $\eps/2$. Combined with \eqref{eq6.5.1}, this leads to \eqref{eq6.3}, as desired. 
\ep

\vskip 0.08in

\noindent {\it Proof of Proposition~\ref{p6.1}.}\, Consider a
solution of the Dirichlet problem
\begin{equation}\label{eq6.8}
(-\Delta)^m u=f \,\,{\mbox{in}}\,\,\Omega', \quad f\in
C_0^{\infty}(\Omega'),\quad u\in \ring W^{m,2}(\Omega'),
\end{equation}

\noindent and take some cut-off function $\eta\in
C_0^\infty(B_{r}(Q))$ equal to 1 on $B_{r/2}(Q)$. Then $\eta u \in
 W^{m,2}(\Omega)$ and
\begin{equation}\label{eq6.9}
(-\Delta)^m (\eta u)=\eta f+[(-\Delta)^m,\eta] u.
\end{equation}

Since $\eta f\in C_0^\infty(\Omega)$ and the commutator above is a differential operator of order $2m-1$ with smooth coefficients supported in $C_{r/2,r}(Q)$, one can write
\begin{equation}\label{eq6.11}
(-\Delta)^m(\eta u)= \sum_{\alpha:\,|\alpha|\leq m} \partial^\alpha
f_\alpha, \quad\mbox{for some}\quad f_\alpha\in L^2(\Omega)\cap
C^\infty(\Omega),
\end{equation}

\noindent with $f_\alpha=0$ in a neighborhood of $Q$ given by the
intersection of $B_{r/2}(Q)$ and the complement to ${\rm
supp}\,f$. Then, by Lemma~\ref{l6.2}, the $\lambda$-th gradient of $\eta u$
(and therefore,  $\nabla^{\lambda}u$) vanishes as $x\to Q$. \ep

\section{The new notion of polyharmonic capacity} \setcounter{equation}{0}

{\it Throughout this Section $n\in [2, 2m+1]$. Assume first that $m\in\NN$ and $n\in [3,2m+1]\cap \NN$ is odd. }

We start with the observation that according to the computations that we recalled in Section~\ref{srecall}, for any bounded domain $\Omega$  in $\RR^n$, $O\in
\RR^n\setminus\Omega$, $u\in C_0^\infty(\Omega)$ and
$v=e^{\left(m-\frac n2+\frac 12\right)t}(u\circ \varkappa^{-1})$ the expression 
\begin{equation}\label{cap1}
\int_{\RR^n}(-\Delta)^m u(x)\,u(x)|x|^{-1}\,dx
\end{equation}

\noindent can be written as 
\begin{eqnarray}\label{cap2}\nonumber
&& \sum_{p=0}^\infty \sum_{l=-p}^{p} \int_{\RR} {\mathcal
L}^{m,n}(\partial_t,-p\,(p+n-2)) v_{pl}(t) \,v_{pl}(t)\,dt\\[4pt]
&& \qquad \nonumber = \sum_{p=0}^\infty \sum_{l=-p}^{p} \int_{\RR} \Re
e\, {\mathcal L}^{m,n}(i\gamma,-p\,(p+n-2))
\left|\widehat{v_{pl}}(\gamma)\right|^2\,d\gamma\\[4pt]
&& \qquad 
=:\sum_{p=0}^\infty \sum_{l=-p}^{p} \int_{\RR} a_{p+\frac n2-\frac 32}(\gamma)
\left|\widehat{v_{pl}}(\gamma)\right|^2\,d\gamma= \sum_{p=0}^\infty \sum_{l=-p}^{p} \int_{\RR} a_{p+\frac n2-\frac 32}(\partial_t)
v_{pl}(t) \,v_{pl}(t)\,dt,
\end{eqnarray}

\noindent where for $m$ and $n$ as above ${\mathcal L}^{m,n}$ are defined in \eqref{eqo2.6}, $v_{pl}$ are the coefficients of the decomposition of $v$ into spherical harmonics \eqref{eqo2.7}, and at the moment the reader can understand $a_{p+\frac n2-\frac 32}(\gamma)$ simply as the notation for the coefficients of $ \Re
e\, {\mathcal L}^{m,n}(i\gamma,-p\,(p+n-2))
$ appearing on the second line of \eqref{cap2}. The notation $a_p(\partial_t)$ used in the last line of \eqref{cap2} stands for the differential operator formed by replacing $i\gamma$ by $\partial_t$ in the representation of $a_p$.
The exact formulas for $a_p$ are computed in the proof of Theorem~\ref{to2.1} in \cite{MayMazInvent2} and we do not copy them here as they play no particular role in the present argument. However, they enjoy a few crucial estimates which we will extensively use. 
All $a_p$, $p\in\NN\cup\{0\}$, are the polynomials of $\gamma^2$ of order $m$. Moreover, according to Step VIII in the proof of Theorem~2.1, in \cite{MayMazInvent2}, {\it loc.cit.}, we have 
\begin{equation}\label{cap2.1}
a_{p +\frac n2-\frac 32}(\gamma)\geq C\sum_{k=1}^m \gamma^{2k} +C \prod_{s=-\frac n2+\frac
32}^{m-\frac n2+\frac 12}\Bigl(p(p+n-2)-s(s+n-2)\Bigr),
\end{equation}
\noindent for some $C>0$ and $p\in\NN\cup\{0\}$. This is the translation of \eqref{eqo2.3.1}. 
Even more precisely, we can write 
\begin{equation}\label{cap3}
a_{p+\frac n2-\frac 32}(\gamma)=\sum_{k=0}^m c_{kp}\gamma^{2k}, \quad \gamma\in \RR, p\in\NN\cup\{0\}, c_{kp}\in \RR, 
\end{equation}

\noindent where $c_{kp}$ are such that, firstly, there exists a constant  $c_0>0$ such that
\begin{equation}\label{cap4.0}
c_{kp}\geq c_0 \,\mbox{ for }\, k\geq 1, p\in \NN\cup\{0\}, \end{equation}
and secondly,
\begin{equation}\label{cap4.0.0}
c_{0p}= C \prod_{s=-\frac n2+\frac
32}^{m-\frac n2+\frac 12}\Bigl(p(p+n-2)-s(s+n-2)\Bigr), \quad p\in \NN\cup\{0\}, \end{equation}
for some $C>0$, so that 
\begin{equation}\label{cap4.0.1}
c_{0p}=0 \,\mbox{ if } \,0\leq p\leq m-\frac n2+\frac 12, \quad c_{0p}\geq c_0 \,\mbox{ if }\, p\geq m-\frac n2+\frac 32,\end{equation}
and moreover, by \eqref{eqoo4.37},
\begin{equation}\label{cap4.1}
c_{0p}\geq c_0 \left(p(p+n-2)\right)^m \,\mbox{ if }\, p\geq m-\frac n2+\frac 32.\end{equation}
The fact that $c_{0p}$ have the exact form \eqref{cap4.0.0} is not explicitly discussed in the proof of Theorem~2.1, in \cite{MayMazInvent2}, but it can be seen directly by writing 
\begin{equation}\label{eqo4.29}{\mathcal
L}^{m,n}(0,\delta)=(-1)^m\prod_{j=0}^{m-1}
\Bigg(\Bigl(m-\frac n2+\frac
12-2j\Bigr)\Bigl(m+\frac n2-\frac
32-2j\Bigr)+\delta_\omega\Bigg)=\prod_{s=-\frac n2+\frac 32}^{m-\frac
n2+\frac 12}\left(-\delta_\omega-s\,(s+n-2)\right).
\end{equation}
We do not claim that the positive constant denoted by $C$ in \eqref{cap2.1} and in \eqref{cap4.0.0} is the same, and it is not important for the discussion. 

All in all, for any bounded domain $\Omega$  in $\RR^n$, $O\in
\RR^n\setminus\Omega$, $u\in C_0^\infty(\Omega)$ and
$v=e^{\left(m-\frac n2+\frac 12\right)t}(u\circ \varkappa^{-1})$ we have 
\begin{equation}\label{cap5}
\int_{\RR^n}(-\Delta)^m u(x)\,u(x)|x|^{-1}\,dx = \sum_{p=0}^\infty \sum_{l=-p}^{p} \sum_{k=0}^m c_{kp}\int_{\RR} 
(\partial_t^k v_{pl}(t))^2\,dt,
\end{equation}

\noindent where the coefficients $c_{kp}\in \RR$ are as above, in particular, they satisfy \eqref{cap4.0}--\eqref{cap4.1}.

{\it Now let us assume that $m\in\NN$ and $n\in [2,2m]\cap \NN$ is even.} As before, let $\Omega$ be a bounded domain in $\RR^n$, $O\in
\RR^n\setminus\Omega$, $u\in C_0^\infty(\Omega)$ and this time,
$v=e^{\left(m-\frac n2\right)t}(u\circ \varkappa^{-1})$. Retain the notation of Theorem~\ref{to5.1} and recall the comments right after its statement, \eqref{eqo5.3}--\eqref{eqo5.14}. 
All this can also be translated into the language of $c_{kp}$, i.e., the coefficients near various terms of decomposition into spherical harmonics on the Fourier transform side. We shall do that now and also combine this with the computations for odd dimensions presented above. 

\begin{definition} \label{defckp}
Let $\Omega$ be a bounded domain in $\RR^n$, $2\leq n\leq 2m+1$, $O\in
\RR^n\setminus\Omega$, $u\in C_0^\infty(\Omega)$ and 
$v=e^{\lambda t}(u\circ \varkappa^{-1})$, with $\lambda$ given by \eqref{eqo7.5}. The coefficients $c_{kp}\in\RR$ are uniquely determined by the condition that for $u$ and $v$ as above
\begin{equation}\label{capdef1}
\int_{\RR^n}(-\Delta)^m u(x)\,u(x)|x|^{2m-n-2\lambda}\,dx = \sum_{p=0}^\infty \sum_{l=-p}^{p} \sum_{k=0}^m c_{kp}\int_{\RR} 
(\partial_t^k v_{pl}(t))^2\,dt.
\end{equation}
Of course, $\{c_{kp}\}$'s depend on $m$ and $n$ as well, but we omit additional indices not to make the notation too cumbersome. 
Let furthermore $Z$ denote the set 
$$Z:=\{p\in \NN\cup\{0\}:\, c_{0p}=0\}.$$ According to the computations above and \eqref{eqo5.8}--\eqref{eqo5.9}, $Z$ coincides exactly with the set defined by \eqref{z1}--\eqref{z2}. For $p\not\in Z$, we have (see \eqref{cap4.1} and \eqref{eqo5.8})
\begin{equation}\label{capdef2}
c_{0p}\geq c_0 \max\{1,\left(p(p+n-2)\right)^m\} \,\mbox{ if }\, p\not\in Z,\end{equation}
and (see \eqref{cap4.0} and \eqref{eqo5.6}--\eqref{eqo5.8}) 
\begin{equation}\label{capdef3}
c_{kp}\geq c_0 \,\mbox{ for }\, k\geq 1, p\in \NN\cup\{0\}, \end{equation}
for some constant $c_0>0$.

To be more specific, we record that 
\begin{equation}\label{capdef4}
c_{op}=\Lmn(0,-p(p+n-2)), \,  \mbox{if $n$ is odd, and }c_{op}=\Lmno(0,-p(p+n-2)), \,\mbox{if $n$ is even},
\end{equation}
and the explicit expressions can be found in \eqref{cap4.0.0}, \eqref{eqo5.13}, \eqref{eqo5.14}. 
\end{definition}

\begin{definition} \label{defPhi}
Given a bounded domain $\Omega$  in $\RR^n$, $2\leq n\leq 2m+1$,  $O\in
\RR^n\setminus\Omega$, $u\in C_0^\infty(\Omega)$ and
$v=u\circ \varkappa^{-1}$ (note a different normalization!), we let 
\begin{equation}\label{cap6}
\Phi[u;\Omega]: = \sum_{p=0}^\infty \sum_{l=-p}^{p} \sum_{k=0}^m c_{kp}\int_{\RR} 
(\partial_t^k v_{pl}(t))^2\,dt.
\end{equation}
\end{definition}
In particular, again, due to the choice of the normalization,
\begin{equation}\label{cap5.1}
\Phi[u;\Omega]=\int_{\RR^n}(-\Delta)^m \left(|x|^{\lambda}u(x)\right)\, |x|^{2m-n-\lambda}u(x)\,dx\quad\mbox{whenever}\quad u\in C_0^\infty(\Omega), 
\end{equation}

\noindent and $\Omega$ is a bounded domain in $\RR^n$ not containing $O$.
 However, $\Phi[u;\Omega]$ is well-defined for a larger class of functions $u$ than the one in \eqref{cap6}--\eqref{cap5.1}, and we will further use the notation \eqref{cap6} in a more general setting. In this vein, for an annulus $C_{a,b}$, $0<a<b<\infty$, we write 
 \begin{equation}\label{cap6.0}
\Phi[u;C_{a,b}]: = \sum_{p=0}^\infty \sum_{l=-p}^{p} \sum_{k=0}^m c_{kp}\int_{\log\frac{1}{b}} ^{\log\frac{1}{a}}
(\partial_t^k v_{pl}(t))^2\,dt, 
\end{equation}

\noindent where $c_{kp}\in \RR$ are the same coefficients as above, $v=u\circ \varkappa^{-1}$ and $u$ is any function for which the expression in \eqref{cap6.0} is finite (with partial derivatives understood in the sense of distributions).

In the case of even dimensions, we shall also make use of the functional $\Phi$ with an additional logarithmic weight, 
\begin{equation}\label{cap6.0-log}
\Phi_R[u;\Omega]: = \sum_{p=0}^\infty \sum_{l=-p}^{p} \sum_{k=0}^m c_{kp}\int_{\RR}
(\partial_t^k v_{pl}(t))^2 (C_R+t)\,dt, \quad C_R:=\log (4R),
\end{equation}
where $R$ is taken so that the support of $u$ is contained in $B_{2R}$, $O\in \RR^n\setminus \Omega$.
Note that the subscript $R$ refers to the choice of weight $C_R+t$ and we loosely refer to the weight as ``logarithmic" subconsciously treating it back in Euclidean coordinates. 

As per \eqref{eqo5.3} and the discussion below it, 
\begin{equation}\label{cap5.1-log}
\Phi_R[u;\Omega]=\int_{\RR^n}(-\Delta)^m \left(|x|^{\lambda}u(x)\right)\, |x|^{2m-n-\lambda} \log (4R |x|^{-1})u(x)\,dx\quad\mbox{for}\quad u\in C_0^\infty(\Omega), \,\, \supp u\subset B_{2R}.
\end{equation}
In particular, it is not difficult to see that the coefficients $\{c_{kp\}}$ coming from the integration with the power weight and with the power-logarithmic weight (cf. \eqref{cap5.1}, \eqref{cap6} and \eqref{cap5.1-log}, \eqref{cap6.0-log}) are indeed the same. The reader can also consult the proof of Theorem~5.1 in \cite{MayMazInvent2} for details. 

Similarly to \eqref{cap6.0}, or an annulus $C_{a,b}$, $0<a<b<\infty$, we set  
 \begin{equation}\label{cap6.0-log1}
\Phi_R[u;C_{a,b}]: = \sum_{p=0}^\infty \sum_{l=-p}^{p} \sum_{k=0}^m c_{kp}\int_{\log\frac{1}{b}} ^{\log\frac{1}{a}}
(\partial_t^k v_{pl}(t))^2\,(C_R+t)\,dt, \quad C_R:=\log (4R),
\end{equation}

\noindent where $c_{kp}\in \RR$ are the same coefficients as above, $v=u\circ \varkappa^{-1}$ and $u$ is any function for which the expression in \eqref{cap6.0} is finite (with partial derivatives understood in the sense of distributions). As before, $R$ is taken so that the support of $u$ is contained in $B_{2R}$.

For future reference we would like to single out the following estimates.
\begin{lemma}\label{lPhi} Let $0<a<b<\infty$. Then 
 \begin{eqnarray}\label{cap6.0.1}
&&\Phi[u;C_{a,b}]\geq C\sum_{\stackrel{k\geq 1,\,i \geq 0}{ i+k\leq
m}}\int_{\log\frac{1}{b}} ^{\log\frac{1}{a}}
\int_{S^{n-1}} (\partial_t^k\nabla_\omega^i
v(t,\omega))^2
\,d\omega dt\nonumber\\[4pt]&&\qquad +\,C\int_{\log\frac{1}{b}} ^{\log\frac{1}{a}}
\int_{S^{n-1}} {\mathcal
L}^{m,n}(0,\delta) v(t,\omega)\,v(t,\omega)\,d\omega dt,\quad \mbox{for $n$ odd},
\end{eqnarray}
and  
\begin{eqnarray}\label{cap6.0-even}
&&\Phi[u;C_{a,b}]\geq C\sum_{\stackrel{k\geq 1,\,i \geq 0}{ i+k\leq
m}}\int_{\log\frac{1}{b}} ^{\log\frac{1}{a}}
\int_{S^{n-1}} (\partial_t^k\nabla_\omega^i
v(t,\omega))^2
\,d\omega dt\nonumber\\[4pt]&&\qquad +\,C\int_{\log\frac{1}{b}} ^{\log\frac{1}{a}}
\int_{S^{n-1}} {\mathcal
L}^{m,n}_o(0,\delta) v(t,\omega)\,v(t,\omega)\,d\omega dt,\quad \mbox{for $n$ even},
\end{eqnarray}
\noindent where $v=u\circ \varkappa^{-1}$ and $u$ is any function for which both sides of \eqref{cap6.0.1}, or, respectively, \eqref{cap6.0-even}, are finite. Moreover, with the same notation, 
\begin{multline}\label{cap6.0.2} s^{-n}\|u\|^2_{
L^2(C_{as,bs})}+\Phi[u; C_{as,bs}]\\[4pt]\approx \sum_{\stackrel{k\geq 0,\,i \geq 0}{ i+k\leq
m}}\int_{\log\frac{1}{b}} ^{\log\frac{1}{a}}
\int_{S^{n-1}} (\partial_t^k\nabla_\omega^i
v(t,\omega))^2
\,d\omega dt \approx \sum_{0\leq k\leq m}s^{2k-n}\int_{C_{as,bs}} |\nabla^k u|^2\,dx,
\end{multline}
with the constants independent of $s$. If, in addition, $u\in C_0^\infty(C_{as,bs})$, then also 
\begin{multline}\label{cap6.0.3} \Phi[u; C_{as,bs}] \approx \sum_{\stackrel{k\geq 0,\,i \geq 0}{ i+k\leq
m}}\int_{\log\frac{1}{b}} ^{\log\frac{1}{a}}
\int_{S^{n-1}} (\partial_t^k\nabla_\omega^i
v(t,\omega))^2
\,d\omega dt \\[4pt]
\approx \sum_{0\leq k\leq m}s^{2k-n} \int_{C_{as,bs}} |\nabla^k u|^2\,dx\approx s^{2m-n} \int_{C_{as,bs}} |\nabla^m u|^2\,dx.
\end{multline}
\end{lemma}

\bp Let us start with \eqref{cap6.0.1}, \eqref{cap6.0-even}. The estimate 
 \begin{equation}\label{cap6.0.4}
\Phi[u;C_{a,b}]\geq C\sum_{k=1}^m\int_{\log\frac{1}{b}} ^{\log\frac{1}{a}}
\int_{S^{n-1}} (\partial_t^k
v(t,\omega))^2
\,d\omega dt,
\end{equation}
follows directly from the definition and \eqref{capdef3}. The bound from below by 
$$\int_{\log\frac{1}{b}} ^{\log\frac{1}{a}}
\int_{S^{n-1}} {\mathcal
L}^{m,n}(0,\delta) v(t,\omega)\,v(t,\omega)\,d\omega dt$$ and $$\int_{\log\frac{1}{b}} ^{\log\frac{1}{a}}
\int_{S^{n-1}} {\mathcal
L}^{m,n}_o(0,\delta) v(t,\omega)\,v(t,\omega)\,d\omega dt,$$ in the case of odd and even dimension, respectively, is due to \eqref{capdef4}. Finally, the mixed derivatives, corresponding to 
$1\leq k$, $0\leq i$, $i+k\leq m$, are bounded by the combination of the right-hand side of \eqref{cap6.0.4} and the two integrals above, 
by the same argument as \eqref{eqoo4.35}--\eqref{eqoo4.40} employing the bound \eqref{capdef2}. 
 This finishes the proof of \eqref{cap6.0.1}, \eqref{cap6.0-even}.

In order to get $\geq$ in the first inequality in \eqref{cap6.0.2}, we now only need to bound $$\sum_{0\leq i\leq m}\int_{\log\frac{1}{b}} ^{\log\frac{1}{a}}
\int_{S^{n-1}} (\nabla_\omega^i
v(t,\omega))^2
\,d\omega dt=\sum_{0\leq i\leq m}\sum_{p=0}^{\infty} \sum_{l=-p}^{p}\int_{\log\frac{1}{b}} ^{\log\frac{1}{a}}
 \left(p(p+n-2)\right)^iv_{pl}(t)^2
\,dt.$$ 
This follows essentially from a simpler version of the argument in  \eqref{eqoo4.35}--\eqref{eqoo4.40}, now with $k=0$. Indeed, 
$$\sum_{p\in Z} \sum_{l=-p}^{p}
\Bigl(p\,(p+n-2)\Bigr)^i \int_{\log\frac{1}{b}} ^{\log\frac{1}{a}}( v_{pl}(t))^2
\,dt\leq C_{m,n}
\int_{\log\frac{1}{b}} ^{\log\frac{1}{a}}\int_{S^{n-1}} (v(t,\omega))^2 \,d\omega dt
 \leq C_{m,n} s^{-n} \|u\|^2_{L^2(C_{as, bs})},$$
\noindent and by \eqref{capdef2} we  have 
\begin{eqnarray}\label{eqoo4.40.1}
&&\sum_{p\not\in Z} \sum_{l=-p}^{p} \Bigl(p\,(p+n-2)\Bigr)^i
\int_{\log\frac{1}{b}} ^{\log\frac{1}{a}}(v_{pl}(t))^2 \,dt \nonumber\\[4pt] &&\qquad \lesssim \sum_{p\not\in Z}\sum_{l=-p}^{p}
\,\,c_{0p}
\int_{\log\frac{1}{b}} ^{\log\frac{1}{a}}(v_{pl}(t))^2 \,dt \leq C
\Phi[u;C_{as,bs}],
\end{eqnarray}
\noindent as desired. 

The $\lesssim$ inequality in \eqref{cap6.0.2} is straightforward from definitions, and the second equivalence in \eqref{cap6.0.2} is a change of variables.

The passage from \eqref{cap6.0.2} to \eqref{cap6.0.3} is due to Hardy inequality. Let us furnish some details.  We only have to show $\lesssim$ direction in the last equivalence in \eqref{cap6.0.3} and $\gtrsim$ part in the first equivalence in \eqref{cap6.0.3}. 
First, for every $u\in \ring W^{m,2}(C_{as, bs})$ we have 
\begin{equation}\label{cap12.1}\nonumber
s^{2m-n}\sum_{k=0}^m\frac{1}{s^{2m-2k}}\int_{C_{as, bs}} |\nabla^ku(x)|^2\,dx \leq Cs^{2m-n}\int_{C_{as, bs}} |\nabla^m u(x)|^2\,dx,
\end{equation}

\noindent using the Hardy inequality. This proves $\lesssim$ direction in the last equivalence \eqref{cap6.0.3}. On the other hand, by one-dimensional Hardy inequality in $t\in \left(\log \frac{1}{bs},\log \frac {1}{as}\right)$, we can bound   
\begin{eqnarray}\label{cap16.1}
\int_{\log\frac{1}{b}} ^{\log\frac{1}{a}}\int_{S^{n-1}} (v(t,\omega))^2 \,d\omega dt\leq C \int_{\log\frac{1}{b}} ^{\log\frac{1}{a}}\int_{S^{n-1}} (\partial^t v(t,\omega))^2 \,d\omega dt,
\end{eqnarray}
which in turn is bounded by $C\Phi[u;C_{as,bs}]$, as desired. \ep

\begin{remark}\label{rcap1} Let $s>0$, $0<a<b<\infty$. Motivated by \eqref{cap6.0.2}, we define the Sobolev-type space
\begin{multline}\label{defV}
V^{m,2}(C_{as,bs}):=\Big\{u\in L^2(C_{as,bs}):\,\|u\|_{V^{m,2}(C_{as,bs})}^2:=s^{-n}\|u\|^2_{
L^2(C_{as,bs})}+\Phi[u; C_{as,bs}]\\[4pt]\approx \sum_{\stackrel{k\geq 0,\,i \geq 0}{ i+k\leq
m}}\int_{\log\frac{1}{b}} ^{\log\frac{1}{a}}
\int_{S^{n-1}} (\partial_t^k\nabla_\omega^i
v(t,\omega))^2
\,d\omega dt \approx \sum_{0\leq k\leq m}s^{2k-n}\int_{C_{as,bs}} |\nabla^k u|^2\,dx<\infty\Big\},
\end{multline}

\noindent where the derivatives are, as usually, understood in the weak sense, and $v=u\circ \varkappa^{-1}$.

\end{remark}

\begin{remark}\label{rcap2} We remark that the first equivalence in \eqref{cap6.0.3} does not require the full power of the condition $u\in C_0^\infty(C_{as,bs})$ but rather a weak vanishing on the boundary sufficient for application of the Hardy inequality in \eqref{cap16.1}.
\end{remark}

Now let us recall from the introduction the space of linear combinations of spherical harmonics from \eqref{cap6.1}--\eqref{cap6.2}. Then, given $P\in\Pi_1$, for any bounded open set $\Omega$ such that $O\notin\overline\Omega$  we  define $\Phi$-capacity of a compactum $K\subset
\Omega$  as
\begin{eqnarray}\label{cap7}
{\rm Cap}_P^{\Phi}\,(K,\Omega):=\inf\Bigg\{\Phi[u;\Omega]:\,\,u\in \ring W^{m,2}(\Omega),\,\,u=P\mbox{ in a
neighborhood of }K\Bigg\}.
\end{eqnarray}
In the introduction we utilized a somewhat different polyharmonic capacity, namely, \eqref{cap8}--\eqref{cap9}. 

We shall use ${\rm Cap}_P^{\Phi}$ in place of ${\rm Cap}_P$ whenever convenient, and employ the following relation between the two of them.

\begin{lemma}\label{lcap1} Let $\Omega$ be an open set in an annulus $C_{as, bs}$ for some constants $0<a<b<\infty$. Then for every compactum $K\subset
\Omega$ and any $P\in \Pi_1$
\begin{equation}\label{cap10}
{\rm Cap}_P\,(K,\Omega)\approx s^{n-2m}\,{\rm Cap}_P^{\Phi}\,(K,\Omega),
\end{equation}

\noindent with the implicit constants depending on $n,m$ and $a,b$ only.
\end{lemma}

\bp  The proof is a direct application of \eqref{cap6.0.3}. For every $\eps>0$ we can choose $u\in \ring W^{m,2}(\Omega)$ such that $u=P$ in a neighborhood of $K$ and 
\begin{equation}\label{cap11}
\int_\Omega |\nabla^m u(x)|^2\,dx\leq {\rm Cap}_P\,(K,\Omega)+\eps.
\end{equation}

\noindent However, for every $u\in \ring W^{m,2}(C_{as, bs})$ (in our case, extended by zero from $\ring W^{m,2}(\Omega)$ to $\ring W^{m,2}(C_{as, bs}$) we have 
\begin{equation}\label{cap12}
\Phi[u;\Omega] \leq Cs^{2m-n}\int_{C_{as, bs}} |\nabla^m u(x)|^2\,dx,
\end{equation}

\noindent by \eqref{cap6.0.3}. This gives the $``\geq"$ inequality in \eqref{cap10}. The converse argument is exactly the same, again using  \eqref{cap6.0.3}. \ep

It follows directly from the definitions that all versions of the capacity defined above are monotone in the sense that for every $P\in\Pi_1$
\begin{eqnarray}\label{eq5.4}
&& K_1\subseteq K_2\subset \Omega \quad  \Longrightarrow  \quad{\rm
Cap}_P\,(K_1,\Omega)\leq {\rm Cap}_P\,(K_2,\Omega),\quad {\rm
Cap}_P^\Phi\,(K_1,\Omega)\leq {\rm Cap}_P^\Phi\,(K_2,\Omega),\\[4pt]\label{eq5.5}
&& K\subset \Omega_1\subseteq \Omega_2 \quad  \Longrightarrow 
\quad{\rm Cap}_P\,(K,\Omega_1)\geq {\rm Cap}_P\,(K,\Omega_2), \quad{\rm Cap}_P^\Phi\,(K,\Omega_1)\geq {\rm Cap}_P^\Phi\,(K,\Omega_2),
\end{eqnarray}

\noindent and analogous statements hold for ${\rm Cap}$. In addition, we have the following relations.

\begin{lemma}\label{l5.2}
Let $K$ be a compactum in $\overline{C_{as,bs}}$ for some
$s>0$, $0<a<b<\infty$. Then for every $P\in\Pi_1$
\begin{equation}\label{eq5.12.1}
{\rm Cap}_P(K, C_{as/2,2bs})\leq
Cs^{n-2m},\qquad {\rm Cap}_P^\Phi(K, C_{as/2,2bs})\leq
C,
\end{equation}

\noindent with the constants independent of $s$.
\end{lemma}

\bp The estimates \eqref{eq5.12.1} come from the scaling considerations. Indeed, 
if  $w(x)=u(sx)$, $x\in\RR^n$, then first of all, 
\begin{equation}\label{eq5.18.1}
u\in \ring
W^{m,2}(C_{as/2,2bs})\quad \Longleftrightarrow \quad w\in \ring
W^{m,2}(C_{a/2,2b}),
\end{equation}

\noindent and secondly, 
\begin{equation}\label{eq5.18.2}
\mbox{$u=P$ in a
neighborhood of $K$}\quad \Longleftrightarrow \quad \mbox{$w=P$ in a neighborhood of
$s^{-1}K$},
\end{equation}

\noindent where
$s^{-1}K:=\{x\in\RR^n:\,sx\in K\}$. Also,
\begin{equation}\label{eq5.19}
\int_{C_{a/2,2b}}|\nabla^m w(x)|^2\,dx= \int_{C_{a/2,2b}}|\nabla^m_x
[u(sx)]|^2\,dx=s^{2m-n} \int_{C_{as/2,2bs}}|\nabla^m_y u(y)|^2\,dy,
\end{equation}

\noindent so that
\begin{equation}\label{eq5.20}
s^{2m-n} {\rm Cap}_P(K, C_{as/2,2bs})= {\rm Cap}_P(s^{-1}K,
C_{a/2,2b}).
\end{equation}

\noindent However, $s^{-1} K\subset \overline{C_{a,b}}$ and
therefore by (\ref{eq5.4}) 
\begin{equation}\label{eq5.20.1}
{\rm Cap}_P(s^{-1}K,
C_{a/2,2b})\leq {\rm Cap}_P(\overline{C_{a,b}},
C_{a/2,2b}),
\end{equation}

\noindent uniformly in  $s$. The right-hand side of \eqref{eq5.20.1} is a constant independent of $s$ and therefore, \eqref{eq5.20}--\eqref{eq5.20.1} along with \eqref{eq5.12} 
prove the first part of \eqref{eq5.12.1}. The bound on ${\rm Cap}_P^\Phi$ then follows from Lemma~\ref{lcap1}. \ep

\begin{lemma}\label{l5.2-2}
Let $K$ be a compactum in $\overline{C_{as,bs}}$ for some
$s>0$, $0<a<b<\infty$. If $n\in [3, 2m+1]$ is odd, then for every $P\in\Pi_1$
\begin{equation}\label{eq5.12}
{\rm Cap}_P(K, \RR^n\setminus\{O\})\approx  {\rm Cap}_P(K,
C_{as/2,2bs}),\end{equation}
\noindent with the constants independent of $s$.

Furthermore, for any $n\in [2, 2m+1]$ we have 
\begin{multline}\label{eq5.12-2}
\inf\Bigg\{\sum_{k=0}^m\int_{\RR^n} \frac{|\nabla^ku(x)|^2}{|x|^{2m-2k}}\,dx:\,\,u\in \ring W^{m,2}(\RR^n\setminus\{O\}),\,\,u=P\mbox{ in a
neighborhood of }K\Bigg\}\\[4pt]\approx {\rm Cap}_P(K,
C_{as/2,2bs}).
\end{multline}
\end{lemma}

\noindent In fact, an analogue of \eqref{eq5.12} holds for ${\rm Cap}_P^\Phi$ as well, in the sense that one can harmlessly enlarge the set  $C_{as/2,2bs}$ in ${\rm Cap}_P^\Phi(K,
C_{as/2,2bs})$ if $K\subset \overline{C_{as,bs}}$. However, one would not take $\RR^n\setminus\{O\}$, since for the elements of the space $\ring W^{m,2}(\RR^n\setminus\{O\})$ the functional $\Phi$ need not be finite.

\bp The $``\leq"$ inequality in \eqref{eq5.12} follows from 
the monotonicity property
(\ref{eq5.5}). Similar monotonicity considerations treat the $``\leq"$ inequality in \eqref{eq5.12-2} once we observe that 
\begin{equation}\label{eq5.13.1}\sum_{k=0}^m\int_{C_{as/2,2bs}} \frac{|\nabla^ku(x)|^2}{|x|^{2m-2k}}\,dx \approx \int_{C_{as/2,2bs}} |\nabla^mu(x)|^2\,dx \quad \mbox{for}\quad u\in\ring W^{m,2}(C_{as/2,2bs})
\end{equation}
by Hardy inequality. 

Let us turn to the opposite direction, starting with \eqref{eq5.12}. 
To this end,  fix some $\eps>0$ and take $u\in \ring
W^{m,2}(\RR^n\setminus\{O\})$ such that $u=P$ in a neighborhood of
$K$ and
\begin{equation}\label{eq5.14}
\int_{\RR^n}|\nabla^m
u(x)|^2\,dx < {\rm Cap}_P(K, \RR^n\setminus\{O\})+\eps.
\end{equation}

\noindent Next, take  the cut-off function
$\zeta \in C_0^\infty(a/2,2b)$ such that 
$\zeta=1$ on the set $[3a/4,3b/2]$
and let 
\begin{equation}\label{eq5.14.0}
u_\zeta (x):=\zeta(|x|/s)u(x),\qquad x\in \RR^n.
\end{equation}

\noindent Then
\begin{equation}\label{eq5.16}
u_\zeta\in\ring W^{m,2}(C_{as/2,2bs})\quad\mbox{and}\quad u_\zeta=P\,\,\mbox{in a
neighborhood of }\,\, K,
\end{equation}

\noindent and therefore, 
\begin{equation}\label{eq5.17}
{\rm Cap}_P(K, C_{as/2,2bs})\leq \int_{C_{as/2,2bs}}|\nabla^m
u_\zeta(x)|^2\,dx.
\end{equation}

\noindent However,
\begin{eqnarray}\label{eq5.18}\nonumber
&&\int_{C_{as/2,2bs}}|\nabla^m
u_\zeta(x)|^2\,dx=\int_{C_{as/2,2bs}}\left|\nabla^m
\left(\zeta(|x|/s)u(x)\right)\right|^2\,dx\\[4pt]
&&\qquad \leq C\sum_{k=0}^m\frac{1}{s^{2m-2k}}\int_{C_{as/2,2bs}} |\nabla^ku(x)|^2\,dx\leq C\sum_{k=0}^m\int_{\RR^n} \frac{|\nabla^ku(x)|^2}{|x|^{2m-2k}}\,dx\nonumber\\[4pt]
&&\qquad \leq C \int_{\RR^n}|\nabla^m
u(x)|^2\,dx,
\end{eqnarray}

\noindent using Hardy inequality in odd dimensions (see, e.g., \cite{MazSobSpNew}, (1.3.3)) for the last estimate above. Combining this with \eqref{eq5.14}, we finish the proof of (\ref{eq5.12}).  

The same argument implies \eqref{eq5.12-2}. We will only be missing the last step of \eqref{eq5.18} as the Hardy inequality for even dimensions in our range does not generally apply on $\RR^n\setminus \{O\}$.  This is the reason why one has to modify the definition of capacity to include the full sum of $L^2$ norms of $k$-th gradients, $k=0, ..., m$ as in \eqref{eq5.12-2}.
\ep

\begin{remark}\label{rChOfAn} The same argument demonstrates that in any dimension $n\in [2, 2m+1]$, if $K$ is a compactum with $K\subseteq \overline{C_{as,bs}}$ and $K\subseteq \overline{C_{a's,b's}}$  for some
$s>0$, $0<a<b<\infty$, $0<a'<b'<\infty$, then 
\begin{equation}\label{eqChOfAn}
{\rm Cap}_P(K, C_{as/2,2bs})\approx  {\rm Cap}_P(K,
C_{a's/2,2b's}),\end{equation}
with the implicit constants depending on $m,n,a,b,a',b'$ only.
\end{remark}

\section{Sufficient condition for $\lambda$-regularity}\label{sSuff}
\setcounter{equation}{0}

\subsection{Poincar\'e-type inequalities}

We shall now list some auxiliary lemmas that will ultimately lead to the optimal choice of $P$ minimizing the capacity ${\rm Cap}_P$.

\begin{lemma}\label{l5.3.0}  Let $s>0$, $0<a<b<\infty$.  Then  
\begin{equation}\label{eq5.21.0}
\Phi[u-P; C_{as,bs}]=\Phi[u; C_{as,bs}],
\end{equation}

\noindent for every $P\in \Pi$, $u\in V^{m,2}(C_{as,bs})$.
\end{lemma}

\bp This is a direct consequence of the definitions. Indeed, if 
\begin{equation}\label{eq5.24.0}
P(x):= \sum_{p\in Z} \sum_{l=-p}^{p} b_{pl}
Y_l^p(x/|x|), \qquad  x\in\RR^n\setminus \{O\},
\end{equation}

\noindent then 
\begin{multline}\label{eq5.24.0.1}
\Phi[u-P; C_{as,bs}]\\[4pt] =\sum_{p\in Z} \sum_{l=-p}^{p} \sum_{k=0}^m c_{kp}\int_{\log\frac{1}{bs}} ^{\log\frac{1}{as}}
(\partial_t^k (v_{pl}(t)-b_{pl}))^2\,dt+\sum_{p\not\in Z} \sum_{l=-p}^{p} \sum_{k=0}^m c_{kp}\int_{\log\frac{1}{bs}} ^{\log\frac{1}{as}}
(\partial_t^k v_{pl}(t))^2\,dt\\[4pt]
=\sum_{p\in Z} \sum_{l=-p}^{p} \sum_{k=0}^m c_{kp}\int_{\log\frac{1}{bs}} ^{\log\frac{1}{as}} 
(\partial_t^k v_{pl}(t))^2\,dt+\sum_{p\not\in Z} \sum_{l=-p}^{p} \sum_{k=0}^m c_{kp}\int_{\log\frac{1}{bs}} ^{\log\frac{1}{as}}
(\partial_t^k v_{pl}(t))^2\,dt\\[4pt]
=\Phi[u; C_{as,bs}],
\end{multline}

\noindent where in the second inequality we used the fact that $c_{0p}=0$ for every $p\in Z$, and if $k\geq 1$ then $\partial_t^k$ annihilates the constants $b_{pl}$. \ep

\begin{lemma}\label{l5.3} Let $s>0$, $0<a<b<\infty$. Then for every $u\in V^{m,2}(C_{as,bs})$ there
exists ${\mathcal P}={\mathcal P}(u,s,a,b)\in \Pi$ with the property 
\begin{equation}\label{eq5.21}
\|u-{\mathcal P}\|_{L^2(C_{as,bs})}^2\leq C s^{n} \Phi[u; C_{as,bs}].
\end{equation}
\end{lemma}

\bp Recall the definition of $\Phi$ in \eqref{cap6.0}. Given $u\in
L^2(C_{as,bs})$ and $v=u\circ \varkappa^{-1}$, we decompose $v$ as in \eqref{eqo2.7}. Using Poincar{\'e}'s inequality, we choose for every $p\in Z$, $-p\leq l\leq p$, the constants $b^0_{pl}$ (depending on $u$, $s$, $a$, $b$) such that 
\begin{equation}\label{eq5.23}
\int_{\log\frac{1}{bs}}^{\log\frac{1}{as}}|v_{pl}(t)-b^0_{pl}|^2\,dt \leq C 
\int_{\log\frac{1}{bs}}^{\log\frac{1}{as}}|\partial_tv_{pl}(t)|^2\,dt,
\end{equation}

\noindent and set 
\begin{equation}\label{eq5.24}
{\mathcal P}(x):= \sum_{p\in Z} \sum_{l=-p}^{p} b^0_{pl}
Y_l^p(x/|x|), \qquad  x\in\RR^n\setminus \{O\}.
\end{equation}

\noindent Then  
\begin{eqnarray}\label{eq5.24.1}
\|u-{\mathcal P}\|_{L^2(C_{as,bs})}^2 &= &\sum_{p\in Z} \sum_{l=-p}^{p} \int_{\log\frac{1}{bs}}^{\log\frac{1}{as}}|v_{pl}(t)-b^0_{pl}|^2e^{-tn}\,dt +\sum_{p\not\in Z} \sum_{l=-p}^{p} \int_{\log\frac{1}{bs}}^{\log\frac{1}{as}}|v_{pl}(t)|^2e^{-tn}\,dt \nonumber\\[4pt] &\leq & C s^n\sum_{p\in Z} \sum_{l=-p}^{p} \int_{\log\frac{1}{bs}}^{\log\frac{1}{as}}|\partial_tv_{pl}(t)|^2\,dt +Cs^n\sum_{p\not\in Z} \sum_{l=-p}^{p} \int_{\log\frac{1}{bs}}^{\log\frac{1}{as}}|v_{pl}(t)|^2\,dt \nonumber\\[4pt] &\leq &Cs^n\Phi[u; C_{as,bs}],
\end{eqnarray}

\noindent where the last inequality follows from \eqref{capdef3} and \eqref{capdef2}. More precisely, we use the fact that $c_{1p}>0$ for $p\in Z$ and $c_{0p}>0$ for $p\not\in Z$. This finishes the proof of the Lemma.
\ep

\begin{lemma}\label{lExt} Let $s>0$, $0<a<b<\infty$. There exists an extension operator
\begin{equation}\label{eq5.28}
{\rm Ex}: V^{m,2}(C_{as,bs})\to V^{m,2}(C_{as/2,2bs}),
\end{equation}

\noindent with the operator norm independent of $s$ and satisfying the
properties

\begin{enumerate}

\item  ${\rm Ex}\, u(x) =u(x)$ for all $x\in C_{as,bs}$,

\item ${\rm Ex}\, P=P$  for any $P\in\Pi$,

\item if  ${\rm dist}\, ({\rm supp}\,u,K)>0$ for some compactum $K\subset \overline{C_{as,bs}}$
then ${\rm dist}\, ({\rm supp}\,{\rm Ex}\,u,K)>0$.
\end{enumerate}
\end{lemma}

\noindent {\it Remark}. The condition (iii) above can be reformulated as follows. If $u=0$ in a neighborhood of a compactum $K\subset \overline{C_{as,bs}}$
then ${\rm Ex}\, u$
vanishes in some neighborhood of $K$ contained in $C_{as/2,2bs}$.

\vskip 0.08 in

\bp Using the spherical harmonics decomposition, the matters can be reduced to the construction of a suitable one-dimensional reflection-type operator. Following the ``reflection of the finite order" procedure (see, e.g., \cite{MazyaSobolevSpaces}, \S 1.1.16, p.31), 
for every $m\in\NN$, one can define an operator
\begin{equation}\label{eqExt1}
{\rm Ex}_0: C^\infty([a,b])\mapsto C^m([a/2, 2b]),
\end{equation}

\noindent which extends to a bounded operator on Sobolev spaces, so that, in particular, 
\begin{equation}\label{eqExt2}
\|\partial_r^k\, {\rm Ex}_0\,u\|_{L^2((a/2, 2b))}\leq C \|\partial_r^ku\|_{L^2((a, b))}, \quad \mbox{for}\quad k=0,1,...,m,
\end{equation}

\noindent with the properties
\begin{enumerate}

\item  ${\rm Ex}_0\, u(r) =u(r)$ for all $r\in (a,b)$,

\item ${\rm Ex}_0\, c=c$  for any constant $c$, and moreover,  if for some  $\eps\in\left(0,\frac{b}{a}-1\right)$
\begin{equation}\label{eqExt3}
u(r)=c \quad\mbox{for} \quad r\in \left[a,a(1+\eps)\right] \quad \mbox{then}\quad  
{\rm Ex}_0 u(r)=c\quad\mbox{for} \quad r\in \left[a\left(1-\frac{\eps a}{2(b-a)}\right),a\right]
\end{equation}

\noindent and symmetrically, if for some  $\eps\in\left(0,\frac{b}{a}-1\right)$
\begin{equation}\label{eqExt4}
u(r)=c \quad\mbox{for} \quad r\in \left[b(1-\eps),b\right] \quad \mbox{then}\quad  
{\rm Ex}_0 u(r)=c\quad\mbox{for} \quad r\in \left[b, b\left(1+\frac{\eps a}{(b-a)}\right)\right].
\end{equation}
\end{enumerate}

One can now make a change of coordinates $t=\log \frac {1}{sr}$ to obtain the operator ${\rm Ex}_1$ extending the functions from $\left[\log \frac{1}{bs}, \log\frac{1}{as}\right]$ to $\left[\log \frac{1}{2bs}, \log\frac{2}{as}\right]$ and satisfying an analogue of the set of the properties above. Note that, in particular, 
\begin{equation}\label{eqExt5}
\|\partial_t^k{\rm Ex}_1v\|_{L^2\left(\left(\log \frac{1}{bs}, \log\frac{1}{as}\right)\right)}\leq C \|\partial_t^kv\|_{L^2\left(\left(\log \frac{1}{2bs}, \log\frac{2}{as}\right)\right)}, \quad \mbox{for}\quad k=0,1,...,m,
\end{equation}

\noindent with the norm independent of $s$.

Finally, let 
\begin{equation}\label{eqExt6}
{\rm Ex}_2 v(t,\omega)=\sum_{p=0}^\infty \sum_{l=-p}^{p} {\rm Ex}_1v_{pl}(t)\,
Y_l^p(\omega), \qquad t\in\RR, \,\,\omega\in S^{n-1}, 
\end{equation}

\noindent and ${\rm Ex} \,u:={\rm Ex}_2 v$ for $v=u\circ \varkappa^{-1}$. Now one can directly check that  \eqref{eqExt6}, \eqref{eqExt5} imply \eqref{eq5.28}, and (i)--(ii) for the operator ${\rm Ex}_0$ ultimately give the properties (i)--(iii) for the operator ${\rm Ex}$.

The only non-trivial conclusion is the property (iii). To see that it holds, take any small $\eps>0$ and denote by $U_\eps(K)$ the set $\{x\in\RR^n:\,{\rm dist}\,(x,K)<\eps\}$. Suppose $u=0$ in $U_\eps(K)\cap \overline{C_{as,bs}}$ for some $\eps>0$. We claim that ${\rm Ex}\, u=0$ in $U_{\eps/4}(K)$. Indeed, ${\rm Ex} \,u(x)=0$ for $x\in U_{\eps/4}(K)\cap {C_{as,bs}}$ by definition and property (i). If, on the other hand, $x\in U_{\eps/4}(K)$ does not belong to ${C_{as,bs}}$, then we consider a segment on a ray $x/|x|=const$ with one end $x_1=x$ and the other $x_2\in \overline{C_{as,bs}}$, of the length $|x_1-x_2|=\eps/2$. The entire segment is at the distance less than $\eps$ from $K$. Therefore, $u(z)=0$ for all $z$ lying on the intersection of the segment $ [x_1,x_2]$ with $\overline{C_{as,bs}}$. Hence, by the properties \eqref{eqExt3}--\eqref{eqExt4} we have $u(z)=0$ for  $z$ on the entire segment $[x_1,x_2]$, in particular, $u(x)=0$, as desired. This finishes the proof of the Lemma.
\ep

%\vskip 0.08in  \noindent{\it Remark.} The relation (\ref{eqbc13})
%is a Poincar{\'e}-type inequality. The quasi-polynomial ${\mathcal P}\in\Pi$
%plays a role of approximating function and the dependence of it on
%the original function $u$ and the domain $C_{s,as}$ can not be
%relaxed in this argument.

%GIVE REFERENCE ON SPHERICAL HARMONICS

\begin{lemma}\label{l5.4.0} Let $s>0$, $0<a<b<\infty$ and $\zeta \in C_0^\infty(a/2,2b)$ be such that $\zeta=1$ on the set $[3a/4,3b/2]$. Then 
\begin{equation}\label{eq5.25.0}
\Phi[\zeta(|\cdot |/s)u; C_{as/2,2bs}]\leq C \|u\|_{V^{m,2}(C_{as/2,2bs})},
\end{equation}

\noindent for every $u\in V^{m,2}(C_{as/2,2bs})$.
\end{lemma}

\bp Note that $\zeta(|\cdot |/s)u\in \ring W^{m,2}(C_{as/2,2bs})$. Hence, by \eqref{cap6.0.3},
\begin{eqnarray}\label{eq5.25.1}
\Phi[\zeta(|\cdot |/s)u; C_{as/2,2bs}]
&\leq & C  \sum_{\stackrel{k,i\geq 0}{ i+k\leq
m}}\int_{\RR}\int_{S^{n-1}} \left(\partial_t^k\nabla_\omega^i
\left(\zeta(e^{-t}/s) v(t,\omega)\right)\right)^2 \,d\omega dt
\nonumber\\[4pt]
&\leq & C  \sum_{\stackrel{k,i\geq 0}{ i+k\leq
m}}\int_{\log\frac{1}{2bs}}^{\log\frac{2}{as}}\int_{S^{n-1}} \left(\partial_t^k\nabla_\omega^i
 v(t,\omega)\right)^2 \,d\omega dt,
\end{eqnarray}

\noindent where $v=u\circ\varkappa^{-1}$. The latter expression is, in turn,  bounded by the right-hand side of \eqref{eq5.25.0}, by \eqref{cap6.0.2}, as desired. \ep

\begin{proposition}\label{l5.4} Suppose $s>0$, $0<a<b<\infty$, and $K$ is a
compact subset of $\overline{C_{as,bs}}$ such that ${\rm
Cap}\,(K,C_{as/2,2bs})>0$. Then for every $u\in V^{m,2}(C_{as, bs})$ with ${\rm dist}\, ({\rm supp}\,u,K)>0$ the following estimate holds
\begin{equation}\label{eq5.25}
\frac{1}{s^{2n-2m}}\int_{C_{as,bs}}|u(x)|^2\,dx\leq \frac{C}{{\rm
Cap}\,(K,C_{as/2,2bs})} \,\,\Phi[u; C_{as,bs}],
\end{equation}

\noindent with the constant independent of $s$.
\end{proposition}

\bp Within this particular argument it is convenient to take a different norm in the space $\Pi$, namely,
$\|P\|_{\Pi}:=\|P\|_{L^2(C_{a,b})}$ and 
$\Pi_1=\{P\in\Pi:\,\|P\|_{\Pi}=1\}$ with such a norm.  This is an
equivalent norm in $\Pi$ and hence, the capacities defined by \eqref{cap7}, \eqref{cap8}--\eqref{cap9} with the new normalization for $P$ are equivalent to the original ones. Therefore, it is enough to prove \eqref{eq5.25} assuming the $P$'s implicitly present in ${\rm Cap}$ are such that $\|P\|_{L^2(C_{a,b})}=1$.

Let us now turn to (\ref{eq5.25}).
Without loss of generality we may assume that
$\|u\|_{L^2(C_{as,bs})}=s^{n/2}$. Then \eqref{eq5.25} reduces to
\begin{equation}\label{eq5.32}
s^{2m-n} {\rm Cap}(K, C_{as/2,2bs})\leq C\Phi[u;
C_{as,bs}],
\end{equation}

\noindent for $u$ and (implicit) $P$ normalized as above. 

Recall the result of Lemma~\ref{l5.3}. Let us denote by ${\mathcal P}={\mathcal P}(u,s,a)$ the function in $\Pi$ satisfying
(\ref{eq5.21}), and by $C_0$ the constant $C$ in
(\ref{eq5.21}). We would like to split the discussion according to whether $\Phi[u; C_{as,bs}]\geq 1/(4C_0)$ or $\Phi[u; C_{as,bs}]\leq 1/(4C_0)$. In the first case, one employs Lemma~\ref{l5.2}, \eqref{eq5.12.1}, and immediately gets the desired estimate \eqref{eq5.32}. 

As for the other situation, 
the first step is to show that for every $P\in\Pi$ with $\|P\|_{L^2(C_{a,b})}=1$ and every $u\in V^{m,2}(C_{as, bs})$
\begin{equation}\label{eq5.26}
{\rm Cap}_P^\Phi(K, C_{as/2,2bs})\leq C \|P-u\|^2_{V^{m,2}(C_{as, bs})}.
\end{equation}

To this end, take  the  function
$\zeta \in C_0^\infty(a/2,2b)$ such that 
$\zeta=1$ on the set $[3a/4,3b/2]$
and let $w(x):=\zeta(|x|/s)(P(x)-{\rm Ex}\,
u(x))$, $x\in C_{as/2,2bs}$, where ${\rm Ex}$ is the extension operator from Lemma~\ref{lExt}.
 Then first of all, 
$w\in\ring W^{m,2}(C_{as/2,2bs})$ and secondly, by our assumptions and property (iii) of  ${\rm Ex}$ the function ${\rm Ex}\,u$ vanishes in some neighborhood of $K$, so that $w=P$ in some neighborhood of $K$. Hence,  
\begin{equation}\label{eq5.26.1}
{\rm Cap}_P^\Phi (K, C_{as/2,2bs})\leq \Phi[w; C_{as/2,2bs}]\leq C\|P-{\rm Ex}\, u\|_{V^{m,2}(C_{as/2,2bs})},
\end{equation}

\noindent by Lemma~\ref{l5.4.0}. Furthermore, by Lemma~\ref{lExt}
\begin{equation}\label{eq5.26.2}
\|P-{\rm Ex}\, u\|_{V^{m,2}(C_{as/2,2bs})}=\|{\rm Ex}\, (P-u)\|_{V^{m,2}(C_{as/2,2bs})}\leq C\|P-u\|_{V^{m,2}(C_{as,bs})},
\end{equation}

\noindent which yields (\ref{eq5.26}).

Next, using Lemma~\ref{lcap1} one can see that the inequality \eqref{eq5.32} (which we aim to prove) is equivalent to 
\begin{equation}\label{eq5.32.1}
\inf_{P\in\Pi_1} {\rm Cap}_P^\Phi (K, C_{as/2,2bs})\leq C\Phi[u;
C_{as,bs}].
\end{equation}

\noindent Hence, by \eqref{eq5.26} it is enough to show 
\begin{equation}\label{eq5.32.3}
\inf_{P\in\Pi_1} \|P-u\|^2_{V^{m,2}(C_{as, bs})}= \inf_{P\in\Pi_1} \left(s^{-n}\|P-u\|^2_{
L^2(C_{as,bs})}+\Phi[P-u; C_{as,bs}]\right) \leq C\Phi[u;
C_{as,bs}].
\end{equation}

\noindent However, $\Phi[P-u; C_{as,bs}]=\Phi[u; C_{as,bs}]$ by Lemma~\ref{l5.3.0}. Hence,  the estimate above is, in fact, just a bound on $\inf_{P\in\Pi_1} s^{-n}\|P-u\|^2_{
L^2(C_{as,bs})}$. More precisely, we are left to prove that
\begin{equation}\label{eq5.32.4}
\inf_{P\in\Pi_1} \|P-u\|^2_{
L^2(C_{as,bs})} \leq Cs^n\Phi[u;
C_{as,bs}], 
\end{equation}

\noindent for every $u$ such that $\|u\|_{L^2(C_{as,bs})}=s^{n/2}$.

Lemma~\ref{l5.3} and the normalization assumptions on $u$ imply that 
\begin{equation}\label{eq5.34}
 2\|u-{\mathcal P}\|_{L^2(C_{as,bs})}\leq 2 \left( C_0s^n\Phi[u; C_{as,bs}]\right)^{\frac 12} \leq 2 \left(C_0s^{n}/(4C_0)\right)^{\frac 12}=s^{n/2}=\|u\|_{L^2(C_{as,bs})}.
\end{equation}

\noindent Therefore, 
\begin{equation}\label{eq5.34.1}
\|{\mathcal P}\|_{L^2(C_{as,bs})}\leq \|u-{\mathcal P}\|_{L^2(C_{as,bs})}+\|u\|_{L^2(C_{as,bs})}\leq \|u\|_{L^2(C_{as,bs})}/2+\|u\|_{L^2(C_{as,bs})}= 3s^{n/2}/2,
\end{equation}

\noindent and conversely, 
\begin{equation}\label{eq5.34.2}
\|{\mathcal P}\|_{L^2(C_{as,bs})}\geq \|u\|_{L^2(C_{as,bs})}-\|u-{\mathcal P}\|_{L^2(C_{as,bs})}\geq \|u\|_{L^2(C_{as,bs})}-\|u\|_{L^2(C_{as,bs})}/2=s^{n/2}/2.
\end{equation}

\noindent Now we renormalize ${\mathcal P}$ to get an element of $\Pi_1$. To do this, take
\begin{equation}\label{eq5.36}
P:=\frac{{\mathcal P}}{\|{\mathcal
P}\|_{L^2(C_{a,b})}}=s^{n/2}\,\frac{{\mathcal P}}{\|{\mathcal
P}\|_{L^2(C_{as,bs})}}.
\end{equation}

\noindent Clearly, $P\in\Pi_1$ by definition. Furthermore, 
\begin{eqnarray}\label{eq5.36.1}
&&\|P-{\mathcal P}\|_{L^2(C_{as,bs})} = \left\|\frac{s^{n/2}}{\|{\mathcal
P}\|_{L^2(C_{as,bs})}}
\,{\mathcal P}-{\mathcal P}\right\|_{L^2(C_{as,bs})} = \left|\frac{s^{n/2}}{\|{\mathcal
P}\|_{L^2(C_{as,bs})}}-1\right|\, \|{\mathcal P}\|_{L^2(C_{as,bs})} 
\nonumber\\[4pt]&&\qquad =\left|s^{n/2} - \|{\mathcal P}\|_{L^2(C_{as,bs})}
\right|=\left| \|u\|_{L^2(C_{as,bs})} - \|{\mathcal P}\|_{L^2(C_{as,bs})}
\right| \leq \| u - {\mathcal P} \|_{L^2(C_{as,bs})}.
\end{eqnarray}

\noindent This implies that 
\begin{equation}\label{eq5.36.2}
\| u - P \|_{L^2(C_{as,bs})} \leq  \| u - {\mathcal P} \|_{L^2(C_{as,bs})}+\|P-{\mathcal P}\|_{L^2(C_{as,bs})} 
\leq 2\| u - {\mathcal P} \|_{L^2(C_{as,bs})}.
\end{equation}

\noindent Finally, in concert with the first inequality in \eqref{eq5.34} the estimate \eqref{eq5.36.2} yields
\begin{equation}\label{eq5.37}
\|u-P\|_{L^2(C_{as,bs})}^2\leq 4 \|u-{\mathcal P}\|_{L^2(C_{as,bs})}^2
\leq 4C_0s^n \Phi[u; C_{as,bs}].
\end{equation}

\noindent The latter estimate, indeed, confirms \eqref{eq5.32.4} and finishes the argument.
 \ep

\subsection{Odd dimensions}

%In the following two sections we investigate the necessary and
%sufficient conditions for the 1-regularity of a point $Q\in\po$.
%As before, by applying a simple shift of coordinates, we can
%assume that $O \not\in \Omega$, and thus we can take $Q=O$. This
%will be convenient since the origin plays a special role in our
%definition of the biharmonic capacity.

\begin{proposition}\label{p7.1}
Let $\Omega$ be a bounded domain in $\RR^n$, $3\leq n \leq 2m+1$,  $n$ odd, 
$O\in\RR^n\setminus\Omega$, $R>0$ and
\begin{equation}\label{eq3.13}
(-\Delta)^m u=f \,\,{\mbox{in}}\,\,\Omega, \quad f\in
C_0^{\infty}(\Omega\setminus B_{4R}),\quad u\in \ring
W^{m,2}(\Omega).
\end{equation}

\noindent Fix some $b\geq 4$. Then for every $x\in
B_{R/b^4}\cap\Omega$
\begin{eqnarray}\label{eq7.2}\nonumber
\sum_{i=0}^\lambda \frac{|\nabla^{i}u(x)|^2}{{|x|^{2\lambda-2i}} }&\leq & \frac{C}{R^{n+2\lambda}}\int_{C_{R,4R}\cap\Omega}
|u(y)|^2\,dy\,\,\\[4pt]&&\qquad
\times\exp \left( -c\sum_{j=2}^{l}(Rb^{-2j})^{(2m-n)} {\rm
Cap}\,(\overline{C_{R \,b^{-2j},R\,b^{-2(j-1)}}}\setminus\Omega,C_{R\,b^{-2j}/2,2R \,b^{-2(j-1)}})\right),
\end{eqnarray}

\noindent where $l\geq 2$, $l\in\NN$, is such that $|x|\leq b^{-2l}R$ and $\lambda$ is defined by \eqref{eqo7.5}.

In particular,  when $O$ is a boundary point of $\Omega\subset \RR^n$ 
\begin{equation}\label{eq7.3}
\mbox{if } \sum_{j=1}^{\infty}a^{-j(2m-n)} {\rm
Cap}\,(\overline{C_{a^{-j},a^{-(j-1)}}}\setminus\Omega,C_{a^{-j}/2,2a^{-(j-1)}})=+\infty\,\,\, \mbox{then $O$ is $\lambda$-regular}.
\end{equation}
Here $a$ is any real number greater than 1.
\end{proposition}

\noindent {\it Proof of Proposition~\ref{p7.1}.} To begin, we approximate $\Omega$ by a sequence of domains with smooth boundaries
$\left\{\Omega_n\right\}_{n=1}^\infty$ satisfying
\begin{equation}\label{eqo3.5}
\bigcup_{n=1}^\infty \Omega_n=\Omega\quad\mbox{and}\quad
{\overline{\Omega}}_n\subset\Omega_{n+1} \quad \mbox{for every}
\quad n\in\NN.
\end{equation}

\noindent Choose $n_0\in\NN$ such that ${\rm supp}\,f\subset
\Omega_n$ for every $n\geq n_0$ and denote by $u_n$ a unique
solution of the Dirichlet problem
\begin{equation}\label{eqo3.6}
(-\Delta)^m u_n = f\quad {\rm in} \quad \Omega_n, \quad  u_n\in \ring
W^{m,2}(\Omega_n), \quad n\geq n_0.
\end{equation}

\noindent The sequence  $\{u_n\}_{n=n_0}^\infty$ converges to $u$ in
$\ring W^{m,2}(\Omega)$ (see, e.g., \cite{Necas}, \S{6.6}). 

 Furthermore, for every $s\leq R$ we introduce a cut-off $\eta_s\in C_0^\infty(B_{2s})$ such that
\begin{equation}\label{eqo3.7.1}
0\leq\eta_s\leq 1\,\, {\rm in}\,\,B_{2s}, \quad \eta_s=1\,\, {\rm
in}\,\,B_s,\quad {\rm and} \quad |\nabla^k \eta_s|\leq Cs^{-k}, \quad
k\leq 2m.
\end{equation}

\noindent One can use the property that $u_n$ is polyharmonic in $\Omega_n\cap B_{4R}\supset \Omega_n\cap B_{4s}$ and \eqref{cap5.1} to deduce that 
\begin{equation}\label{eqn7.1}
\Phi\left[\frac{\eta_s u_n}{|x|^{\lambda}}; \Omega\right]\leq \frac{C}{s^{n+2\lambda}}\int_{C_{s,4s}}
|u_n(y)|^2\,dy.
\end{equation}

Indeed, since $u_n\in \ring
W^{m,2}(\Omega_n)$, $\Omega_n$ is bounded and ${\rm dist}\,(O,\Omega_n)>0$, we have according to \eqref{cap5.1}
\begin{equation}\label{cap5.1.1}
\Phi\left[\frac{\eta_s u_n}{|x|^{\lambda}}; \Omega\right]=\int_{\RR^n}(-\Delta)^m \left(\eta_s u_n\right)\,\left(|x|^{2m-n-2\lambda}\eta_s u_n\right)\,dx.
\end{equation}
Since $u_n$ is polyharmonic in $\Omega_n\cap B_{4R}$ and $\eta_s$ is
supported in $B_{2R}$, one can see that $\eta_s\,(-\Delta)^m u_n =0$, and hence, the expression above is equal to 
\begin{equation}\label{eqo3.9}
\int_{\RR^n}\bigg([(-\Delta)^m,\eta_s] u_n(x)\bigg) \bigg(\eta_s(x)u_n(x)\,|x|^{2m-n-2\lambda}\bigg)\,dx,
\end{equation}

\noindent where the brackets denote the commutator, i.e., $$[(-\Delta)^m,\eta_s] u_n(x)=(-\Delta)^m\left(\eta_s(x)u_n(x)\right)-\eta_s(x)(-\Delta)^m u_n(x),$$ the integral in (\ref{eqo3.9}) is understood in the sense of pairing between $\ring W^{m,2}(\Omega_n)$ and its dual. 
Evidently, the
support of the integrand is a subset of ${\rm
supp}\,\nabla\eta\subset C_{s,2s}$, and therefore, the expression
in (\ref{eqo3.9}) is bounded by
\begin{equation}\label{eqo3.10}
C\sum_{i=0}^m \frac{1}{s^{2\lambda+n-2i}}\int_{C_{s,2s}} |\nabla^i
u_n(x)|^2\,dx\leq \frac{C}{s^{2\lambda+n}}\int_{C_{s,4s}} |u_n(x)|^2\,dx,
\end{equation}

\noindent using Lemma~\ref{lo3.1}. Hence, we arrive at \eqref{eqn7.1}.

Furthermore, it follows directly from definitions that 
\begin{equation}\label{eqn7.2}
\Phi\left[\frac{u_n}{|x|^{\lambda}}; B_s\right] \leq \Phi\left[\frac{\eta_s u_n}{|x|^{\lambda}}; \Omega\right],
\end{equation}

\noindent where $\Phi\left[\frac{u_n}{|x|^{\lambda}}; B_s\right]$ is understood in the sense of $\eqref{cap6.0}$ with $b=s$ and $a$ smaller than ${\rm dist}\,(O,\Omega_n)>0$.

Next, denote
\begin{equation}\label{eq7.3.5}
\varphi(s):=\sup_{|x|\leq s}\left(\sum_{i=0}^\lambda \frac{|\nabla^{i}u_n(x)|^2}{{|x|^{2\lambda-2i}} }\right)+\Phi\left[\frac{u_n}{|x|^{\lambda}}; B_s\right],\qquad s>0.
\end{equation}

\noindent Then 
\begin{equation}\label{eq7.3.6}
\varphi(s)\leq \frac{C}{s^{n+2\lambda}}\int_{C_{s,16s}}
|u_n(y)|^2\,dy, \qquad s\leq R/4.
\end{equation}

\noindent using \eqref{eqn7.1}, \eqref{eqn7.2} and Proposition~\ref{po7.2}. 

Now fix $b\geq 4$.  
If  ${{\rm
Cap}\,(\overline{C_{s,b^2s}}\setminus\Omega,\RR^n\setminus
\{O\})}>0$ then according to Proposition~\ref{l5.4} and \eqref{eq7.3.6}
\begin{eqnarray}\label{eq7.3.6.1}
\varphi(s)&\leq & C\,\frac{s^{n-2m}}{s^{2n-2m}}\int_{C_{s,b^2s}}
\frac{|u_n(y)|^2}{|y|^{2\lambda}}\,dy \leq C\,\frac{s^{n-2m}}{{{\rm
Cap}\,(\overline{C_{s,b^2s}}\setminus\Omega,C_{s/2,2b^2s})}}\,\Phi\left[\frac{u_n}{|x|^{\lambda}}; C_{s,b^2s}\right]\nonumber\\[4pt]
&\leq & C\,\frac{s^{n-2m}}{{{\rm
Cap}\,(\overline{C_{s,b^2s}}\setminus\Omega,C_{s/2,2b^2s})}}\,\left(\varphi(b^2s)-\varphi(s)\right), \quad\mbox{for any}\quad s\leq R/4.
\end{eqnarray}

\noindent Note that the condition ${\overline{\Omega}_n}\subset\Omega$ guarantees that the distance from ${\rm supp}\,(u_n/|x|^{\lambda})$ to $\overline{C_{s,b^2s}}\setminus\Omega$ is strictly greater than zero and justifies the use of Proposition~\ref{l5.4}. One can see that there exists $c>0$ such that
\begin{eqnarray}\label{eq7.8}
\varphi(s)&\leq &\frac{1}{1+C^{-1}\, s^{2m-n}{{\rm
Cap}\,(\overline{C_{s,b^2s}}\setminus\Omega,C_{s/2,2b^2s})}}\,\varphi(b^2s)\nonumber\\[4pt]
&\leq & \exp \left( -cs^{2m-n} {\rm
Cap}\,(\overline{C_{s,b^2s}}\setminus\Omega,C_{s/2,2b^2s})\right)\,\varphi(b^2s),
\end{eqnarray}

\noindent since 
\begin{equation}\label{eq7.8.1}
s^{2m-n} {{\rm
Cap}\,(\overline{C_{s,b^2s}}\setminus\Omega,C_{s/2,2b^2s})}\leq C,
\end{equation}

\noindent by \eqref{eq5.12}, \eqref{eq5.12.1}. Now we iterate the process taking $s=R\,b^{-2(j+1)}$, $j\in\NN$. Then
\begin{equation}\label{eq7.9.0}
\varphi\left(R \,b^{-2(j+1)}\right)\leq
\exp \left( - cR^{2m-n}b^{-2(j+1)(2m-n)} {\rm
Cap}\,(\overline{C_{R \,b^{-2(j+1)},R \,b^{-2j}}}\setminus\Omega,C_{R \,b^{-2(j+1)}/2,2R \,b^{-2j}})\right)\,\varphi\left(R\,b^{-2j}\right),
\end{equation}

\noindent and hence, 

\begin{equation}\label{eq7.9}
\varphi\left(R\,b^{-2l}\right)\leq
\exp \left( -cR^{2m-n}\sum_{j=2}^{l}b^{-2j(2m-n)} {\rm
Cap}\,(\overline{C_{R \,b^{-2j},R\,b^{-2(j-1)}}}\setminus\Omega,C_{R\,b^{-2j}/2,2R \,b^{-2(j-1)}})\right)\,\varphi\left(Rb^{-2}\right),
\end{equation}

\noindent for all for $l=2,3,...$.

Pick $l=2,3,...$ so that
\begin{equation}\label{eq7.10}
b^{-2l-2}R\leq |x|\leq b^{-2l}R.
\end{equation}
Then \eqref{eq7.9} and \eqref{eq7.3.6} imply that for $x$ as in \eqref{eq7.10} (or, in fact, any $x$ such that $|x|\leq b^{-2l}R$), we have
\begin{equation}\label{eq7.12.0.2} \sum_{i=0}^\lambda \frac{|\nabla^{i}u_n(x)|^2}{{|x|^{2\lambda-2i}} }
\leq
\exp \left( -cR^{2m-n}\sum_{j=2}^{l}b^{-2j(2m-n)} {\rm
Cap}\,(\overline{C_{R \,b^{-2j},R\,b^{-2(j-1)}}}\setminus\Omega,C_{R\,b^{-2j}/2,2R \,b^{-2(j-1)}})\right)\,\varphi\left(Rb^{-2}\right).
\end{equation}

\noindent  Moreover, we note that 
\begin{equation}\label{eq7.3.6.0}
\varphi\left(b^{-2}R\right)\leq \frac{C}{R^{n+2\lambda}}\int_{C_{b^{-2}R,16b^{-2}R}}
|u_n(y)|^2\,dy \leq \frac{C}{R^{n+2\lambda}}\int_{C_{R,4R}}
|u_n(y)|^2\,dy,
\end{equation}

\noindent by Proposition~\ref{po7.2}.
The combination of \eqref{eq7.12.0.2} with \eqref{eq7.3.6.0} then yields \eqref{eq7.2} with $u_n$ in place of $u$, and the limiting procedure finishes the argument for \eqref{eq7.2}. 

Now let us turn to \eqref{eq7.3}. The estimate \eqref{eq7.2} directly leads to the following conclusion. When $O$ is a boundary point of $\Omega\subset \RR^n$ 
\begin{equation}\label{eq7.3-aux}
\mbox{if } \sum_{j=1}^{\infty}(a^{-j}R)^{(2m-n)} {\rm
Cap}\,(\overline{C_{a^{-j}R,a^{-(j-1)R}}}\setminus\Omega,C_{a^{-j}R/2,2a^{-(j-1)}R})=+\infty\,\,\, \mbox{then $O$ is $\lambda$-regular},
\end{equation}
where $a=b^2\geq 16$. Next, the condition $a=b^2\geq 16$ can be substituted by any $a>1$ using monotonicity of capacity to shrink $\overline{C_{R \,b^{-2j},R\,b^{-2(j-1)}}}\setminus\Omega$ as necessary (starting with some $b\geq 4$ such that $b^2/a\in\NN$),  and then Remark~\ref{rChOfAn} to adjust the ambient annulus. The exact constant in intervening inequalities would depend on the ratio of $b^2\geq 16$ and $a>1$, but that does not affect the final result \eqref{eq7.3-aux}. Finally, there exists $N\in\ZZ$ such that $R\approx a^{-N}$, so that the series in \eqref{eq7.3-aux} can be rewritten as the series in \eqref{eq7.3}, with the summation over $j=N+1, N+2, ...$, but that again does not affect the question of convergence. Hence, we arrive at \eqref{eq7.3}.
 \ep

The results of Proposition~\ref{p7.1} can be turned into the estimates on polyharmonic functions at infinity, respectively, still being restricted to the case of the odd dimension. 

\begin{proposition}\label{p7.2}Let $\Omega$ be a bounded domain in $\RR^n$, $3\leq n \leq 2m+1$, $n$ odd, 
$O\in\RR^n\setminus\Omega$, $r>0$ and assume that
\begin{equation}\label{eq7.15}
(-\Delta)^m u=f \,\,{\mbox{in}}\,\,\Omega, \quad f\in
C_0^{\infty}(B_{r/4}\cap\Omega),\quad u\in \ring W^{m,2}(\Omega).
\end{equation}

\noindent Fix some $b\geq 4$. Then  for any $x\in\Omega\setminus
B_{b^4r}$ 
\begin{eqnarray}\label{eq7.16}
\nonumber \sum_{i=0}^{\lambda} |\nabla^i u(x)|^2 |x|^{2\lambda+2n-4m+2i}&  \leq &
 C\,r^{n+2\lambda-4m}\int_{C_{\frac r4, r}}
|u_n(y)|^2\,dy\,\,\\[4pt]&&\qquad
\times\exp \left( -c\sum_{j=2}^{l}(rb^{2j})^{(2m-n)} {\rm
Cap}\,(\overline{C_{rb^{2(j-1)},rb^{2j}}}\setminus\Omega_n,C_{rb^{2(j-1)}/2,2 rb^{2j}})\right),
\end{eqnarray}

\noindent where $l\geq 2$, $l\in\NN$, is such that $|x|\geq b^{2l}r$ and 
$\lambda$ given by \eqref{eqo7.5}. \end{proposition}

\bp Retain the
approximation of $\Omega$ with the sequence of smooth domains
$\Omega_n$ satisfying (\ref{eqo3.5}) and define $u_n$ according to
(\ref{eqo3.6}). We denote by ${\mathcal I}$ the inversion $x\mapsto
y=x/|x|^2$ and by $U_n$ the Kelvin transform of $u_n$,
\begin{equation}\label{eqo3.26}
U_n(y):=|y|^{2m-n}\, u_n(y/|y|^2), \quad y\in {\mathcal I}(\Omega_n).
\end{equation}

\noindent Then
\begin{equation}\label{eqo3.27}
(-\Delta)^m U_n(y)=|y|^{-n-2m} ((-\Delta)^m u_n)(y/|y|^2),
\end{equation}

\noindent and therefore, $U_n$ is polyharmonic in ${\mathcal
I}(\Omega_n)\cap B_{4/r}$. Moreover, 
\begin{equation}\label{eqo3.29}
U_n\in \ring W^{m,2}({\mathcal I}(\Omega_n))
\quad\Longleftrightarrow\quad u_n\in \ring W^{m,2}(\Omega_n).
\end{equation}

\noindent Observe also that $\Omega_n$ is a bounded domain with
$O\not\in \overline{\Omega}_n$, hence, so is ${\mathcal I}(\Omega_n)$
and $O\not\in \overline{{\mathcal I}(\Omega_n)}$.

For any $x\in \Omega\setminus B_{4r}$ 
\begin{equation}\label{eqo3.32}
|\nabla^i u_n(x)|\leq C \sum_{k=0}^i |x|^{2m-n-i-k}\,(\nabla^k
U_n)(x/|x|^2),
\end{equation}

\noindent since $u_n(x)=|x|^{2m-n}\, U_n(x/|x|^2)$. Hence,
\begin{equation}\label{eq3.32.0}
\sum_{i=0}^\lambda |\nabla^i u_n(x)||x|^{\lambda+n-2m+i} \leq C \sum_{i=0}^\lambda |x|^{\lambda-i}\,(\nabla^i
U_n)(x/|x|^2),
\end{equation}

\noindent where $U_n$ comes from the Kelvin transform of $u_n$ and falls under the scope of Proposition~\ref{p7.1} with $R=1/r$ and ${\mathcal I}(\Omega_n)$ in place of $\Omega$. It follows that
for any $b\geq 4$, any $l\geq 2$, $l\in\NN$, and any $x$ such that $|x|\geq b^{2l}r$ we have 
\begin{eqnarray}\label{eq3.32.0.1}
&& \sum_{i=0}^\lambda \left(\frac{(\nabla^i
U_n)(x/|x|^2)}{(x/|x|^2)^{\lambda-i}}\right)^2  \leq Cr^{n+2\lambda}\int_{C_{\frac 1r,\frac 4r}}
|U_n(y)|^2\,dy\,\,\\[4pt]&&\qquad
\times\exp \left( -c\sum_{j=2}^{l}(b^{-2j}/r)^{(2m-n)} {\rm
Cap}\,(\overline{C_{b^{-2j}/r,b^{-2(j-1)}/r}}\setminus{\mathcal I}(\Omega_n),C_{b^{-2j}/(2r),2 b^{-2(j-1)}/r})\right).\nonumber
\end{eqnarray}

\noindent Next, we would like to express the capacity of the set ${\overline{C_{\frac 1{b^2s},\frac 1s}}}\setminus{\mathcal
I}(\Omega_n)$  in terms of the capacity of the set $\overline{C_{s,b^2s}}\setminus\Omega_n$, with teh goal of using this for $s=r/b^{-2(j-1)}$, as above.

Fix any $\eps>0$ and choose $u\in \ring W^{m,2}(C_{\frac 1{2b^2s},\frac 2s})$ with $u=P\in \Pi_1$ in a neighborhood of
${\overline{C_{\frac 1{b^2s},\frac 1s}}}\setminus{\mathcal I}(\Omega_n)$ so that 

\begin{equation}\label{eq7.18.1}
\int_{C_{\frac 1{2b^2s},\frac 2s}}|\nabla^mu(x)|^2\,dx <{\rm
Cap}\,(\overline{C_{\frac 1{b^2s},\frac 1s}}\setminus{\mathcal
I}(\Omega_n),C_{\frac 1{2b^2s},\frac 2s})+\eps.
\end{equation}

\noindent Then the function given by 
 $U(y):=u(y/|y|^2)$, $y\in \RR^n\setminus \{O\}$,  belongs to  $\ring W^{m,2}(C_{\frac s2,2b^2s})$, and   $U(y)=P(y/|y|^2)=P(y)$ for all $y$ in a
neighborhood of ${\overline{C_{ s, b^2s}}}\setminus\Omega_n$.
Moreover, 
\begin{equation}\label{eq7.19}
\int_{C_{\frac 1{2b^2s},\frac 2s}} |\nabla^m u(x)|^2\,
dx\approx s^{4m-2n} \int_{C_{\frac s{2},2b^2s}} |\nabla^m U(y)|^2\, dy,
\end{equation}

\noindent analogously to \eqref{eqo3.32}. Therefore, 
\begin{equation}\label{eqn7.17}
 s^{4m-2n} \,{\rm
Cap}\,(\overline{C_{s,b^2s}}\setminus\Omega_n,C_{\frac{s}{2},2b^2s})\leq C\, {\rm
Cap}\,(\overline{C_{\frac 1{b^2s},\frac 1s}}\setminus{\mathcal I}(\Omega_n),C_{\frac 1{2b^2s},\frac {2}{s}}).
\end{equation}

One can see that the
opposite inequality also holds, by reduction  to the previous case with $1/(b^2s)$ in
place of $s$ and ${\mathcal I}(\Omega_n)$ in place of $\Omega_n$. Therefore, 
\begin{equation}\label{eqn7.17.1}
 s^{4m-2n} \,{\rm
Cap}\,(\overline{C_{s,b^2s}}\setminus\Omega_n,C_{s/2,2b^2s})\approx {\rm
Cap}\,(\overline{C_{\frac 1{b^2s},\frac 1s}}\setminus{\mathcal I}(\Omega_n),C_{\frac 1{2b^2s},\frac 2s}),
\end{equation}

\noindent with the implicit constant independent of $s$.

Finally, the combination of \eqref{eq3.32.0} with \eqref{eqn7.17.1} gives for any $l\geq 2$, $l\in\NN$, and any $x$ such that $|x|\geq b^{2l}r$
\begin{eqnarray}\label{eq7.20}\nonumber
&&\hskip -1cm \sum_{i=0}^\lambda |\nabla^i u_n(x)|^2|x|^{2\lambda+2n-4m+2i} \leq Cr^{n+2\lambda}\int_{C_{\frac 1r,\frac 4r}}
|U_n(y)|^2\,dy\,\,\\[4pt]&&\qquad
\times\exp \left( -c\sum_{j=2}^{l}(b^{-2j}/r)^{(2m-n)} {\rm
Cap}\,(\overline{C_{b^{-2j}/r,b^{-2(j-1)}/r}}\setminus{\mathcal I}(\Omega_n),C_{b^{-2j}/(2r),2 b^{-2(j-1)}/r})\right)
\nonumber\\[4pt]\nonumber
&& \leq Cr^{n+2\lambda}\int_{C_{\frac 1r,\frac 4r}}
|U_n(y)|^2\,dy\,\,\\[4pt]&&\qquad
\times\exp \left( -c\sum_{j=2}^{l}(rb^{2j})^{(2m-n)} {\rm
Cap}\,(\overline{C_{rb^{2(j-1)},rb^{2j}}}\setminus\Omega_n,C_{rb^{2(j-1)}/2,2 rb^{2j}})\right)
\nonumber\\[4pt]\nonumber
 &&\hskip -1cm \qquad \leq C\,r^{n+2\lambda-4m}\int_{C_{\frac r4, r}}
|u_n(y)|^2\,dy\,\,\\[4pt]&&\qquad
\times\exp \left( -c\sum_{j=2}^{l}(rb^{2j})^{(2m-n)} {\rm
Cap}\,(\overline{C_{rb^{2(j-1)},rb^{2j}}}\setminus\Omega_n,C_{rb^{2(j-1)}/2,2 rb^{2j}})\right)\end{eqnarray}
\noindent since $\overline{\Omega}_n\subset\Omega$. Now the
argument can be finished using the limiting procedure. \ep

\subsection{Even dimensions}
The results for even dimensions, while yielding sufficiency of the capacitory condition in our Wiener test, do not quite take the form of \eqref{eq7.2} or \eqref{eq7.16}. This is unfortunate but it is in line with the existing literature on the subject pertaining to the Wiener-type conditions for continuity of solutions (see, e.g., \cite{M2}). 

\begin{proposition}\label{p7.1-even}
Let $\Omega$ be a bounded domain in $\RR^n$, $2\leq n \leq 2m$,  $n$ even, 
$O\in\RR^n\setminus\Omega$, $R>0$ and
\begin{equation}\label{eq3.13-even}
(-\Delta)^m u=f \,\,{\mbox{in}}\,\,\Omega, \quad f\in
C_0^{\infty}(\Omega\setminus B_{4R}),\quad u\in \ring
W^{m,2}(\Omega).
\end{equation}

\noindent If
\begin{equation}\label{eq7.3-even}
 \sum_{j=1}^{\infty}j\,a^{-j(2m-n)} {\rm
Cap}\,(\overline{C_{a^{-j},a^{-(j-1)}}}\setminus\Omega,C_{a^{-j}/2,2a^{-(j-1)}})=+\infty,
\end{equation}
for some $a>1$, then 
\begin{equation}\label{eq7.3.1-even}
\lim_{s\to 0} \sup_{|x|\leq s} \sum_{i=0}^\lambda \frac{|\nabla^{i}u(x)|^2}{{|x|^{2\lambda-2i}} }=0,
\end{equation}
\noindent where $\lambda$ is defined by \eqref{eqo7.5}.
Hence, the condition \eqref{eqo1.9-even} implies that $O$ is $\lambda$-regular.
\end{proposition}

\bp We start off the same way as in the proof of Proposition~\ref{p7.1} with $\Phi_s$ in place of $\Phi$, $s\leq R$, (recall the definition in \eqref{cap6.0-log}--\eqref{cap6.0-log1}), and follow \eqref{eqo3.5}--\eqref{eqo3.10}, in particular, establishing
\begin{equation}\label{eqn7.2-even}
\Phi_s\left[\frac{u_n}{|x|^{\lambda}}; B_s\right] \leq \Phi_s\left[\frac{\eta_s u_n}{|x|^{\lambda}}; \Omega\right]\leq \frac{C}{s^{2\lambda+n}}\int_{C_{s,4s}} |u_n(x)|^2\,dx.
\end{equation}

\noindent Indeed, as before, the definition of, e.g., $\Phi_s\left[\frac{u_n}{|x|^{\lambda}}; B_s\right]$ is 
well-justified for $B_s$ in place of an annulus since the support of $u_n$ is away from the origin, we use \eqref{cap5.1-log} in place of \eqref{cap5.1} to obtain analogue of \eqref{cap5.1.1} with a logarithmic weight, and the remaining argument is the same upon observing that due to the fact that $\eta_s (-\Delta)^m u_n=0$ the integration is restricted to the annulus $C_{s, 2s}$, that is, $\ln (2s)^{-1}\leq t \leq \ln s^{-1}$ and hence, $\ln 2 \leq \ln 4s +t\leq \ln 4$. Then for any $b\geq 4$
\begin{multline*} 
\frac{1}{R^{2\lambda+n}}\int_{C_{R,4R}} |u_n(x)|^2\,dx \gtrsim\Phi_R\left[\frac{u_n}{|x|^{\lambda}}; B_R\right]=\sum_{j=0}^\infty \Phi_R\left[\frac{u_n}{|x|^{\lambda}}; C_{b^{-j-2}R, b^{-j}R}\right] \gtrsim \sum_{j=0}^\infty j \,\Phi \left[\frac{u_n}{|x|^{\lambda}}; C_{b^{-j-2}R, b^{-j}R}\right] \\[4pt]
\gtrsim \sum_{j=0}^\infty j \,{{\rm
Cap}\,(\overline{C_{b^{-j-2}R, b^{-j}R}}\setminus\Omega,C_{b^{-j-2}R/2, 2b^{-j}R})} \,\frac{1}{(b^{-j}R)^{2\lambda+n}}\int_{C_{b^{-j-2}R, b^{-j}R}} |u_n(x)|^2\,dx\\[4pt]
\gtrsim \sum_{j=0}^\infty j \,{{\rm
Cap}\,(\overline{C_{b^{-j-2}R, b^{-j}R}}\setminus\Omega,C_{b^{-j-2}R/2, 2b^{-j}R})} \sup_{|x|\leq b^{-j-2}R} \sum_{i=0}^\lambda \frac{|\nabla^{i}u(x)|^2}{{|x|^{2\lambda-2i}} }
\end{multline*}
where we used Proposition~\ref{l5.4} in the third inequality and Proposition~\ref{po7.2} for the fourth one. Hence, if the limit in \eqref{eq7.3.1-even} is strictly positive, the series 
$$\sum_{j=0}^\infty j \,{{\rm
Cap}\,(\overline{C_{b^{-j-2}R, b^{-j}R}}\setminus\Omega,C_{b^{-j-2}R/2, 2b^{-j}R})} <\infty, $$
which yields the statement of Proposition~\ref{p7.1-even} with a slightly modified version of  \eqref{eq7.3-even}. And, similarly to the end of the proof of Proposition~\ref{p7.1} we can adjust the annuli as convenient to obtain the desired result with any $a>1$.\ep

\section{Necessary condition for $\lambda$-regularity}
\setcounter{equation}{0}

This section will be entirely devoted to the proof of the second
part of Theorem~\ref{to1.2}, i.e. the necessary condition for
1-regularity. 
\subsection{Fine estimates on the quadratic forms, general discussion}
To start, let us set the notation and discuss some fine estimates on the involved quadratic forms. 

{\it Assume first that $n\in [3,2m+1]$ is odd.}

Let us recall the result of Theorem~\ref{to4.2} and the functions $g$ and $h$ from \eqref{eqo4.12}. We denote by $\Q_{g,\tau}(v)$ the quadratic form associated to the right-hand side of \eqref{eqo4.11}, that is, 
\begin{multline}\label{eqn1}
\Q_{g,\tau}(v)=\sum_{\stackrel{i\geq 0, \,k\geq 0}{0\leq i+k\leq m}}\int_{\RR}\int_{S^{n-1}} \A_{ik}(\partial_t) [e^{-t}g(t,\tau)] \left(\partial_t^k\nabla_\omega^iv\right)^2 \,d\omega dt\\[4pt]
=c_0 \int_{S^{n-1}} v^2(\tau, \omega) \,d\omega+ \sum_{\stackrel{i\geq 0, \,k\geq 0}{0<i+k\leq m}}\int_{\RR}\int_{S^{n-1}} \A_{ik}(\partial_t) [e^{-t}g(t,\tau)](t,\tau) \left(\partial_t^k\nabla_\omega^iv\right)^2 \,d\omega dt,
\end{multline}
where $\A_{ik}(\partial_t)$ are appropriate operators given by polynomials of $\partial_t$ and $c_0>0$.  The form is uniquely defined by the requirement that when $u$ and $v$ are as in the statement of Theorem~\ref{to4.2}, we have 
\begin{equation}\label{eqn2}
\int_{\RR^n}(-\Delta)^m u(x)\,u(x)g(\log |x|^{-1},\log |\xi|^{-1})\,dx=\Q_{g,\tau}(v),
\end{equation}
\noindent but clearly one can in principle consider $\Q_{g,\tau}(v)$ for more general functions $v$.  An interested reader can look at (4.28) in \cite{MayMazInvent2}, {\it loc.cit.}, substituting $h(t-\tau)$ by $e^t g(t,\tau)$, to read off the exact representations of $\A_{ik}(\partial_t)$. Also, according to Theorem~\ref{to4.2} (or rather its proof) we have 
\begin{multline}\label{eqn2.1}
\Q_{g,\tau}(v)  \gtrsim \int_{S^{n-1}}v^2(\tau,\omega)\,d\omega \\[4pt] + \sum_{\stackrel{k\geq 1, \,i\geq 0}{i+k\leq m}} \int_{\RR}\int_{S^{n-1}}
\left(\partial_t^k \nabla_\omega^i v\right)^2\,d\omega dt+
\int_{\RR}\int_{S^{n-1}} v\,\prod_{p=-\frac n2+\frac 32}^{m-\frac
n2+\frac 12}\left(-\delta-p\,(p+n-2)\right) v\,d\omega dt,
\end{multline}
\noindent and the inequality holds for any $v$ such that $\Q_{g,\tau}(v)$ is finite.
We remark that the second equality in \eqref{eqn1}, i.e., the fact that the only term with $v^2$ in the quadratic form reduces to $c_0 \int_{S^{n-1}} v^2(\tau, \omega) \,d\omega$ for some $c_0>0$, can also be extracted from the proof of Theorem~\ref{to4.2} -- this consideration basically dictated the choice of $h$ and, respectively, of $g$. 

Analogously, we shall denote by $\Q_{e^th,\tau}(v)$ the corresponding quadratic form with the function $C_1e^t h(t-\tau)$ in place of $g(t,\tau)$ (cf. \eqref{eqo4.12}). When $u$ and $v$ are as in the statement of Theorem~\ref{to4.2}, 
\begin{equation}\label{eqn3.2}
C \int_{\RR^n}(-\Delta)^m u(x)\,u(x)|x|^{-1}h(\log |x|^{-1},\log |\xi|^{-1})\,dx= \Q_{e^th,\tau}(v),
\end{equation}
for some $C>0$.

We do not know if the resulting quadratic form $\Q_{e^th,\tau}(v)$ is positive, however, once again, examination of the proof of Theorem~\ref{to4.2} reveals that  one can write 
\begin{multline}\label{eqn3}
\Q_{e^th,\tau}(v)=\sum_{\stackrel{i\geq 0, \,k\geq 0}{0\leq i+k\leq m}}\int_{\RR}\int_{S^{n-1}} \A_{ik}(\partial_t) h(t-\tau) \left(\partial_t^k\nabla_\omega^iv\right)^2 \,d\omega dt\\[4pt]
=c_0 \int_{S^{n-1}} v^2(\tau, \omega) \,d\omega+ \sum_{\stackrel{i\geq 0, \,k\geq 0}{0<i+k\leq m}}\int_{\RR}\int_{S^{n-1}} \A_{ik}(\partial_t) h(t-\tau) \left(\partial_t^k\nabla_\omega^iv\right)^2 \,d\omega dt.
\end{multline}
and, even more precisely, that 
\begin{multline}\label{eqn3.1}
\Q_{e^th,\tau}(v)=c_0 \int_{S^{n-1}} v^2(\tau, \omega) \,d\omega+ \sum_{\stackrel{i\geq 0, \,k\geq 1}{i+k\leq m}}\int_{\RR}\int_{S^{n-1}} \A_{ik}(\partial_t) h(t-\tau) \left(\partial_t^k\nabla_\omega^iv\right)^2 \,d\omega dt\\[4pt]
+ c_0 \int_{\RR}\int_{S^{n-1}} {\mathcal
L}^{m,n}(0,\delta)v(t,\omega)\,v(t,\omega)h(t-\tau)\,d\omega
dt\\[4pt]
 +\,c_0\sum_{p=0}^\infty \sum_{l=-p}^{p} \int_{\RR}
v_{pl}^2(t) \left({\mathcal L}^{m,n}(-\partial_t,-p(p+n-2))-{\mathcal L}^{m,n}(-\partial_t,0)\right) h(t-\tau)\, dt,
\end{multline}
\noindent (see (4.28) in \cite{MayMazInvent2}, {\it loc.cit.}). Here, 
\begin{equation}\label{eqo4.26}
\left({\mathcal L}^{m,n}(-\partial_t,-p(p+n-2))-{\mathcal L}^{m,n}(-\partial_t,0)\right) h(t)\geq 0, \quad \mbox{ for all } \quad t\neq 0,\quad 0\leq p\leq  m-\frac n2+\frac 12,
\end{equation}
\noindent (see (4.26) in \cite{MayMazInvent2}, {\it loc.cit.}), and for the remaining part 
\begin{multline}\label{eqn3.1.1}
 \int_{\RR}\int_{S^{n-1}} {\mathcal
L}^{m,n}(0,\delta)v(t,\omega)\,v(t,\omega)|h(t-\tau)|\,d\omega
dt\\[4pt]
+\sum_{p=m-\frac n2 +\frac 32}^\infty \sum_{l=-p}^{p} \int_{\RR}
v_{pl}^2(t) \left|\bigg({\mathcal L}^{m,n}(-\partial_t,-p(p+n-2))-{\mathcal L}^{m,n}(-\partial_t,0)\bigg) h(t-\tau)\right|\, dt \\[4pt] \lesssim  \int_{\RR}\int_{S^{n-1}} {\mathcal
L}^{m,n}(0,\delta)v(t,\omega)\,v(t,\omega)\,d\omega
dt,
\end{multline}
(see (4.30) and (4.33), (4.41)--(4.42) in \cite{MayMazInvent2}, {\it loc.cit.}). In fact, these bounds together with Theorem~\ref{to2.1} dictate the choice of $g$ and give \eqref{eqn2.1} and Theorem~\ref{to4.2} (see \cite{MayMazInvent2}). In particular, we also have, by construction, 
\begin{multline}\label{eqn2.2}
\Q_{g,\tau}(v) \gtrsim \int_{S^{n-1}}v^2(\tau,\omega)\,d\omega + \sum_{\stackrel{k\geq 1, \,i\geq 0}{i+k\leq m}} \int_{\RR}\int_{S^{n-1}}
\left(\partial_t^k \nabla_\omega^i v\right)^2\,d\omega dt\\[4pt] +
\int_{\RR}\int_{S^{n-1}} v\,\prod_{p=-\frac n2+\frac 32}^{m-\frac
n2+\frac 12}\left(-\delta-p\,(p+n-2)\right) v\,d\omega dt\\[4pt]
+ \sum_{p=0}^{m-\frac n2+\frac 12} \sum_{l=-p}^{p} \int_{\RR}
v_{pl}^2(t) \left({\mathcal L}^{m,n}(-\partial_t,-p(p+n-2))-{\mathcal L}^{m,n}(-\partial_t,0)\right) h(t-\tau)\, dt.
\end{multline}
\noindent This will be useful later on. 

\vskip 0.08 in

{\it Suppose now that $n\in [2,2m]$ is even. }

We take $g=g_R$ as defined in \eqref{eqo6.15} and let $\Q_{g,\tau}(v)$ be the quadratic form associated to the right-hand side of \eqref{eqo6.14}, that is, 
\begin{multline}\label{eqn1.1}
\Q_{g,\tau}(v)=\Q_{g,\tau}^R(v)=\sum_{\stackrel{i\geq 0, \,k\geq 0}{0\leq i+k\leq m}}\int_{\RR}\int_{S^{n-1}} \A_{ik}^o(\partial_t) g_R(t,\tau) \left(\partial_t^k\nabla_\omega^iv\right)^2 \,d\omega dt\\[4pt]
=c_0 \int_{S^{n-1}} v^2(\tau, \omega) \,d\omega+ \sum_{\stackrel{i\geq 0, \,k\geq 0}{0<i+k\leq m}}\int_{\RR}\int_{S^{n-1}} \A_{ik}^o(\partial_t) g_R(t,\tau)(t,\tau) \left(\partial_t^k\nabla_\omega^iv\right)^2 \,d\omega dt,
\end{multline}
where $\A_{ik}^o(\partial_t)$ are appropriate operators given by polynomials of $\partial_t$ and $c_0>0$. Just as in the case of odd dimensions, the form is uniquely defined by \eqref{eqn2} (but note that $v=e^{\lambda t}(u\circ \varkappa^{-1})$, with $\lambda$ different in the even and odd case according to \eqref{eqo7.5}, $\xi\in B_{2R}$,  $\tau =\log |\xi|^{-1}$, and, of course,  $g$ is different). 
The parameter $R$ is only relevant for even dimensions and will be underlined when important in the proof.

Once again, according to (the proof of) Theorem~\ref{to6.3} we have 
\begin{multline}\label{eqn2.1.0}
\Q_{g,\tau}(v)  \gtrsim \int_{S^{n-1}}v^2(\tau,\omega)\,d\omega + \sum_{k=1}^m \sum _{i=0}^{m-k}\int_{\RR}\int_{S^{n-1}}
\left(\partial_t^k \nabla_\omega^i v\right)^2\,(C+C_R+t)\,d\omega dt
\\[4pt]
 +\,C
\int_{\RR}\int_{S^{n-1}} v\,{\mathcal L}^{m,n}_o(0,\delta) v\,(C+C_R+t)\,d\omega dt,\end{multline}
that is, whenever  $m$ is even,
\begin{multline}\label{eqn2.1.1}
\Q_{g,\tau}(v)  \gtrsim \int_{S^{n-1}}v^2(\tau,\omega)\,d\omega + \sum_{k=1}^m \sum _{i=0}^{m-k}\int_{\RR}\int_{S^{n-1}}
\left(\partial_t^k \nabla_\omega^i v\right)^2\,(C+C_R+t)\,d\omega dt
\\[4pt]
 +\,C
\int_{\RR}\int_{S^{n-1}} v\,\prod_{p}\left(-\delta-p\,(p+n-2)\right)^2 v\,(C+C_R+t)\,d\omega dt,\end{multline}
\noindent where the product  is over $p=-n/2+2, -n/2+4,..., m-n/2-2, m-n/2$, that is, $p=-n/2+2j$ with $j=1,2,...,m/2$. If $m$ is odd,
\begin{multline}\label{eqn2.1.2}
\Q_{g,\tau}(v)  \gtrsim \int_{S^{n-1}}v^2(\tau,\omega)\,d\omega + \sum_{k=1}^m \sum _{i=0}^{m-k}\int_{\RR}\int_{S^{n-1}}
\left(\partial_t^k \nabla_\omega^i v\right)^2\,(C+C_R+t)\,d\omega dt
\\[4pt]
+\,C
\int_{\RR}\int_{S^{n-1}} v\,\prod_{p}\left(-\delta-p\,(p+n-2)\right)^2 \left(-\delta_{\omega}+(n/2-1)^2\right) v\,(C+C_R+t)\,d\omega dt,
\end{multline}
\noindent where the product is over $p=-n/2+3, -n/2+5,..., m-n/2-2, m-n/2$, that is, $p=-n/2+1+2j$ with $j=1,2,...,(m-1)/2$. In both cases, $C>0$ is a positive constant depending on $m$ and $n$ only, $C_R=\log(4R)$ depends on the support of $u$ and will be chosen appropriately below, $\xi\in B_{2R}$,  $\tau =\log |\xi|^{-1}$, and 
the inequality holds for any $v$ such that $\Q_{g,\tau}(v)$ is finite. 

Similarly to \eqref{eqn3.2}--\eqref{eqn2.2}, some finer lower estimates are available for $\Q_{g,\tau}(v)$ if we carefully examine the proof of Theorem~\ref{to6.3} in \cite{MayMazInvent2}. Specifically, with $h$ as in Lemmas~\ref{lo6.1} and \ref{lo6.2}, we have 
\begin{multline}\label{eqo6.28}\quad {\mathcal L}^{m,n}_o(-\partial_t,-p(p+n-2))h(t)\geq 0, \quad t\neq 0,  \\[4pt] \quad \mbox{for all}\,\,0\leq p\leq m-n/2 \,\,\mbox{even when $m-n/2$ is even},
\\[4pt] \quad \mbox{and for all}\,\,0\leq p\leq m-n/2 \,\,\mbox{odd when $m-n/2$ is odd}.\end{multline}
This corresponds to (6.28) and (6.38) in \cite{MayMazInvent2}, {\it loc.cit.} As discussed right after (6.28) and (6.38) in \cite{MayMazInvent2}, {\it loc.cit.} , the same is true for the function $\widetilde g$ in place of $h$, 
\begin{equation}\label{eqo6.16}
\widetilde g (t,\tau):=h(t-\tau)+\mu^{(4)}(C_R+\tau).\qquad t,\tau\in\RR.
\end{equation}
since for the relevant values of indices ${\mathcal L}^{m,n}_o(-\partial_t,-p(p+n-2))$, viewed as a polynomial in $\partial_t$, has a double root at zero for all $0\leq p\leq m-n/2$ even and hence, the result of its action (as an operator) on $\widetilde g$ is the same as the result of its action on $h$. We record for the future reference: 
\begin{multline}\label{eqo6.28-g}\quad {\mathcal L}^{m,n}_o(-\partial_t,-p(p+n-2))\widetilde g (t,\tau)\geq 0, \quad t\neq \tau,  \\[4pt] \quad \mbox{for all}\,\,0\leq p\leq m-n/2 \,\,\mbox{even when $m-n/2$ is even},
\\[4pt] \quad \mbox{and for all}\,\,0\leq p\leq m-n/2 \,\,\mbox{odd when $m-n/2$ is odd}.\end{multline}
The refined version of the lower estimates on $\Q_{g,\tau}(v)$ can respectively be written as 
\begin{multline}\label{eqn2.1.0-refined}
\Q_{g,\tau}(v)  \gtrsim \int_{S^{n-1}}v^2(\tau,\omega)\,d\omega + \sum_{k=1}^m \sum _{i=0}^{m-k}\int_{\RR}\int_{S^{n-1}}
\left(\partial_t^k \nabla_\omega^i v\right)^2\,(C+C_R+t)\,d\omega dt
\\[4pt]
 +\,C
\int_{\RR}\int_{S^{n-1}} v\,{\mathcal L}^{m,n}_o(0,\delta) v\,(C+C_R+t)\,d\omega dt
\\[4pt] + \sum_p \sum_{l=-p}^p \int_{\RR}
v_{pl}^2(t) \,{\mathcal L}^{m,n}_o(-\partial_t,-p(p+n-2)) \widetilde g (t,\tau)\, dt, \end{multline}
where the last sum in $p$ runs over all $0\leq p\leq m-n/2$ even when $m-n/2$ is even, and all $0\leq p\leq m-n/2$ odd when $m-n/2$ is odd.

\subsection{The scheme of the proof}

With these preliminaries at hand, we now start the core of the argument for the necessity of the capacitory condition, the estimates on the polyharmonic potential. Assume first that the dimension is odd. Suppose that for some $P\in\Pi_1$  the series
in (\ref{eqo1.10}) is convergent. Then  for every $\eps>0$ there exists $N\in\NN$ such that
\begin{equation}\label{eq8.3}
\sum_{j=N-2}^\infty 2^{-j(2m-n)} \,{{\rm
Cap}_P\,(\overline{C_{2^{-j},2^{-j+2}}}\setminus\Omega, C_{2^{-j-2},2^{-j+4}})}<\eps.
\end{equation}

By the same token, whenever the dimension is even, we fix $P\in\Pi_1$ such that the series
in (\ref{eqo1.10-even}) is convergent and conclude that for every $\eps>0$ there exists $N\in\NN$ such that
\begin{equation}\label{eq8.3.even}
\sum_{j=N-2}^\infty 2^{-j(2m-n)} \,j\,{{\rm
Cap}_P\,(\overline{C_{2^{-j},2^{-j+2}}}\setminus\Omega, C_{2^{-j-2},2^{-j+4}})}<\eps.
\end{equation}

Now let $K:=\overline{B_{2^{-N}}}\setminus \Omega$ and
$D:=B_{10}\setminus K$. We shall prove that the point $O$ is not
$\lambda$-regular with respect to $D$, and therefore with respect to
$\Omega$, since $D$ coincides with $\Omega$ in a fixed
neighborhood of $O$ (see Proposition~\ref{p6.1}).

Roughly speaking, the counterexample will be furnished by the function $V$ defined as follows. Let us fix $P\in \Pi_1$ from \eqref{eq8.3} or \eqref{eq8.3.even} in the case of the odd and even dimension, respectively, and let $\PP(x):=|x|^\lambda P(x)$, $x\in\rn$. Take some cut-off $\eta\in C_0^\infty(B_{2})$
equal to 1 on $B_{3/2}$ and denote 
\begin{equation}\label{eq8.4.0}
f:=-[(-\Delta)^m,\eta] \PP. 
\end{equation} 
Let $V$ be a solution
of the boundary value problem
\begin{equation}\label{eq8.4}
(-\Delta)^m V=f \,\,{\mbox{in}}\,\,D, \quad V\in \ring W^{m,2}(D).
\end{equation}

\noindent  We shall prove the necessity part of Theorem~\ref{to1.2} by showing that $|\nabla^\lambda V|$ does not vanish as $x\to O$, $x\in D$. 

The function $V$ is basically built from a suitable potential, that is, $V=U-\eta \PP$, where $U$ is a solution to 
\begin{equation}\label{eq8.4.4}
(-\Delta)^m U=0 \,\,{\mbox{in}}\,\,D, \quad U\in \ring W^{m,2}(B_{10}), \quad U= \PP \,\,\mbox{on }K.
\end{equation}
Roughly speaking, we prove that the capacitory conditions imply that $U$ together with its derivatives is small in a suitable sense near $O$, and hence, $\nabla^\lambda V$ is large. 

At this moment, \eqref{eq8.4.0}--\eqref{eq8.4.4} is only a sketch of the idea. Detailed justifications of all the steps will be provided in the course of the proof below. Quite unexpectedly, the exact structure of $\PP$, peculiar to the parity of dimension and of $m$, $m-n/2$, becomes vital for the argument. Let us turn to the details.

\subsection{Main estimates. Bounds for auxiliary functions $T$ and $W$ related to polyharmonic potentials on the annuli}

The actual argument will not directly address $U$ and $V$ but rather their approximations. To this end, consider the covering
of $K=\overline{B_{2^{-N}}}\setminus \Omega$ by the sets $K\cap
C_{2^{-j},2^{-j+2}}$, $j\geq N$, and observe that
\begin{eqnarray}\label{eq8.25}
&& K\cap {\overline{C_{2^{-j},2^{-j+2}}}}=
{\overline{C_{2^{-j},2^{-j+2}}}}\setminus \Omega, \quad j\geq
N+2,\\[4pt]\label{eq8.26}
&&  K\cap {\overline{C_{2^{-j},2^{-j+2}}}}\subseteq
{\overline{C_{2^{-j},2^{-j+2}}}}\setminus \Omega,\quad j=N, N+1.
\end{eqnarray}

\noindent Let $\{\eta^j\}_{j=N-2}^\infty$ be the corresponding
partition of unity such that
\begin{equation}\label{eq8.27}
\eta^j\in C_0^\infty(C_{2^{-j},2^{-j+2}}),\quad |\nabla^k
\eta^j|\leq C2^{kj}, \quad k=0,1,2,\quad\mbox{ and }\quad
\sum_{j=N-2}^\infty \eta^j=1 \mbox{ on } {\overline{B_{2^{-N+1}}}}.
\end{equation}

\noindent By $U^j$ we denote the capacitary potential of $K\cap
\overline{C_{2^{-j},2^{-j+2}}}$ with the boundary data $P$, i.e.
the minimizer for the optimization problem
\begin{eqnarray}\label{eq8.28}\nonumber
&&\inf\Bigg\{\int_{C_{2^{-j-2},2^{-j+4}}}(\nabla^m
u(x))^2\,dx:\,\,u\in \ring W^{m,2}(C_{2^{-j-2},2^{-j+4}}),\\[4pt]
&&\qquad\qquad\qquad\qquad\,u=P\mbox{ in a neighborhood of }K\cap
{\overline{C_{2^{-j},2^{-j+2}}}}\Bigg\}.
\end{eqnarray}

\noindent Such $U^j$ always exists and belongs to $\ring
W^{m,2}(C_{2^{-j-2},2^{-j+4}})$ since $P$ is an infinitely
differentiable function in a neighborhood of $K\cap
{\overline{C_{2^{-j},2^{-j+2}}}}$. The infimum above is equal to
\begin{equation}\label{eq8.29}
{\rm Cap}_P \{K\cap
{\overline{C_{2^{-j},2^{-j+2}}}},C_{2^{-j-2},2^{-j+4}}\}.
\end{equation}

Let us now define the function
\begin{equation}\label{eq8.30}
T(x):=\sum_{j=N-2}^\infty |x|^\lambda\, \eta^j(x) U^j(x),\quad x\in \RR^n,
\end{equation}
\noindent and let $\vartheta:=e^{\lambda t}(T\circ \varkappa^{-1})$. We claim that 
\begin{equation}\label{eqn10}
\Q_{g,\tau} (\vartheta)<C\eps. 
\end{equation}
To be more precise, the following statement holds.

\begin{lemma}\label{lSt4} Given $P\in \Pi_1$ and $\PP(x):=|x|^\lambda P(x)$, $x\in\rn$, assume that the corresponding capacity satisfies \eqref{eq8.3}, \eqref{eq8.3.even}, in the case of the odd and the even dimension, respectively, and retain the definition of $T$, $\vartheta$ as above.  Fix some point $\xi\in B_{2^{-N}}$ and $\tau:=\log
|\xi|^{-1}$. Then 
\begin{equation}\label{eqn7}
\Q_{g,\tau} (\vartheta)<C\eps,
\end{equation}
is valid with a constant $C$ depending on $m$ and $n$ only. In the case when the dimension is even, we assume that the parameter $R=10$ in the definition of $\Q_{g,\tau}$ (see \eqref{eqn1.1}).
\end{lemma}
\bp To start, let us record the estimates on the weight functions (they can be read off the definitions or see \cite{MayMazInvent2}, (7.11)--(7.13), {\it loc. cit.}, for a more detailed discussion). We have 

We observe that when $n$ is odd, the formula \eqref{eqo4.12} yields 
\begin{equation}\label{eqo7.8.1}
\left|\nabla_x^k\, g(\log|x|^{-1},\log \rho^{-1})\right| \leq C |x|^{-k-1}, \quad 0\leq k\leq 2m, \quad x\in\RR^n\setminus\{0\}, \quad\rho\in (0,\infty),
\end{equation}

\noindent while for $n$ even 
\begin{equation}\label{eqo7.8.2}
\left|g(\log|x|^{-1},\log \rho^{-1})\right| \leq C_1+C_2 (C_R+\log |x|^{-1}), \quad 0<|x|, \rho<2R,
\end{equation}

\noindent and 
\begin{equation}\label{eqo7.8.3}
\left|\nabla_x^k\, g(\log|x|^{-1},\log \rho^{-1})\right| \leq C |x|^{-k}, \quad 1\leq k\leq 2m, \quad 0<|x|, \rho<2R.
\end{equation}

\noindent Here, as usually, we assume $|x|\neq \rho$ when $k=2m$, and lower derivatives of $g$ as well as $g$ itself are defined at $x$ such that $|x|= \rho$ by continuity. Hence, in particular, 
when $n$ is odd, $N\gg 1$,
\begin{equation}\label{eqo7.8.4}
\left|\nabla_x^k\, [|x|\,g(\log|x|^{-1},\log \rho^{-1})]\right| \lesssim |x|^{-k}\approx 2^{jk}, \quad 0\leq k\leq 2m, \quad x\in C_{2^{-j-2},2^{-j+4}}, \quad\rho<R,
\end{equation}

\noindent and when $n$ is even, 
\begin{equation}\label{eqo7.8.5}
\left|\nabla_x^k\, g(\log|x|^{-1},\log \rho^{-1})\right| \lesssim (C_R+\log |x|^{-1})|x|^{-k}\approx  j\, 2^{jk}, \quad 0\leq k\leq 2m, \quad x\in C_{2^{-j-2},2^{-j+4}}, \quad\rho<R, \end{equation}

\noindent since for $x\in C_{2^{-j-2},2^{-j+4}}$ we have $$C_R+\log |x|^{-1}=\log(4R)+\log |x|^{-1}=\log \frac{40}{|x|}\approx Cj.$$

Let us now concentrate on the case of the odd dimensions. Recall that $\Q_{g,\tau} (\vartheta)$ splits as in \eqref{eqn1}. Due to \eqref{eqo7.8.4} the second term gives 
\begin{multline}\label{eqn11}
\left|\sum_{\stackrel{i\geq 0, \,k\geq 0}{0<i+k\leq m}}\int_{\RR}\int_{S^{n-1}} \A_{ik}(\partial_t) [e^{-t}g(t,\tau)](t,\tau) \left(\partial_t^k\nabla_\omega^i\vartheta\right)^2 \,d\omega dt\right|\\[4pt]
\lesssim 
\sum_{j=N-2}^{\infty} \sum_{0\leq k\leq m}
\int_{C_{2^{-j},2^{-j+2}}}\frac{|\nabla^k
(U^{j}(x))|^2}{|x|^{n-2k}}\,dx\lesssim \sum_{j=N-2}^{\infty}
2^{-j(2m-n)}\int_{C_{2^{-j-2},2^{-j+4}}}|\nabla^m U^{j}(x)|^2\,dx,
\end{multline}
\noindent  using the Hardy's inequality for $U^j\in \ring W^{m,2}(C_{2^{-j-2},2^{-j+4}})$.
The last expression above is, in turn, bounded by 
\begin{multline}\label{eq8.32}
 C\sum_{j=N-2}^{\infty} 2^{-j(2m-n)}\,{\rm Cap}_P \{K\cap
{\overline{C_{2^{-j},2^{-j+2}}}},C_{2^{-j-2}, 2^{-j+4}}\} \\[4pt]\lesssim \sum_{j=N-2}^{\infty} 2^{-j(2m-n)}\,{\rm Cap}_P
\{{\overline{C_{2^{-j},2^{-j+2}}}}\setminus
\Omega,C_{2^{-j-2}, 2^{-j+4}}\} <C\eps,
\end{multline}

\noindent by (\ref{eq8.29}), (\ref{eq8.25})--(\ref{eq8.26}), the
monotonicity property (\ref{eq5.4}), and (\ref{eq8.3}). As for the second term in $\Q_{g,\tau} (\vartheta)$, we have
\begin{multline}\label{eq8.45}
\int_{S^{n-1}}\vartheta^2(\tau,\omega)\,d\omega=\frac{1}{|\xi|^{2\lambda}}\fint_{S_{|\xi|}}T^2(\xi)\,d\sigma_\xi\leq C
\sum_{j:\,2^{-j}\leq |\xi|\leq 2^{-j+2}}
\fint_{S_{|\xi|}}(U^j(\xi))^2\,d\sigma_\xi\\[4pt]
 \lesssim C  \sum_{j:\,2^{-j}\leq |\xi|\leq 2^{-j+2}} 2^{-2j\lambda}
\int_{\RR^n}(-\Delta)^m U^j(x)\,U^j(x)\,g(\log |x|^{-1}, \log |\xi|^{-1})\,dx\\[4pt]\lesssim \sum_{j=N-2}^\infty
2^{-j(2m-n)}\int_{C_{2^{-j-2},2^{-j+4}}}|\nabla^m U^{j}(x)|^2\,dx,
\end{multline}

\noindent using Theorem~\ref{to4.2} for the function $U^j\in \ring W^{m,2}(C_{2^{-j-2},2^{-j+4}})$. Finally, the
right-hand side of (\ref{eq8.45}) is bounded by $C\eps$ following \eqref{eqn11}--\eqref{eq8.32}.

This finishes the proof of \eqref{eqn7} when the underlying dimension is odd. 

In the case of the even dimension, we can follow the same argument with the only difference that the use of \eqref{eqo7.8.5} in place of \eqref{eqo7.8.4} introduces an extra logarithmic factor in the second expression in \eqref{eqn11} and, respectively, an extra factor $j$ in the series on the right-hand side of \eqref{eqn11}, in \eqref{eq8.32}, and on the right hand side of \eqref{eq8.45}. For the latter, we are using Theorem~\ref{to6.3}, and throughout $R=10$. We remark that the conclusion is still the same, that is, \eqref{eqn7} holds, since for the even dimensions we are using \eqref{eq8.3.even} which takes care of an extra factor $j$ in the series. 
\ep

Next, we build and estimate yet another auxiliary function. 
\begin{lemma}\label{lSt4.w} Given $P\in \Pi_1$ and $\PP(x):=|x|^\lambda P(x)$, $x\in\rn$, assume that the corresponding capacity satisfies \eqref{eq8.3}, \eqref{eq8.3.even}, in the case of the odd and the even dimension, respectively, and retain the definition of $T$, $\vartheta$  as above.  Furthermore, let $W_M$ be the solution of the problem 
\begin{equation}\label{def.w}
(-\Delta)^m W_M=(-\Delta)^m T \,\,{\mbox{in}}\,\,D_M, \quad W_M\in \ring W^{m,2}(D_M), \quad D_M=D\setminus B_{2^{-M}},
\end{equation}
\noindent with $M\gg N$.
Finally, let $w_M:=e^{\lambda t}(W_M\circ \varkappa^{-1})$. 
Fix some point $\xi\in B_{2^{-N}}$ and $\tau:=\log
|\xi|^{-1}$. Then 
\begin{equation}\label{eqn7.w}
\Q_{g,\tau} (w_M)<C\eps,
\end{equation}
is valid with a constant $C$ independent of $M$. In the case when the dimension is even, we assume that the parameter $R=10$ in the definition of $\Q_{g,\tau}$ (see \eqref{eqn1.1}).
\end{lemma}

We remark that $T$ belongs to $\ring W^{m,2}(B_{10})$. This can be seen following the same argument as in the Lemma above, since the $\ring W^{m,2}(B_{10})$ norm of $T$ is strictly smaller than the series appearing in the second expression in \eqref{eqn11} and than its analogue in the case of even dimensions. (One gains an extra power of $|x|$ for odd dimensions and a logarithmic factor in even dimensions compared to the situation in \eqref{eqn11}). Thus, $(-\Delta)^m T \in W^{-m,2}(D)$, the dual to $\ring W^{m,2}(D)$, and its resteriction to $D_M$ belongs to $W^{-m,2}(D_M)$. Hence,  
the boundary value problem \eqref{def.w} has a unique solution in $\ring W^{m,2}(D_M)$.

\vskip 0.08 in
\noindent {\it Proof of Lemma~\ref{lSt4.w}.} Throughout the proof of the Lemma we denote $W_M$ by $W$ and $w_M$ by $w$, respectively. 

Since $(-\Delta)^m W=-(-\Delta)^m T$ in $D$,  by \eqref{eqn2} we have
the formula
\begin{eqnarray}\label{eq8.33}\nonumber
\Q_{g,\tau}(w)&=&\int_{\RR^n}(-\Delta)^m
W(x)\, W(x) \,g(\log |x|^{-1}, \log |\xi|^{-1})\,dx\\[4pt]
& =& -\int_{\RR^n} (-\Delta)^m T(x)\, W(x)\,g(\log |x|^{-1}, \log |\xi|^{-1})\,dx.
\end{eqnarray}

\noindent
 In what follows we will show that
\begin{equation}\label{eq8.34}
-\int_{\RR^n} (-\Delta)^m T(x)\, W(x)\,g(\log |x|^{-1}, \log |\xi|^{-1})\,dx\leq C\eps^{1/2}(\Q_{g,\tau}(w))^{1/2}.
\end{equation}

\vskip 0.08 in 
\noindent {\bf Case I: odd dimension}. 

Once again, we start with the case when the dimension is odd. Let us pass to spherical coordinates and write 
\begin{multline}\label{eqn13}
-\int_{\RR^n} (-\Delta)^m T(x)\, W(x)\,g(\log |x|^{-1}, \log |\xi|^{-1})\,dx\\[4pt]
= 
-  \int_{\RR}\int_{S^{n-1}} {\mathcal
L}^{m,n}(\partial_t,\delta)\vartheta (t,\omega)\,w(t,\omega)[e^{-t}g(t,\tau)]\,d\omega dt.
\end{multline}
There are many ways to write the result as a bilinear form of $\vartheta$ and $w$. Here, we will follow step-by-step the procedure in the proof of Theorem~\ref{to4.2} in \cite{MayMazInvent2}, although considering bilinear rather than quadratic form leads to new technical difficulties, and the goal is different. 
%We will admit the convention that we ``integrate by parts" (in quotation marks, since \eqref{eqn13} is already understood in the sense of weak derivatives) from left to right, that is, taking the derivatives off $\vartheta$ and passing them on to $w(t,\omega)[e^{-t}g(t,\tau)]$ until the corresponding term has at most $m$ derivatives on $\vartheta$. 
The key point is that the terms where at least one derivative in $t$ lands on $w$ are harmless, as they conveniently fall under the scope of the second sum on the right-hand side of \eqref{eqn2.2}, for $w$ in place of $v$, and can be fairly directly estimated by $\Q_{g,\tau}(w)$ (we will provide the details below). The terms which do not possess any $t$ - derivatives of $w$ are considerably more delicate and have to be split in parts and treated carefully. 

Let us concentrate on the integral with  $h(t-\tau)$ in place of $[e^{-t}g(t,\tau)]$ first:
\begin{eqnarray}\label{eqo4.15}
&& \int_{\RR}\int_{S^{n-1}} {\mathcal
L}^{m,n}(\partial_t,\delta)\vt(t,\omega)\,w(t,\omega)h(t-\tau)\,d\omega
dt\nonumber\\[4pt]
&&\qquad = \int_{\RR}\int_{S^{n-1}} {\mathcal
L}^{m,n}(0,\delta)\vt(t,\omega)\,w(t,\omega)h(t-\tau)\,d\omega
dt\nonumber\\[4pt]
&&\qquad\quad + \int_{\RR}\int_{S^{n-1}} \left({\mathcal
L}^{m,n}(\partial_t,\delta)-{\mathcal
L}^{m,n}(0,\delta)\right)\vt(t,\omega)\,w(t,\omega)h(t-\tau)\,d\omega
dt =:J_1+J_2.
\end{eqnarray}

\noindent The term $J_1$ is left alone for the moment.  All terms of the operator $\left({\mathcal
L}^{m,n}(\partial_t,\delta)-{\mathcal
L}^{m,n}(0,\delta)\right)$  contain $\partial_t^k$ for some $k\geq 1$, and therefore, we can write

\begin{equation}\label{eqo4.16}
{\mathcal L}^{m,n}(\partial_t,\delta)-{\mathcal
L}^{m,n}(0,\delta)=\sum_{\stackrel{k\geq 1,\,i\geq 0}{2i+k\leq
2m}}d_{ik}(-\delta)^i\partial_t^k,\qquad \mbox{for some}\quad d_{ik}\in\RR.
\end{equation}

\noindent  Hence,
\begin{equation}\label{eqo4.17}
J_2 = \sum_{\stackrel{k\geq 1,\,i\geq 0}{2i+k\leq
2m}}d_{ik}\int_{\RR}\int_{S^{n-1}}
\partial_t^k\nabla_\omega^i \vt(t,\omega)\,\nabla_\omega^i
w(t,\omega)h(t-\tau)\,d\omega dt.
\end{equation}

We claim that
\begin{equation}\label{eqo4.18}
\int_{\RR}\partial_t^k\vt\,wh\,dt= 
\int_{\RR}\vt w\,(-\partial_t)^kh\,dt+\sum_{\stackrel{i\geq
0,\,j\geq 0,\,l\geq 1}{i+l+j\leq
k}}b_{ijl}^k\int_{\RR}\partial_t^i\vt \,\partial_t^l w\,\partial_t^jh\,dt,\quad
b_{ijl}^k\in\RR,
\end{equation}

\noindent for any $k\geq 1$. This will be proved by induction. For
$k=1$ we have

\begin{equation}\nonumber
 \int_{\RR}\partial_t\vt\,wh\,dt=-
\int_{\RR}\vt \partial_t w\,h\,dt-\int_{\RR} \vt w\,\partial_t h\,dt,
\end{equation}

\noindent which is the desired form.
Let us now assume that (\ref{eqo4.18}) holds for
$k=1,2,...,r-1$, and prove it for $k=r$. Indeed, using (\ref{eqo4.18}) for $r-1$,
\begin{multline}\label{eqo4.19}
\int_{\RR}\partial_t^r\vt\,w\,h\,dt=-
\int_{\RR}\partial_t^{r-1}\vt\,\partial_t w\,h \,dt -
\int_{\RR}\partial_t^{r-1}\vt \,w\,\partial_th\,dt\\[4pt]
=-\int_{\RR}\partial_t^{r-1}\vt\,\partial_t w\,h \,dt
 - \int_{\RR}\vt \,w\,(-\partial_t)^{r-1}\partial_t h\,dt
+\sum_{\stackrel{i\geq
0,\,j\geq 0,\,l\geq 1}{i+j+l\leq r-1}}b_{ijl}^{r-1}\int_{\RR}\partial_t^i\vt \partial_t^lw\,\partial_t^{j+1}h\,dt,
\end{multline}
\noindent which can be written in the form
(\ref{eqo4.18}) for $k=r$. This finishes the proof of \eqref{eqo4.18}.

Then, using \eqref{eqo4.18} for $v$ and $\nabla_\omega^i v$, \eqref{eqo4.17} leads to the representation
\begin{eqnarray}\label{eqo4.21}\nonumber
&& J_2 = \sum_{\stackrel{k\geq 1,\,i\geq 0}{2i+k\leq
2m}}d_{ik}\int_{\RR}\int_{S^{n-1}}
\nabla_\omega^i
\vt(t,\omega) \nabla_\omega^i w(t, \omega)\,(-\partial_t)^k h(t-\tau)\,d\omega dt\nonumber \\[4pt]
&&\qquad + \sum_{\stackrel{k\geq 1,\,i\geq 0}{2i+k\leq
2m}}\sum_{\stackrel{i'\geq 0,\,j'\geq 0,\,l'\geq 1}{i'+j'+l'\leq
k}}b_{i'j'l'}^{ki}\int_{\RR}\int_{S^{n-1}}
\partial_t^{i'}\nabla_\omega^i\vt \,\partial_t^{l'}\nabla_\omega^i w\,\partial_t^{j'}h(t-\tau)\,d\omega dt=:J_2'+J_2''.
\end{eqnarray}

Due to the condition $l'\geq 1$, the term $J_2''$, after possibly some more integration by parts taking the derivatives off $\vt$, is bounded as
\begin{multline}\label{eqn14}
|J_2''|\lesssim \sum_{\stackrel{i\geq 0, \,k\geq 0}{i+k\leq
m}} \left( \int_{\RR}\int_{S^{n-1}}
|\partial_t^{k}\nabla_\omega^i\vt|^2 d\omega dt\right)^{1/2} \sum_{\stackrel{i\geq 0, \,k\geq 1}{i+k\leq
m}} \left( \int_{\RR}\int_{S^{n-1}}
|\partial_t^{k}\nabla_\omega^i w|^2 d\omega dt\right)^{1/2}\\[4pt] \lesssim \eps^{1/2}\, \Q_{g,\tau} (w)^{1/2},
\end{multline} 
\noindent similarly to \eqref{eqn11}--\eqref{eq8.32}.

Turning to $J_2'$, 
\begin{eqnarray}\label{eqo4.22}\nonumber
&& J_2'= \sum_{\stackrel{k\geq 1,\,i\geq 0}{2i+k\leq
2m}}d_{ik}\sum_{p=0}^\infty \sum_{l=-p}^{p} \int_{\RR}  p^i(p+n-2)^i
\vt_{pl}(t) w_{pl}(t)(-\partial_t)^k h(t-\tau)\,  dt
\nonumber \\[4pt]
&&\qquad = \sum_{p=0}^\infty \sum_{l=-p}^{p} \int_{\RR} 
\vt_{pl}(t) w_{pl}(t) \sum_{\stackrel{k\geq 1,\,i\geq 0}{2i+k\leq
2m}}d_{ik} p^i(p+n-2)^i (-\partial_t)^k h(t-\tau)\,  dt
\nonumber \\[4pt]
&&\qquad =\sum_{p=0}^\infty \sum_{l=-p}^{p} \int_{\RR} 
\vt_{pl}(t) w_{pl}(t) \left({\mathcal L}^{m,n}(-\partial_t,-p(p+n-2))-{\mathcal
L}^{m,n}(0,-p(p+n-2))\right) h(t-\tau)\,  dt, \end{eqnarray}

\noindent where we employed \eqref{eqo4.16}. 

Now, 
\begin{eqnarray}\label{eqo4.23} 
&&-\sum_{p=0}^\infty \sum_{l=-p}^{p} \int_{\RR} 
\vt_{pl}(t) w_{pl}(t) {\mathcal
L}^{m,n}(0,-p(p+n-2))h(t-\tau)\,  dt\nonumber \\[4pt]
&& = -\int_{\RR}\int_{S^{n-1}}  {\mathcal
L}^{m,n}(0,\delta)
\vt(t,\omega)w(t,\omega)h(t-\tau)\,d\omega dt=-J_1,
\end{eqnarray}
\noindent hence, this term will be cancelled.  Also, by \eqref{eqo4.1}, we see that 
\begin{eqnarray}\label{eqo4.24}\sum_{p=0}^\infty \sum_{l=-p}^{p} \int_{\RR} 
\vt(t,\omega)w(t,\omega) {\mathcal L}^{m,n}(-\partial_t,0) h(t-\tau)\,d\omega dt= \int_{S^{n-1}} 
\vt(\tau, \omega)w(\tau, \omega) \,d\omega.
 \end{eqnarray}

\noindent Hence, 
\begin{multline}\label{eqo4.25}\nonumber
J_2' = \sum_{p=0}^\infty \sum_{l=-p}^{p} \int_{\RR}  
\vt_{pl}(t) w_{pl}(t)\left({\mathcal L}^{m,n}(-\partial_t,-p(p+n-2))-{\mathcal L}^{m,n}(-\partial_t,0)\right)  h(t-\tau)\,  dt \\[4pt]
-J_1+ \int_{S^{n-1}} 
\vt(\tau, \omega)w(\tau, \omega) \,d\omega\\[4pt]
=: J_{2,1}'-J_1+\int_{S^{n-1}} 
\vt(\tau, \omega)w(\tau, \omega) \,d\omega.
\end{multline}
Directly from \eqref{eqn2.1}, 
\begin{equation}\label{eqn15}
\int_{S^{n-1}} 
\vt(\tau, \omega)w(\tau, \omega) \,d\omega \lesssim \Q_{g,\tau} (\vt)^{1/2}\, \Q_{g,\tau} (w)^{1/2},
\end{equation} 
and we are left with $J_{2,1}'$.

Recall \eqref{eqo4.26}. It follows that for $0\leq p\leq m-\frac n2+\frac 12$ the integrals in $J_2'$ have positive weights and hence, 
\begin{multline}\label{eqn16}
\sum_{p=0}^{m-\frac n2+\frac 12} \sum_{l=-p}^{p} \int_{\RR}  
\vt_{pl}(t) w_{pl}(t)\left({\mathcal L}^{m,n}(-\partial_t,-p(p+n-2))-{\mathcal L}^{m,n}(-\partial_t,0)\right)  h(t-\tau)\,  dt \\[4pt]
\lesssim \left(\sum_{p=0}^{m-\frac n2+\frac 12} \sum_{l=-p}^{p} \int_{\RR}  
\vt_{pl}(t)^2\left({\mathcal L}^{m,n}(-\partial_t,-p(p+n-2))-{\mathcal L}^{m,n}(-\partial_t,0)\right)  h(t-\tau)\,  dt \right)^{1/2}\\[4pt]
\times \left(\sum_{p=0}^{m-\frac n2+\frac 12} \sum_{l=-p}^{p} \int_{\RR}  
w_{pl}(t)^2\left({\mathcal L}^{m,n}(-\partial_t,-p(p+n-2))-{\mathcal L}^{m,n}(-\partial_t,0)\right)  h(t-\tau)\,  dt \right)^{1/2}\\[4pt]
\lesssim \Q_{g,\tau} (\vt)^{1/2}\, \Q_{g,\tau} (w)^{1/2},
\end{multline}
by \eqref{eqn2.2}. On the other hand, 
\begin{multline}\label{eqn17}
\sum_{p=m-\frac n2+\frac 32}^\infty \sum_{l=-p}^{p} \int_{\RR}  
\vt_{pl}(t) w_{pl}(t)\left({\mathcal L}^{m,n}(-\partial_t,-p(p+n-2))-{\mathcal L}^{m,n}(-\partial_t,0)\right)  h(t-\tau)\,  dt \\[4pt]
\lesssim \left(\sum_{p=m-\frac n2+\frac 32}^\infty \sum_{l=-p}^{p} \int_{\RR}  
\vt_{pl}(t)^2\Bigg|\left({\mathcal L}^{m,n}(-\partial_t,-p(p+n-2))-{\mathcal L}^{m,n}(-\partial_t,0)\right)  h(t-\tau)\Bigg|\,  dt \right)^{1/2}\\[4pt]
\times \left(\sum_{p=m-\frac n2+\frac 32}^\infty \sum_{l=-p}^{p} \int_{\RR}  
w_{pl}(t)^2\Bigg|\left({\mathcal L}^{m,n}(-\partial_t,-p(p+n-2))-{\mathcal L}^{m,n}(-\partial_t,0)\right)  h(t-\tau)\Bigg|\,  dt \right)^{1/2},
\end{multline}
\noindent while 
\begin{multline}\label{eqo4.33}  \sum_{p=m-\frac n2+\frac 32}^\infty \sum_{l=-p}^{p} \int_{\RR} 
\vt_{pl}^2(t) \left|\left({\mathcal L}^{m,n}(-\partial_t,-p(p+n-2))-{\mathcal L}^{m,n}(-\partial_t,0)\right) h(t-\tau)\right|\,  dt\\[4pt]
\lesssim \sum_{i=0}^m \sum_{p=m-\frac n2+\frac 32}^{\infty} \sum_{l=-p}^{p} \Bigl(p\,(p+n-2)\Bigr)^i
\int_{\RR}\vt_{pl}^2(t) \,dt 
\\[4pt]\lesssim  
\int_{\RR}\int_{S^{n-1}} \vt\,\prod_{p=-\frac n2+\frac 32}^{m-\frac
n2+\frac 12}\left(-\delta-p\,(p+n-2)\right) \vt\,d\omega dt \lesssim \Q_{g,\tau}(\vt),
\end{multline}
\noindent (see \eqref{eqn3.1.1}), and an analogous estimate holds for $w$.  Hence, 
\begin{equation}\label{eqn18}
J_{2,1}' \lesssim \Q_{g,\tau} (\vt)^{1/2}\, \Q_{g,\tau} (w)^{1/2},
\end{equation} 

All in all, \eqref{eqo4.15}--\eqref{eqn18} gives 
\begin{multline}\label{eqn19}
 \int_{\RR}\int_{S^{n-1}} {\mathcal
L}^{m,n}(\partial_t,\delta)\vt(t,\omega)\,w(t,\omega)h(t-\tau)\,d\omega
dt \\[4pt]
\lesssim \Q_{g,\tau} (\vt)^{1/2}\, \Q_{g,\tau} (w)^{1/2}+ \eps^{1/2}\, \Q_{g,\tau} (w)^{1/2}.
\end{multline} 

A considerably simpler argument handles the case of a positive constant in place of $h$ in \eqref{eqo4.15}--\eqref{eqn18}: after writing an analogue of \eqref{eqo4.15} for a constant weight, we see that an analogue of $J_1$ is bounded by the term in the last line of \eqref{eqo4.33} to the power $1/2$ times the same term for $w$ (by \eqref{eqn3.1.1}), and hence, satisfies the bound by $ \Q_{g,\tau} (\vt)^{1/2}\, \Q_{g,\tau} (w)^{1/2}$, while the analogue of $J_2$ (which again has $\partial_t^k$, $k\geq 1$) lets us pass the $t$-derivative directly to $w$ and satisfies the same bound as \eqref{eqn14}. 

Collecting all these estimates, we arrive at 
\begin{multline}\label{eqn20}
-  \int_{\RR}\int_{S^{n-1}} {\mathcal
L}^{m,n}(\partial_t,\delta)\vartheta (t,\omega)\,w(t,\omega)[e^{-t}g(t,\tau)]\,d\omega dt  \\[4pt] \lesssim \Q_{g,\tau} (\vt)^{1/2}\, \Q_{g,\tau} (w)^{1/2}+ \eps^{1/2}\, \Q_{g,\tau} (w)^{1/2}, 
\end{multline} 
and therefore, due to \eqref{eq8.33}, \eqref{eqn10},
\begin{equation}\label{eqn21}
 \Q_{g,\tau} (w)^{1/2}  \lesssim \Q_{g,\tau} (\vt)^{1/2}+\eps^{1/2}\lesssim \eps^{1/2}. 
\end{equation} 
This yields, again by combination with \eqref{eqn7}, the validity of \eqref{eqn7.w}, and finishes the proof of Lemma~\ref{lSt4.w} in the case when the dimension is odd. 

\vskip 0.08 in \noindent {\bf Case II: even dimension}. 
We start again by writing 
\begin{multline}\label{eqn13-even}
-\int_{\RR^n} (-\Delta)^m T(x)\, W(x)\,g(\log |x|^{-1}, \log |\xi|^{-1})\,dx\\[4pt]
= 
-  \int_{\RR}\int_{S^{n-1}} {\mathcal
L}_o^{m,n}(\partial_t,\delta)\vartheta (t,\omega)\,w(t,\omega) g(t,\tau)\,d\omega dt, 
\end{multline}
and, recalling the definition of the function $\widetilde g$ from \eqref{eqo6.16}, 
\begin{eqnarray}\label{eqo4.15-even}
&& \int_{\RR}\int_{S^{n-1}} {\mathcal
L}_0^{m,n}(\partial_t,\delta)\vt(t,\omega)\,w(t,\omega)\widetilde g(t,\tau)\,d\omega
dt\nonumber\\[4pt]
&&\qquad = \int_{\RR}\int_{S^{n-1}} {\mathcal
L}^{m,n}_0(0,\delta)\vt(t,\omega)\,w(t,\omega)\widetilde g(t,\tau)\,d\omega
dt\nonumber\\[4pt]
&&\qquad\quad + \int_{\RR}\int_{S^{n-1}} \left({\mathcal
L}_0^{m,n}(\partial_t,\delta)-{\mathcal
L}^{m,n}(0,\delta)\right)\vt(t,\omega)\,w(t,\omega)\widetilde g(t,\tau)\,d\omega
dt =:J_1+J_2.
\end{eqnarray}
\noindent It is important to observe that 
\begin{equation}\label{eqo6.19-even} 
|\widetilde g(t,\tau)|\leq C_0(m,n)+|\mu^{(4)}|\,(C_R+t), \quad t,\tau\geq \log (2R)^{-1}, \end{equation}

\noindent and 
\begin{equation}\label{eqo6.20-even} 
|\partial_t^l\,\widetilde g(t,\tau)|\leq C_1(m,n),\quad t,\tau\geq \log (2R)^{-1}, \quad 1\leq l\leq 2m,
\end{equation}

\noindent for some constants $C_0(m,n), C_1(m,n)>0$ depending on $m,n$ only. We note that $\partial_t^l\,\widetilde g$ can be defined at $t=\tau$ for all $l<2m$ by continuity, and for $l=2m$ one assumes $t\neq \tau$ in \eqref{eqo6.20}.

Analogously to \eqref{eqo4.16}--\eqref{eqn14} and using \eqref{eqo6.19}--\eqref{eqo6.20}, the term $J_2$ splits into $J_2=J_2'+J_2''$, where 
\begin{multline}\label{eqn14-even}
|J_2''|\lesssim \sum_{\stackrel{i\geq 0, \,k\geq 0}{i+k\leq
m}} \left( \int_{\RR}\int_{S^{n-1}}
|\partial_t^{k}\nabla_\omega^i\vt|^2 (C_1+C_2(C_R+t))d\omega dt\right)^{1/2} \sum_{\stackrel{i\geq 0, \,k\geq 1}{i+k\leq
m}} \left( \int_{\RR}\int_{S^{n-1}}
|\partial_t^{k}\nabla_\omega^i w|^2 (C_1+C_2(C_R+t)) d\omega dt\right)^{1/2}\\[4pt] \lesssim \eps^{1/2}\, \Q_{g,\tau} (w)^{1/2}.
\end{multline} 
Here we used analogue of \eqref{eqn11}--\eqref{eq8.32} for even dimensions (which, in particular, accommodates the logarithmic weight) to handle the term with $\theta$ and \eqref{eqn2.1.0} to handle the term with $w$ above. 

As for $J_2'$, we have, similarly to the case when the dimension is odd (cf. \eqref{eqo4.22}), 
\begin{multline}\label{eqo4.22-even}
J_2' =\sum_{p=0}^\infty \sum_{l=-p}^{p} \int_{\RR} 
\vt_{pl}(t) w_{pl}(t) \left({\mathcal L}_0^{m,n}(-\partial_t,-p(p+n-2))-{\mathcal
L}^{m,n}_0(0,-p(p+n-2))\right) \widetilde g(t,\tau)\,  dt\\[4pt]
=\sum_{p=0}^\infty \sum_{l=-p}^{p} \int_{\RR} 
\vt_{pl}(t) w_{pl}(t) {\mathcal L}_0^{m,n}(-\partial_t,-p(p+n-2)) \widetilde g(t,\tau)\,  dt-J_1, \end{multline}
and $J_1$ cancels out (see \eqref{eqo4.15}). Turning to the remaining portion of $J_2'$, we split
\begin{multline}\label{eqo4.22.1-e}
\sum_{p=0}^\infty \sum_{l=-p}^{p} \int_{\RR} 
\vt_{pl}(t) w_{pl}(t) {\mathcal L}_0^{m,n}(-\partial_t,-p(p+n-2))\widetilde g(t,\tau)\,  dt\\[4pt]
= \sum_{p=0}^\infty \sum_{l=-p}^{p} \int_{\RR} 
\vt_{pl}(t) w_{pl}(t) \left({\mathcal L}_0^{m,n}(-\partial_t,-p(p+n-2))-{\mathcal L}_0^{m,n}(-\partial_t,-p_0(p_0+n-2))\right)\widetilde g(t,\tau)\,  dt\\[4pt]+\sum_{p=0}^\infty \sum_{l=-p}^{p} \int_{\RR} 
\vt_{pl}(t) w_{pl}(t)  {\mathcal L}_0^{m,n}(-\partial_t,-p_0(p_0+n-2))\widetilde g(t,\tau)\,  dt, \end{multline}
where
\begin{equation}\label{eqo6.23}
p_0:=\left\{\begin{array}{l}
0,\qquad\,\, \mbox{when } m-n/2 \mbox{ is even},\\[8pt]
1,\qquad\,\, \mbox{when } m-n/2 \mbox{ is odd}.\\[6pt]
\end{array}
\right.
\end{equation}
Since
\begin{equation}\label{eqo6.25} {\mathcal
L}^{m,n}_o(-\partial_t,-p_0(p_0+n-2)) \widetilde g(t,\tau)={\mathcal
L}^{m,n}_o(-\partial_t,-p_0(p_0+n-2))h(t-\tau)=\delta(t-\tau).
\end{equation}
(see, (6.25) in \cite{MayMazInvent2}, {\it loc.cit.}, for details), the second term on the right-hand side of \eqref{eqo4.22.1-e} reduces to the left-hand side of \eqref{eqn15} and enjoys the same estimates, using an analogue of \eqref{eq8.45} for even dimensions and \eqref{eqn2.1.0}. We are left with 
\begin{multline}\label{eqo4.22.2-e}
\sum_{p=0}^\infty \sum_{l=-p}^{p} \int_{\RR} 
\vt_{pl}(t) w_{pl}(t) \left({\mathcal L}_0^{m,n}(-\partial_t,-p(p+n-2))-{\mathcal L}_0^{m,n}(-\partial_t,-p_0(p_0+n-2))\right)\widetilde g(t,\tau)\,  dt\\[4pt]
= \sum_{p=0}^{m-n/2} \sum_{l=-p}^{p} \int_{\RR} 
\vt_{pl}(t) w_{pl}(t) \left({\mathcal L}_0^{m,n}(-\partial_t,-p(p+n-2))-{\mathcal L}_0^{m,n}(-\partial_t,-p_0(p_0+n-2))\right)\widetilde g(t,\tau)\,  dt\\[4pt]
+ \sum_{p=m-n/2+1}^\infty \sum_{l=-p}^{p} \int_{\RR} 
\vt_{pl}(t) w_{pl}(t) \left({\mathcal L}_0^{m,n}(-\partial_t,-p(p+n-2))-{\mathcal L}_0^{m,n}(-\partial_t,-p_0(p_0+n-2))\right)\widetilde g(t,\tau)\,  dt. \end{multline}
Much as in the case of odd dimensions, 
\begin{equation}\label{eqo6.41}
{\mathcal L}^{m,n}_o(0,-p(p+n-2))\geq Cp^m(p+n-2)^m, \quad\mbox{for $p>m-n/2$,}
\end{equation}
(cf. (6.42) in \cite{MayMazInvent2}, {\it loc.cit.}) and hence,
\begin{multline}\label{eqo4.22.3-e}
\sum_{p=m-n/2+1}^\infty \sum_{l=-p}^{p} \int_{\RR} 
\vt_{pl}(t) w_{pl}(t) \left({\mathcal L}_0^{m,n}(-\partial_t,-p(p+n-2))-{\mathcal L}_0^{m,n}(-\partial_t,-p_0(p_0+n-2))\right)\widetilde g(t,\tau)\,  dt\\[4pt]
\lesssim \left(\sum_{p=m-n/2+1}^\infty \sum_{l=-p}^{p} \int_{\RR} 
\vt_{pl}(t)^2  p^m(p+n-2)^m (C_1+C_2(C_R+t)) \,  dt\right)^{1/2}
\\[4pt]
\times\left(\sum_{p=m-n/2+1}^\infty \sum_{l=-p}^{p} \int_{\RR} 
w_{pl}(t)^2  p^m(p+n-2)^m (C_1+C_2(C_R+t)) \,  dt\right)^{1/2}\\[4pt]
\lesssim \left(\sum_{p=m-n/2+1}^\infty \sum_{l=-p}^{p} \int_{\RR} 
{\mathcal L}_0^{m,n}(0,-p(p+n-2)) \vt_{pl}(t) \,   \vt_{pl}(t) (C_1+C_2(C_R+t)) \,  dt\right)^{1/2}
\\[4pt]
\times\left(\sum_{p=m-n/2+1}^\infty \sum_{l=-p}^{p} \int_{\RR} 
{\mathcal L}_0^{m,n}(0,-p(p+n-2)) w_{pl}(t) \,   w_{pl}(t) (C_1+C_2(C_R+t)) \,  dt\right)^{1/2} \\[4pt]\lesssim \eps^{1/2}\, \Q_{g,\tau} (w)^{1/2},\end{multline}
using by now habitual considerations for the term with $\theta$ and \eqref{eqn2.1.0} for the term with $w$. This takes care of the large $p$'s in \eqref{eqo4.22.2-e}. 

Turning to $0\leq p\leq m-n/2$, we split the discussion according to the case when $m-n/2$ is even and the case when $m-n/2$ is odd. 

For $m-n/2$ even and $0\leq p\leq m-n/2$ even we use \eqref{eqo6.28-g} and proceed as in \eqref{eqn16} using \eqref{eqn2.1.0-refined} at the last step. The same considerations apply to the case when $m-n/2$ is odd and $0\leq p\leq m-n/2$ is odd. In the complementary scenario, we recall that 
\begin{multline}\label{eqo6.35}
{\mathcal L}^{m,n}_o(0,-p(p+n-2))\geq C, \quad\mbox{when $m-n/2$ is even and $0\leq p\leq m-n/2$ is odd,}\\[4pt]
\quad\mbox{or when $m-n/2$ is odd and $0\leq p\leq m-n/2$ is even,}
\end{multline}
(see (6.35) and (6.38) in \cite{MayMazInvent2}, {\it loc.cit.}) and hence, for all such $p$
\begin{multline}\label{eqo6.36}
\left|\Bigg({\mathcal L}^{m,n}_o(-\partial_t,-p(p+n-2))
-{\mathcal L}^{m,n}_o(-\partial_t,0) \Bigg) \,\,\widetilde g(t,\tau)\right|\leq C_1+C_2(C_R+t)\\[4pt] \leq {\mathcal L}^{m,n}_o(0,-p(p+n-2))\left(C_3+C_4(C_R+t)\right),\quad \quad t,\tau\geq \log (2R)^{-1}.
\end{multline}
Now the argument is finished as in \eqref{eqo4.22.3-e}.

Finally, we observe that one can treat $C'+C'' (C_R+t)$ in place of $\widetilde g$ via the same argument. In fact, the situation  is even simpler, as some tricky terms requiring positivity get annihilated. All in all, we arrive again at \eqref{eqn21}, as desired.
\ep

\begin{corollary}\label{corSt4.w} Given $P\in \Pi_1$ and $\PP(x):=|x|^\lambda P(x)$, $x\in\rn$, assume that the corresponding capacity satisfies \eqref{eq8.3}, \eqref{eq8.3.even}, in the case of the odd and the even dimension, respectively, and retain the definition of $T$, $\vartheta$  as above.  Furthermore, let $W$ be the solution of the problem 
\begin{equation}\label{def.w2}
(-\Delta)^m W=(-\Delta)^m T \,\,{\mbox{in}}\,\,D, \quad W\in \ring W^{m,2}(D).\end{equation}
\noindent Fix some point $\xi\in B_{2^{-N+1}}$ and $\tau:=\log
|\xi|^{-1}$. Then 
\begin{equation}\label{eqn7.w2}
\fiint_{C_{|\xi|, 2|\xi|}} \left(\frac{|W(x)|}{|x|^\lambda}\right)^2 \,dx
<C\eps,
\end{equation}
is valid with a constant $C$ depending on $m$ and $n$ only. 
\end{corollary}

\bp First of all, according to Lemma~\ref{lSt4.w}, we have 
\begin{equation}\label{eqn7.w2.1}
\fiint_{C_{|\xi|, 2|\xi|}} \left(\frac{|W_M(x)|}{|x|^\lambda}\right)^2 
<C\eps,
\end{equation}
with a constant $C$ independent of $M$.

On the other hand, the same argument as in \eqref{eqLim1.2}--\eqref{eqLim1.6} shows that $W_M$ converges to $W$ as $M\to\infty$ in $\ring W^{m,2}(D)$ norm. Since $W_M$ and $W$ are both in $\ring W^{m,2}(D)$, we have 
\begin{multline}\label{eqn7.w2.2}
\fiint_{C_{|\xi|, 2|\xi|}} \left(\frac{|W_M(x)-W(x)|}{|x|^\lambda}\right)^2 \,dx\approx \frac{1}{|\xi|^{2\lambda}}\fiint_{C_{|\xi|, 2|\xi|}} |W_M(x)-W(x)|^2 \,dx \\[4pt] \leq C_{\xi, D}  \,\|W_M-W\|_{\ring W^{m,2}(D)},
\end{multline}
\noindent by Poincar\'e inequality (see, e.g., \cite{MazSobSpNew}, (5.4.1)).
Due to the aforementioned convergence, there exists $M$, depending on $\eps$ and $\xi$, such that the right-hand side of \eqref{eqn7.w2.2} is smaller than $\eps$. Since the constant in \eqref{eqn7.w2.1} is independent of $M$ and $\xi$, this yields \eqref{eqn7.w2}, as desired.
\ep

\subsection{Conclusion of the proof}

At this stage we are ready to construct function $V$ and to finish the proof of the Theorem. Recall the definition of the cut-off $\eta$ from \eqref{eq8.4.0}. We claim that 
\begin{equation}\label{def.v} V:= T-W-\eta \PP
\end{equation} 
is exactly the solution to \eqref{eq8.4} on $D$. Let us be somewhat more precise. The function $V$ is defined by \eqref{def.v} on all of $B_{10}$ and is equal to zero in the complement of $D$ in the $\ring W^{m,2}$ sense. First of all, let us show that $V \in \ring W^{m,2}(D)$. 
To this end, observe that $W \in \ring W^{m,2}(D)$ by definition and hence, we can concentrate on $T-\eta \PP$. 

The fact that the $\ring W^{m,2}(D)$ norm of $V$ is finite follows directly from the fact that the norm estimates are valid for $T$ and $\eta \PP$ individually. For $T$ we have discussed this fact before the proof of Lemma~\ref{lSt4.w}. To analyze $\PP=|x|^\lambda\,P$, for some $P\in \Pi_1$, we separate the case of odd and even dimensions. When the dimension is odd, one simply observes that $|\nabla^k \PP(x)|$ is bounded by $|x|^{\lambda-k}$, hence, it has a finite  $W^{m,2}$ norm in any bounded set (possibly containing $O$). When the dimension is even, it is important to recall that spherical harmonics $Y^p_l$ are restrictions to the unit sphere of homogeneous polynomials of degree $p$, that is, $|x|^p Y^p_l$ is a homogeneous polynomial of degree $p$. Since the set $\Pi_1$ is made of linear combinations of spherical harmonics of degrees $m-n/2-2j$, for all $j=0,1,...$ such that $m-n/2-2j\geq 0$, the function $|x|^{m-n/2}P(x)$ is still a polynomial of degree $m-n/2$. Thus,  $|\nabla^k \PP(x)|$, for all $k=0,...,m-n/2$ are again polynomials (and higher derivatives are zero). Thus, once again, $\PP$ has a finite  $W^{m,2}$ norm in any bounded set (possibly containing $O$).

Furthermore, the boundary of $D$ consists of two parts,
$$\partial D= \partial B_{10}\cap \partial (K\cap \overline{B_{2^{-N}}}).$$
The latter portion has a part of $\partial K \cap \overline{B_{2^{-N}}}$ and a part of $K \cap \partial (\overline{B_{2^{-N}}})$, but as we shall see soon, this separation is irrelevant: what is important is that it is a portion of compactum $K$ lying in $\overline{B_{2^{-N}}}$. Indeed, due to the support properties of $\eta^j$ and the fact that they form a partition of unity we have 
$$ \Big(\sum_{j=1}^{\infty} \eta^j\Big)\Big|_{\overline{B_{2^{-N}}}}=1.$$
On the other hand, by definition, 
$$U_j\Big|_{K\cap \overline{C_{2^{-j}, 2^{-j+2}}}}=P$$
and $\supp (\eta^j)\subset C_{2^{-j}, 2^{-j+2}}$. Hence, 
$$\eta^j U_j\Big|_K=\eta_j P,$$
and 
$$ T\Big|_{K\cap \overline{B_{2^{-N}}}}=\PP.$$
All these equalities are taken in the $\ring W^{m,2}$ sense. Hence, since $\eta \PP=\PP$ on $K\cap \overline{B_{2^{-N}}}$, we have $V=0$ on $K\cap \overline{B_{2^{-N}}}$ in the $\ring W^{m,2}$ sense, as desired. Since $\supp \eta \subset B_2$ and $T$ is supported in $B_{2}$ as well, we have $V\in \ring W^{m,2}(D)$ as desired. 

The fact that $(-\Delta)^mV=f$ in $D$, with $f$ given by \eqref{eq8.4.0}, is now a consequence of definitions and polyharmonicity of $\PP$ in $\RR^n\setminus \{O\}.$

According to Lemma~\ref{lSt4} and Corollary~\ref{corSt4.w}, we have, for every $\xi \in B_{2^{-N-1}},$
\begin{equation}\label{eq.small} 
\fiint_{C_{|\xi|, 2|\xi|}} \left(\frac{|W(x)|}{|x|^\lambda}\right)^2 \,dx+ \fiint_{C_{|\xi|, 2|\xi|}} \left(\frac{|T(x)|}{|x|^\lambda}\right)^2 \,dx
\leq C\eps,
\end{equation}
with the constant $C$ depending on $m$ and $n$ only (and independent of $M$, $\Omega$, and $\xi$). Since 
$$ \fiint_{C_{|\xi|, 2|\xi|}} \left(\frac{|\PP(x)|}{|x|^\lambda}\right)^2 \,dx =  \fiint_{C_{|\xi|, 2|\xi|}} |P(x)|^2 \,dx  \geq C,$$ 
where $C$ is a positive constant independent of $\xi$, it follows that  
$$ \fiint_{C_{|\xi|, 2|\xi|}} \left(\frac{|V(x)|}{|x|^\lambda}\right)^2 \,dx \geq C,$$ 
where $C$ is a positive constant independent of $\xi$.
This implies that $\nabla^\lambda V$ does not vanish at $O$, as desired, and finishes the proof of the Theorem. \ep
% For self: the last line above is just using pointwise the mean value theorem


\begin{thebibliography}{999}

\bibitem{AdamsPotAnal} D.R.\,Adams, {\it $L\sp p$ potential theory
techniques and nonlinear PDE}, Potential theory (Nagoya, 1990),
1--15, de Gruyter, Berlin, 1992.

\bibitem{AWiener} D.\,Adams, {\it Potential and capacity before and
after Wiener}, Proceedings of the Norbert Wiener Centenary
Congress, 1994 (East Lansing, MI, 1994), 63--83, Proc. Sympos.
Appl. Math., 52, Amer. Math. Soc., Providence, RI, 1997.

\bibitem{AHe} D.\,Adams, L.\,Hedberg, {\it
Function spaces and potential theory}, Grundlehren der
Mathematischen Wissenschaften [Fundamental Principles of
Mathematical Sciences], 314. Springer-Verlag, Berlin, 1996.

\bibitem{Agmon} S.\,Agmon, {\it Maximum theorems for solutions of higher order
elliptic equations}, Bull. Amer. Math. Soc. 66 1960 77--80.

\bibitem{ADN} S.\,Agmon, A.\,Douglis and L.\,Nirenberg, {\it Estimates
near the boundary for solutions of elliptic partial differential
equations satisfying general boundary conditions. I}, Comm. Pure
Appl. Math., 12 (1959), 623--727.

\bibitem{Antman2005} S.S.\,Antman, {\it Nonlinear problems of elasticity.} Second edition. Applied Mathematical Sciences, 107. Springer, New York, 2005.

\bibitem{AliceChang} S.-Y. A.\,Chang, {\it
Conformal invariants and partial differential equations},
Bull. Amer. Math. Soc. (N.S.) 42 (2005), no. 3, 365--393 (electronic). 

\bibitem{AliceChang2} S.-Y. A.\,Chang,
P.C.\,Yang, {\it Non-linear partial differential equations in conformal geometry,}  Proceedings of the International Congress of Mathematicians, Vol. I (Beijing, 2002),  189--207, Higher Ed. Press, Beijing, 2002.

\bibitem{CiarletElasticity123} Ph.\,Ciarlet, {\it Mathematical elasticity.} Vol. II. Theory of plates. Studies in Mathematics and its Applications, 27. North-Holland Publishing Co., Amsterdam, 1997. 

\bibitem{DJK} B.\,Dahlberg, C.\,Kenig, G.\,Verchota, 
{\it The Dirichlet problem for the biharmonic equation in a Lipschitz domain.} (French summary) 
Ann. Inst. Fourier (Grenoble) 36 (1986), no. 3, 109Ð135. 

\bibitem{DMM} G.\,Dal Maso,
U.\,Mosco, {\it Wiener criteria and energy decay for relaxed
Dirichlet problems}, Arch. Rational Mech. Anal. 95 (1986), no. 4,
345--387.

\bibitem{EG}L.C.\,Evans, R.F.\,Gariepy, {\it Wiener's criterion
for the heat equation}, Arch. Rational Mech. Anal. 78 (1982), no.
4, 293--314.

\bibitem{FGL} E.\,Fabes, N.\,Garofalo, E.\,Lanconelli, {\it
Wiener's criterion for divergence form parabolic operators with
$C^1$-Dini continuous coefficients}, Duke Math. J. 59 (1989),
no. 1, 191--232.

\bibitem{FJK} E.\,Fabes, D.\,Jerison, C.\,Kenig, {\it The Wiener test
for degenerate elliptic equations}, Ann. Inst. Fourier (Grenoble)
32 (1982), no. 3, vi, 151--182.

\bibitem{KMR} V.\,Kozlov, V.\,Maz'ya, J.\,Rossmann,
{\it Spectral problems associated with corner singularities of
solutions to elliptic equations}, Mathematical Surveys and
Monographs, 85. American Mathematical Society, Providence, RI,
2001.

\bibitem{La} D.\,Labutin, {\it Potential estimates for a class of
fully nonlinear elliptic equations}, Duke Math. J. 111 (2002), no.
1, 1--49.

\bibitem{L} H.\,Lebesgue, {\it Sur le cas d'impossibilit\'e du probl\`eme de Dirichlet ordinaire}, C.R. des S\'eances de la Soci\'et\'e Math\'ematique de France, 17 (1913). 

\bibitem{LSW} W.\,Littman, G.\,Stampacchia,
H.F.\,Weinberger, {\it Regular points for elliptic equations with
discontinuous coefficients}, Ann. Scuola Norm. Sup. Pisa (3), 17
(1963), 43--77.

\bibitem{MZ}  J.\,Mal\'y, W.P.\,Ziemer, {\it Fine regularity of
solutions of elliptic partial differential equations},
Mathematical Surveys and Monographs, 51. American Mathematical
Society, Providence, RI, 1997.

\bibitem{MayMaz2} S.\,Mayboroda, V.\,Maz'ya, {\it Boundedness of the gradient of a solution and
Wiener test of order one for the biharmonic equation}, Invent. Math. 175 (2009), no. 2, 287--334.

\bibitem{MayMazGr} S.\,Mayboroda, V.\,Maz'ya, {\it Pointwise estimates for the polyharmonic Green function in general domains}, Cialdea, Alberto (ed.) et al., Analysis, partial differential equations and applications. The Vladimir Maz'ya anniversary volume. Selected papers of the international workshop, Rome, Italy, June 30--July 3, 2008. Basel: BirkhŠuser. Operator Theory: Advances and Applications 193, 143-158 (2009).

\bibitem{MayMazInvent2} S.\,Mayboroda, V.\,Maz'ya, {\it Regularity of solutions to the polyharmonic equation in general domains}, Invent. Math. 196 (2014), no. 1, 1Ð68.

\bibitem{MazyaSobolevSpaces} V.\,Maz'ya, {\it Prostranstva S. L. Soboleva}. (Russian) [Sobolev spaces] Leningrad. Univ., Leningrad, 1985. 416 pp.

\bibitem{MazSobSpNew}  V.\,Maz'ya, {\it Sobolev spaces with applications to elliptic partial differential equations. Second, revised and augmented edition.} Grundlehren der Mathematischen Wissenschaften [Fundamental Principles of Mathematical Sciences], 342. Springer, Heidelberg, 2011.

\bibitem{MWiener} V.\,Maz'ya, {\it Unsolved problems connected with
the Wiener criterion}, The Legacy of Norbert Wiener: A Centennial
Symposium (Cambridge, MA, 1994), 199--208, Proc. Sympos. Pure
Math., 60, Amer. Math. Soc., Providence, RI, 1997.

\bibitem{M2} V.\,Maz'ya, {\it
The Wiener test for higher order elliptic equations.}, Duke Math.
J. 115 (2002), no. 3, 479--512.

\bibitem{MazNP83}
V.\,Maz'ya, S.A.\,Nazarov, and B.A.\,Plamenevski{\u\i}, \emph{Singularities
  of solutions of the {D}irichlet problem in the exterior of a thin cone}, Mat.
  Sb. (N.S.) 122 (164) (1983), no.~4, 435--457. English translation:
  Math. USSR-Sb. 50 (1985), no.~2, 415--437.


\bibitem{MR1} V.\,Maz'ya, J.\,Rossmann, {\it  On the
Agmon-Miranda maximum principle for solutions of elliptic
equations in polyhedral and polygonal domains}, Ann. Global Anal.
Geom. 9 (1991), no. 3, 253--303.

\bibitem{MR} V.\,Maz'ya, J.\,Rossmann, {\it On the Agmon-Miranda maximum principle for solutions of strongly elliptic equations in domains of $\RR^n$ with conical points}, Ann. Global Anal. Geom. 10
(1992), no. 2, 125--150.

\bibitem{Meleshko} V.V.\,Meleshko, {\it Selected topics in the history of the
two-dimensional biharmonic problem},  Appl. Mech. Rev. 56 (2003), 
33--85.

\bibitem{Mir48}
C.\,Miranda, \emph{Formule di maggiorazione e teorema di esistenza per le
  funzioni biarmoniche de due variabili}, Giorn. Mat. Battaglini (4)
  2(78) (1948), 97--118. 

\bibitem{Mir58}
C.\,Miranda, \emph{Teorema del massimo modulo e teorema di esistenza e di unicit\`a
  per il problema di {D}irichlet relativo alle equazioni ellittiche in due
  variabili}, Ann. Mat. Pura Appl. (4) 46 (1958), 265--311.


\bibitem{Necas}  J.\,Ne\v cas, {\it Les m\'ethodes directes en th\'eorie
des \'equations elliptiques}, Masson et Cie, \'Editeurs, Paris;
Academia, \'Editeurs, Prague 1967.

\bibitem{PVLp} J.\,Pipher, G.\,Verchota, {\it The Dirichlet problem in
$L\sp p$ for the biharmonic equation on Lipschitz domains}, Amer.
J. Math. 114 (1992), no. 5, 923--972.

\bibitem{PVmax}  J.\,Pipher, G.\,Verchota, {\it A maximum principle
for biharmonic functions in Lipschitz and $C^1$ domains},
Comment. Math. Helv. 68 (1993), no. 3, 385--414.

\bibitem{PVpoly}  J.\,Pipher, G.\,Verchota, {\it Maximum principles for the polyharmonic
equation on Lipschitz domains}, Potential Anal. 4 (1995), no. 6,
615--636.

\bibitem{P} H.\,Poincar\'e, {\it Sur les \'equations aux deriv\'ees partielle de la physique math\'ematique,}, Amer. J. Math., 12 (1890), 211--299.

\bibitem{Kamenev} T.\,Sedrakyan, L.\,Glazman, and A.\,Kamenev, {\it Absence of Bose condensation in certain frustrated lattices}, arXiv:1303.7272 (2014)

\bibitem{Shen1} Z.\,Shen, {\it On estimates of biharmonic functions on Lipschitz and convex domains,} J. Geom. Anal. 16 (2006), no. 4, 721Ð734.

\bibitem{Shen2} Z.\,Shen, {\it  Necessary and sufficient conditions for the solvability of the Lp Dirichlet problem on Lipschitz domains}, Math. Ann. 336 (2006), no. 3, 697Ð725.

\bibitem{Skrypnik} I.V.\,Skrypnik, {\it Methods for analysis of nonlinear elliptic boundary value problems.} Translated from the 1990 Russian original by Dan D. Pascali. Translations of Mathematical Monographs, 139. American Mathematical Society, Providence, RI, 1994.

\bibitem{TW} N.\,Trudinger, X.-J.\,Wang,
{\it On the weak continuity of elliptic operators and applications
to potential theory}, (English summary) Amer. J. Math. 124 (2002),
no. 2, 369--410.

\bibitem{Wiener} N.\,Wiener,  {\it The Dirichlet problem},  J. Math. Phys.  3  (1924),  pp.
127--146.

\bibitem{Z} S.C.\,Zaremba, {\it Sur le principe du minimum}, Bull. Acad. Sci. Cracovie, (1909).

\end{thebibliography}
\end{document}